\documentclass{article}

\PassOptionsToPackage{numbers, compress}{natbib}
\usepackage[preprint]{neurips_2026} % arXiv v2 (extended version): no workshop notice

\usepackage[utf8]{inputenc} % allow utf-8 input
\usepackage[T1]{fontenc}    % use 8-bit T1 fonts
\usepackage{hyperref}       % hyperlinks
\usepackage{url}            % simple URL typesetting
\usepackage{booktabs}       % professional-quality tables
\usepackage{amsfonts}       % blackboard math symbols
\usepackage{amsmath}
\usepackage{amssymb}
\usepackage{nicefrac}       % compact symbols for 1/2, etc.
\usepackage{microtype}      % microtypography
\usepackage{xcolor}         % colors
\usepackage{colortbl}       % shaded table cells (appendix per-mode table)
\usepackage{graphicx}
\usepackage[section]{placeins} % arXiv version: flush floats at section boundaries

\title{A Response Theory Probe for Learned Stochastic AI Simulators, Tested on Lorenz-63}

\author{%
  João Böger\\
  Department of Technology,\\ 
  Management and Economics\\
  Technical University of Denmark\\
  Kongens Lyngby, 2800, DK \\
  \texttt{jpade@dtu.dk} \\
  \And
  Simon Driscoll\\
  Department of Applied\\
  Mathematics and \\
  Theoretical Physics\\
  University of Cambridge\\
  Cambridge, CB3 0WA, UK \\
  \texttt{sd2136@cam.ac.uk} \\
  \And
  Niccolò Zagli\\
  School of Computing and \\Mathematical Sciences\\
  University of Leicester\\
  Leicester, LE1 7RH, UK \\
  \texttt{nz68@leicester.ac.uk} \\
  \AND
  Valerio Lucarini\\
  School of Computing and \\Mathematical Sciences\\
  University of Leicester\\
  Leicester, LE1 7RH, UK \\
  School of Sciences \\
  Great Bay University \\
  Dongguan, P.R. China \\
  \texttt{v.lucarini@leicester.ac.uk} \\
  \And
  Francisco Camara Pereira\\
  Department of Technology,\\ 
  Management and Economics\\
  Technical University of Denmark\\
  Kongens Lyngby, 2800, DK \\
  \texttt{camara@dtu.dk} \\
}

\newcommand{\GzeroSigchishift}{0.0549} % results/camera_ready/camera_ready_metrics.json:claims.GzeroSigchishift
\newcommand{\Gzerochishift}{0.0443} % results/camera_ready/camera_ready_metrics.json:claims.Gzerochishift
\newcommand{\Gzeron}{10} % results/camera_ready/camera_ready_metrics.json:claims.Gzeron
\newcommand{\Gzeroslowflag}{0.10} % results/camera_ready/camera_ready_metrics.json:claims.Gzeroslowflag
\newcommand{\MLPdriftSecantMaxDevPct}{7.21} % results/camera_ready/camera_ready_metrics.json:claims.MLPdriftSecantMaxDevPct
\newcommand{\MLPdriftn}{80} % results/camera_ready/camera_ready_metrics.json:claims.MLPdriftn
\newcommand{\MLPdriftonewing}{33} % results/camera_ready/camera_ready_metrics.json:claims.MLPdriftonewing
\newcommand{\MLPflowacc}{0.33} % results/camera_ready/camera_ready_metrics.json:claims.MLPflowacc
\newcommand{\MLPflowfplus}{0.474} % results/camera_ready/camera_ready_metrics.json:claims.MLPflowfplus
\newcommand{\MLPflowslow}{-0.146} % results/camera_ready/camera_ready_metrics.json:claims.MLPflowslow
\newcommand{\NODEinvflag}{0.93} % results/camera_ready/camera_ready_metrics.json:claims.NODEinvflag
\newcommand{\NODEn}{80} % results/camera_ready/camera_ready_metrics.json:claims.NODEn
\newcommand{\NODEslowflag}{0.24} % results/camera_ready/camera_ready_metrics.json:claims.NODEslowflag
\newcommand{\NODEslowflagTruthThr}{0.07} % results/camera_ready/camera_ready_metrics.json:claims.NODEslowflagTruthThr
\newcommand{\NODEswitchAboveOne}{80} % results/camera_ready/camera_ready_metrics.json:claims.NODEswitchAboveOne
\newcommand{\NSDEgapXPct}{4.33} % results/camera_ready/camera_ready_metrics.json:claims.NSDEgapXPct
\newcommand{\NSDEinvflag}{1.00} % results/camera_ready/camera_ready_metrics.json:claims.NSDEinvflag
\newcommand{\NSDEn}{80} % results/camera_ready/camera_ready_metrics.json:claims.NSDEn
\newcommand{\NSDEsigmaRelDiffPctMax}{0.204} % results/camera_ready/camera_ready_metrics.json:claims.NSDEsigmaRelDiffPctMax
\newcommand{\NSDEslowflag}{0.23} % results/camera_ready/camera_ready_metrics.json:claims.NSDEslowflag
\newcommand{\NSDEswitchAboveOne}{80} % results/camera_ready/camera_ready_metrics.json:claims.NSDEswitchAboveOne
\newcommand{\RCexhibitFisherp}{4.9\times10^{-10}} % results/camera_ready/camera_ready_metrics.json:claims.RCexhibitFisherp
\newcommand{\RCfisherSwitchSlowp}{0.424} % results/camera_ready/camera_ready_metrics.json:claims.RCfisherSwitchSlowp
\newcommand{\RCfreqOffsetMaxPct}{0.133} % results/camera_ready/camera_ready_metrics.json:claims.RCfreqOffsetMaxPct
\newcommand{\RCn}{78} % results/camera_ready/camera_ready_metrics.json:claims.RCn
\newcommand{\RCnX}{26} % results/camera_ready/camera_ready_metrics.json:claims.RCnX
\newcommand{\RCnXslowfail}{20} % results/camera_ready/camera_ready_metrics.json:claims.RCnXslowfail
\newcommand{\RCrealTwoFasterFrac}{0.77} % results/camera_ready/camera_ready_metrics.json:claims.RCrealTwoFasterFrac
\newcommand{\RCslowfailInSwitch}{0.70} % results/camera_ready/camera_ready_metrics.json:claims.RCslowfailInSwitch
\newcommand{\RCslowfailOutSwitch}{0.80} % results/camera_ready/camera_ready_metrics.json:claims.RCslowfailOutSwitch
\newcommand{\RCslowflag}{0.76} % results/camera_ready/camera_ready_metrics.json:claims.RCslowflag
\newcommand{\RCslowflagTruthThr}{0.55} % results/camera_ready/camera_ready_metrics.json:claims.RCslowflagTruthThr
\newcommand{\SINDyhfplus}{0.416} % results/camera_ready/camera_ready_metrics.json:claims.SINDyhfplus
\newcommand{\SindyFulln}{10} % results/camera_ready/camera_ready_metrics.json:claims.SindyFulln
\newcommand{\TrueRealOneRecv}{0.0508} % results/camera_ready/camera_ready_metrics.json:claims.TrueRealOneRecv
\newcommand{\TrueRealOneRemed}{-0.389} % results/camera_ready/camera_ready_metrics.json:claims.TrueRealOneRemed
\newcommand{\Truen}{50} % results/camera_ready/camera_ready_metrics.json:claims.Truen
\newcommand{\TruenX}{48} % results/camera_ready/camera_ready_metrics.json:claims.TruenX
\newcommand{\TruenXslowfail}{3} % results/camera_ready/camera_ready_metrics.json:claims.TruenXslowfail
\newcommand{\accOneNODEmax}{1.54} % results/camera_ready/camera_ready_metrics.json:appendix.accumulation.node.1.0.max
\newcommand{\accOneNODEmed}{0.87} % results/camera_ready/camera_ready_metrics.json:appendix.accumulation.node.1.0.median
\newcommand{\accOneNODEmin}{0.62} % results/camera_ready/camera_ready_metrics.json:appendix.accumulation.node.1.0.min
\newcommand{\accOneNSDEmax}{1.51} % results/camera_ready/camera_ready_metrics.json:appendix.accumulation.nsde.1.0.max
\newcommand{\accOneNSDEmed}{0.87} % results/camera_ready/camera_ready_metrics.json:appendix.accumulation.nsde.1.0.median
\newcommand{\accOneNSDEmin}{0.62} % results/camera_ready/camera_ready_metrics.json:appendix.accumulation.nsde.1.0.min
\newcommand{\accOneRCmax}{1.65} % results/camera_ready/camera_ready_metrics.json:appendix.accumulation.rc_flowmap.1.0.max
\newcommand{\accOneRCmed}{0.87} % results/camera_ready/camera_ready_metrics.json:appendix.accumulation.rc_flowmap.1.0.median
\newcommand{\accOneRCmin}{0.65} % results/camera_ready/camera_ready_metrics.json:appendix.accumulation.rc_flowmap.1.0.min
\newcommand{\accOneRefmed}{0.65} % results/camera_ready/camera_ready_metrics.json:appendix.accumulation.reference.1.0.median
\newcommand{\accOneSindymax}{1.30} % results/camera_ready/camera_ready_metrics.json:appendix.accumulation.sindy.1.0.max
\newcommand{\accOneSindymed}{0.97} % results/camera_ready/camera_ready_metrics.json:appendix.accumulation.sindy.1.0.median
\newcommand{\accOneSindymin}{0.76} % results/camera_ready/camera_ready_metrics.json:appendix.accumulation.sindy.1.0.min
\newcommand{\accOneTruemax}{1.32} % results/camera_ready/camera_ready_metrics.json:appendix.accumulation.true.1.0.max
\newcommand{\accOneTruemed}{0.89} % results/camera_ready/camera_ready_metrics.json:appendix.accumulation.true.1.0.median
\newcommand{\accOneTruemin}{0.60} % results/camera_ready/camera_ready_metrics.json:appendix.accumulation.true.1.0.min
\newcommand{\bandEonebhi}{0.266} % results/camera_ready/camera_ready_metrics.json:criteria.bands.e1b.hi
\newcommand{\bandEoneblo}{0.163} % results/camera_ready/camera_ready_metrics.json:criteria.bands.e1b.lo
\newcommand{\bandEoneblooOut}{0.04} % results/camera_ready/camera_ready_metrics.json:criteria.bands.e1b.loo_inside (1 - x)
\newcommand{\bandNEighteenShi}{0.324} % results/camera_ready/camera_ready_metrics.json:criteria.bands.n18.hi
\newcommand{\bandNEighteenSlo}{0.18} % results/camera_ready/camera_ready_metrics.json:criteria.bands.n18.lo
\newcommand{\bandNEighteenSn}{18} % results/camera_ready/camera_ready_metrics.json:criteria.bands.n18.n
\newcommand{\bandNullkOut}{6} % results/camera_ready/camera_ready_metrics.json:criteria.band.null_test.k_out
\newcommand{\bandNullnTest}{50} % results/camera_ready/camera_ready_metrics.json:criteria.band.null_test.n_test
\newcommand{\bandNullp}{0.0144} % results/camera_ready/camera_ready_metrics.json:criteria.band.null_test.binom_p_upper
\newcommand{\bandPooledhi}{0.309} % results/camera_ready/camera_ready_metrics.json:criteria.bands.pooled.hi
\newcommand{\bandPooledlo}{0.162} % results/camera_ready/camera_ready_metrics.json:criteria.bands.pooled.lo
\newcommand{\bandPooledlooOut}{0.02} % results/camera_ready/camera_ready_metrics.json:criteria.bands.pooled.loo_inside (1 - x)
\newcommand{\bandPooledn}{100} % results/camera_ready/camera_ready_metrics.json:criteria.bands.pooled.n
\newcommand{\bandRunnerhi}{0.309} % results/camera_ready/camera_ready_metrics.json:criteria.bands.runner.hi
\newcommand{\bandRunnerlo}{0.162} % results/camera_ready/camera_ready_metrics.json:criteria.bands.runner.lo
\newcommand{\bandVarF}{1.46} % results/camera_ready/camera_ready_metrics.json:criteria.band.variance_test.F_runner_over_e1b
\newcommand{\bandVarp}{0.191} % results/camera_ready/camera_ready_metrics.json:criteria.band.variance_test.F_p
\newcommand{\chiShareLeadTruthcv}{1.09} % results/camera_ready/camera_ready_metrics.json:claims.chiShareLeadTruthcv
\newcommand{\chiShareTruthn}{10} % results/camera_ready/camera_ready_metrics.json:claims.chiShareTruthn
\newcommand{\chiZeroBandloo}{0.02} % results/camera_ready/camera_ready_metrics.json:criteria.band.loo_inside (1 - x: LOO out-rate)
\newcommand{\chiZeroBandn}{100} % results/camera_ready/camera_ready_metrics.json:bands.primary.n
\newcommand{\chiZeroTruemed}{0.221} % results/camera_ready/camera_ready_metrics.json:models.true.metrics.chi0.median
\newcommand{\chiaccTruth}{0.89} % results/camera_ready/camera_ready_metrics.json:claims.chiaccTruth
\newcommand{\chiflagmax}{0.09} % results/camera_ready/camera_ready_metrics.json:claims.chiflagmax
\newcommand{\coefActiveFull}{18} % results/camera_ready/camera_ready_metrics.json:appendix.inputs.sindy.sindy_full.active terms
\newcommand{\coefActiveHundred}{30} % results/camera_ready/camera_ready_metrics.json:appendix.inputs.sindy.sindy100.active terms
\newcommand{\coefActiveTrue}{7} % results/camera_ready/camera_ready_metrics.json:appendix.inputs.sindy.true.active terms
\newcommand{\coefBetaFull}{2.667} % results/camera_ready/camera_ready_metrics.json:appendix.inputs.sindy.sindy_full.β (−z in dz)
\newcommand{\coefBetaHundred}{2.402} % results/camera_ready/camera_ready_metrics.json:appendix.inputs.sindy.sindy100.β (−z in dz)
\newcommand{\coefBetaTrue}{2.667} % results/camera_ready/camera_ready_metrics.json:appendix.inputs.sindy.true.β (−z in dz)
\newcommand{\coefBigFull}{8} % results/camera_ready/camera_ready_metrics.json:appendix.inputs.sindy.sindy_full.terms with |coef| ≥ 0.1
\newcommand{\coefBigHundred}{14} % results/camera_ready/camera_ready_metrics.json:appendix.inputs.sindy.sindy100.terms with |coef| ≥ 0.1
\newcommand{\coefBigTrue}{7} % results/camera_ready/camera_ready_metrics.json:appendix.inputs.sindy.true.terms with |coef| ≥ 0.1
\newcommand{\coefConstDxFull}{0.000} % results/camera_ready/camera_ready_metrics.json:appendix.inputs.sindy.sindy_full.const dx
\newcommand{\coefConstDxHundred}{-3.423} % results/camera_ready/camera_ready_metrics.json:appendix.inputs.sindy.sindy100.const dx
\newcommand{\coefConstDxTrue}{0.000} % results/camera_ready/camera_ready_metrics.json:appendix.inputs.sindy.true.const dx
\newcommand{\coefConstDyFull}{0.307} % results/camera_ready/camera_ready_metrics.json:appendix.inputs.sindy.sindy_full.const dy
\newcommand{\coefConstDyHundred}{3.139} % results/camera_ready/camera_ready_metrics.json:appendix.inputs.sindy.sindy100.const dy
\newcommand{\coefConstDyTrue}{0.000} % results/camera_ready/camera_ready_metrics.json:appendix.inputs.sindy.true.const dy
\newcommand{\coefConstDzFull}{0.000} % results/camera_ready/camera_ready_metrics.json:appendix.inputs.sindy.sindy_full.const dz
\newcommand{\coefConstDzHundred}{-2.169} % results/camera_ready/camera_ready_metrics.json:appendix.inputs.sindy.sindy100.const dz
\newcommand{\coefConstDzTrue}{0.000} % results/camera_ready/camera_ready_metrics.json:appendix.inputs.sindy.true.const dz
\newcommand{\coefMaxRelErrPct}{0.116} % results/camera_ready/camera_ready_metrics.json:claims.coefMaxRelErrPct
\newcommand{\coefMinusXFull}{10.003} % results/camera_ready/camera_ready_metrics.json:appendix.inputs.sindy.sindy_full.−(x in dx)
\newcommand{\coefMinusXHundred}{10.755} % results/camera_ready/camera_ready_metrics.json:appendix.inputs.sindy.sindy100.−(x in dx)
\newcommand{\coefMinusXTrue}{10.000} % results/camera_ready/camera_ready_metrics.json:appendix.inputs.sindy.true.−(x in dx)
\newcommand{\coefRhoFull}{27.967} % results/camera_ready/camera_ready_metrics.json:appendix.inputs.sindy.sindy_full.ρ (x in dy)
\newcommand{\coefRhoHundred}{27.734} % results/camera_ready/camera_ready_metrics.json:appendix.inputs.sindy.sindy100.ρ (x in dy)
\newcommand{\coefRhoTrue}{28.000} % results/camera_ready/camera_ready_metrics.json:appendix.inputs.sindy.true.ρ (x in dy)
\newcommand{\coefSigmaFull}{10.001} % results/camera_ready/camera_ready_metrics.json:appendix.inputs.sindy.sindy_full.σ (y in dx)
\newcommand{\coefSigmaHundred}{10.520} % results/camera_ready/camera_ready_metrics.json:appendix.inputs.sindy.sindy100.σ (y in dx)
\newcommand{\coefSigmaTrue}{10.000} % results/camera_ready/camera_ready_metrics.json:appendix.inputs.sindy.true.σ (y in dx)
\newcommand{\eFoldRealOneRef}{2.73} % results/camera_ready/camera_ready_metrics.json:claims.eFoldRealOneRef
\newcommand{\eFoldRealTwoRef}{1.97} % results/camera_ready/camera_ready_metrics.json:claims.eFoldRealTwoRef
\newcommand{\eigRealTwoReRCmed}{-0.639} % results/camera_ready/camera_ready_metrics.json:models.rc_flowmap.final.groups.real2.re.median
\newcommand{\eigRealTwoReTruemax}{-0.483} % results/camera_ready/camera_ready_metrics.json:models.true.final.groups.real2.re.max
\newcommand{\eigRealTwoReTruemin}{-0.611} % results/camera_ready/camera_ready_metrics.json:models.true.final.groups.real2.re.min
\newcommand{\fitFPlusFull}{0.500} % results/camera_ready/camera_ready_metrics.json:appendix.inputs.sindy.f_plus_window.sindy_full
\newcommand{\fitFPlusHundred}{0.464} % results/camera_ready/camera_ready_metrics.json:appendix.inputs.sindy.f_plus_window.sindy100
\newcommand{\fitPairsFull}{10000000} % results/camera_ready/camera_ready_metrics.json:appendix.inputs.sindy.n_pairs.sindy_full
\newcommand{\fitPairsHundred}{100000} % results/camera_ready/camera_ready_metrics.json:appendix.inputs.sindy.n_pairs.sindy100
\newcommand{\fitTFull}{10000} % results/camera_ready/camera_ready_metrics.json:appendix.inputs.sindy.T_fit_tu.sindy_full
\newcommand{\fitTHundred}{100} % results/camera_ready/camera_ready_metrics.json:appendix.inputs.sindy.T_fit_tu.sindy100
\newcommand{\flagInvTransformer}{0.79} % results/camera_ready/camera_ready_metrics.json:exploratory.transformer.final.flags.rate_inv
\newcommand{\flagInvTrue}{0.04} % results/camera_ready/camera_ready_metrics.json:models.true.final.flags.rate_inv
\newcommand{\flagOscMergedTransformer}{0.15} % results/camera_ready/camera_ready_metrics.json:exploratory.transformer.appendix.osc_merged.rate
\newcommand{\flagOscTransformer}{0.73} % 58/80 = 0.725, rounded half-up to match \pctOscTransformer (Table 8); results/camera_ready/camera_ready_metrics.json:exploratory.transformer.final.flags.rate_osc
\newcommand{\flagOscTrue}{0.02} % results/camera_ready/camera_ready_metrics.json:models.true.final.flags.rate_osc
\newcommand{\flagSlowTransformer}{0.19} % results/camera_ready/camera_ready_metrics.json:exploratory.transformer.final.flags.rate_slow
\newcommand{\flagSlowTrue}{0.06} % results/camera_ready/camera_ready_metrics.json:models.true.final.flags.rate_slow
\newcommand{\floorAcfkOut}{0} % results/camera_ready/camera_ready_metrics.json:criteria.floors.acf_distance.null_test.k_out
\newcommand{\floorAcflooOut}{0.02} % results/camera_ready/camera_ready_metrics.json:criteria.floors.acf_distance.loo_inside (1 - x)
\newcommand{\floorAcfp}{1} % results/camera_ready/camera_ready_metrics.json:criteria.floors.acf_distance.null_test.binom_p_upper
\newcommand{\floorFPluskOut}{0} % results/camera_ready/camera_ready_metrics.json:criteria.floors.f_plus.null_test.k_out
\newcommand{\floorFPluslooOut}{0.04} % results/camera_ready/camera_ready_metrics.json:criteria.floors.f_plus.loo_inside (1 - x)
\newcommand{\floorFPlusp}{1} % results/camera_ready/camera_ready_metrics.json:criteria.floors.f_plus.null_test.binom_p_upper
\newcommand{\floorJointLooOut}{0.10} % results/camera_ready/camera_ready_metrics.json:criteria.floor_joint_loo_inside_e1b (1 - x)
\newcommand{\floorN}{50} % results/camera_ready/camera_ready_metrics.json:floor.n_pairs
\newcommand{\floorSwitchkOut}{2} % results/camera_ready/camera_ready_metrics.json:criteria.floors.switch_ratio.null_test.k_out
\newcommand{\floorSwitchlooOut}{0.04} % results/camera_ready/camera_ready_metrics.json:criteria.floors.switch_ratio.loo_inside (1 - x)
\newcommand{\floorSwitchp}{0.6} % results/camera_ready/camera_ready_metrics.json:criteria.floors.switch_ratio.null_test.binom_p_upper
\newcommand{\floorWOnekOut}{0} % results/camera_ready/camera_ready_metrics.json:criteria.floors.w1_x.null_test.k_out
\newcommand{\floorWOnelooOut}{0.02} % results/camera_ready/camera_ready_metrics.json:criteria.floors.w1_x.loo_inside (1 - x)
\newcommand{\floorWOnep}{1} % results/camera_ready/camera_ready_metrics.json:criteria.floors.w1_x.null_test.binom_p_upper
\newcommand{\gateBoundX}{31.7} % results/camera_ready/camera_ready_metrics.json:criteria.gate.bounds[0]
\newcommand{\gateBoundY}{43.6} % results/camera_ready/camera_ready_metrics.json:criteria.gate.bounds[1]
\newcommand{\gateBoundZ}{76.1} % results/camera_ready/camera_ready_metrics.json:criteria.gate.bounds[2]
\newcommand{\gateFactor}{1.5} % results/camera_ready/camera_ready_metrics.json:criteria.gate.factor
\newcommand{\gridN}{20000} % results/camera_ready/camera_ready_metrics.json:appendix.inputs.drift_rmse.grid.n
\newcommand{\gridSeed}{909} % results/camera_ready/camera_ready_metrics.json:appendix.inputs.drift_rmse.grid.seed
\newcommand{\inBandEonebGzeroH}{5} % results/camera_ready/camera_ready_metrics.json:models.true_rk4_h.final.bands_in.e1b
\newcommand{\inBandEonebGzeroHSig}{3} % results/camera_ready/camera_ready_metrics.json:models.true_rk4_h_sig.final.bands_in.e1b
\newcommand{\inBandEonebMLPdrift}{5} % results/camera_ready/camera_ready_metrics.json:models.mlp_drift.final.bands_in.e1b
\newcommand{\inBandEonebMLPflow}{0} % results/camera_ready/camera_ready_metrics.json:models.mlp_flow.final.bands_in.e1b
\newcommand{\inBandEonebNODE}{61} % results/camera_ready/camera_ready_metrics.json:models.node.final.bands_in.e1b
\newcommand{\inBandEonebNSDE}{62} % results/camera_ready/camera_ready_metrics.json:models.nsde.final.bands_in.e1b
\newcommand{\inBandEonebRC}{67} % results/camera_ready/camera_ready_metrics.json:models.rc_flowmap.final.bands_in.e1b
\newcommand{\inBandEonebSindy}{7} % results/camera_ready/camera_ready_metrics.json:models.sindy.final.bands_in.e1b
\newcommand{\inBandEonebSindyFull}{10} % results/camera_ready/camera_ready_metrics.json:diagnostic.sindy_full.final.bands_in.e1b
\newcommand{\inBandEonebTrue}{44} % results/camera_ready/camera_ready_metrics.json:models.true.final.bands_in.e1b
\newcommand{\inBandNEighteenSGzeroH}{10} % results/camera_ready/camera_ready_metrics.json:models.true_rk4_h.final.bands_in.n18
\newcommand{\inBandNEighteenSGzeroHSig}{8} % results/camera_ready/camera_ready_metrics.json:models.true_rk4_h_sig.final.bands_in.n18
\newcommand{\inBandNEighteenSMLPdrift}{8} % results/camera_ready/camera_ready_metrics.json:models.mlp_drift.final.bands_in.n18
\newcommand{\inBandNEighteenSMLPflow}{2} % results/camera_ready/camera_ready_metrics.json:models.mlp_flow.final.bands_in.n18
\newcommand{\inBandNEighteenSNODE}{66} % results/camera_ready/camera_ready_metrics.json:models.node.final.bands_in.n18
\newcommand{\inBandNEighteenSNSDE}{67} % results/camera_ready/camera_ready_metrics.json:models.nsde.final.bands_in.n18
\newcommand{\inBandNEighteenSRC}{59} % results/camera_ready/camera_ready_metrics.json:models.rc_flowmap.final.bands_in.n18
\newcommand{\inBandNEighteenSSindy}{5} % results/camera_ready/camera_ready_metrics.json:models.sindy.final.bands_in.n18
\newcommand{\inBandNEighteenSSindyFull}{9} % results/camera_ready/camera_ready_metrics.json:diagnostic.sindy_full.final.bands_in.n18
\newcommand{\inBandNEighteenSTrue}{44} % results/camera_ready/camera_ready_metrics.json:models.true.final.bands_in.n18
\newcommand{\inBandPooledGzeroH}{10} % results/camera_ready/camera_ready_metrics.json:models.true_rk4_h.final.bands_in.pooled
\newcommand{\inBandPooledGzeroHSig}{7} % results/camera_ready/camera_ready_metrics.json:models.true_rk4_h_sig.final.bands_in.pooled
\newcommand{\inBandPooledMLPdrift}{7} % results/camera_ready/camera_ready_metrics.json:models.mlp_drift.final.bands_in.pooled
\newcommand{\inBandPooledMLPflow}{1} % results/camera_ready/camera_ready_metrics.json:models.mlp_flow.final.bands_in.pooled
\newcommand{\inBandPooledNODE}{73} % results/camera_ready/camera_ready_metrics.json:models.node.final.bands_in.pooled
\newcommand{\inBandPooledNSDE}{73} % results/camera_ready/camera_ready_metrics.json:models.nsde.final.bands_in.pooled
\newcommand{\inBandPooledRC}{72} % results/camera_ready/camera_ready_metrics.json:models.rc_flowmap.final.bands_in.pooled
\newcommand{\inBandPooledSindy}{7} % results/camera_ready/camera_ready_metrics.json:models.sindy.final.bands_in.pooled
\newcommand{\inBandPooledSindyFull}{10} % results/camera_ready/camera_ready_metrics.json:diagnostic.sindy_full.final.bands_in.pooled
\newcommand{\inBandPooledTrue}{50} % results/camera_ready/camera_ready_metrics.json:models.true.final.bands_in.pooled
\newcommand{\inBandRunnerGzeroH}{10} % results/camera_ready/camera_ready_metrics.json:models.true_rk4_h.final.bands_in.runner
\newcommand{\inBandRunnerGzeroHSig}{7} % results/camera_ready/camera_ready_metrics.json:models.true_rk4_h_sig.final.bands_in.runner
\newcommand{\inBandRunnerMLPdrift}{7} % results/camera_ready/camera_ready_metrics.json:models.mlp_drift.final.bands_in.runner
\newcommand{\inBandRunnerMLPflow}{1} % results/camera_ready/camera_ready_metrics.json:models.mlp_flow.final.bands_in.runner
\newcommand{\inBandRunnerNODE}{73} % results/camera_ready/camera_ready_metrics.json:models.node.final.bands_in.runner
\newcommand{\inBandRunnerNSDE}{73} % results/camera_ready/camera_ready_metrics.json:models.nsde.final.bands_in.runner
\newcommand{\inBandRunnerRC}{72} % results/camera_ready/camera_ready_metrics.json:models.rc_flowmap.final.bands_in.runner
\newcommand{\inBandRunnerSindy}{7} % results/camera_ready/camera_ready_metrics.json:models.sindy.final.bands_in.runner
\newcommand{\inBandRunnerSindyFull}{10} % results/camera_ready/camera_ready_metrics.json:diagnostic.sindy_full.final.bands_in.runner
\newcommand{\inBandRunnerTrue}{50} % results/camera_ready/camera_ready_metrics.json:models.true.final.bands_in.runner
\newcommand{\inWOneRC}{76} % results/camera_ready/camera_ready_metrics.json:models.rc_flowmap.final.invariant.w1_x.k_inside
\newcommand{\mFreqAOne}{9.85} % results/camera_ready/camera_ready_metrics.json:models.mlp_drift.final.per_mode_full.merged[0].freqs[0]
\newcommand{\mFreqBOne}{9.98} % results/camera_ready/camera_ready_metrics.json:models.mlp_drift.final.per_mode_full.merged[0].freqs[1]
\newcommand{\maxGdevAllPct}{1.92} % results/camera_ready/camera_ready_metrics.json:claims.maxGdevAllPct
\newcommand{\maxGdevTruthPct}{0.83} % results/camera_ready/camera_ready_metrics.json:claims.maxGdevTruthPct
\newcommand{\mnullMeanOne}{0.9982} % results/camera_ready/camera_ready_metrics.json:models.mlp_drift.final.per_mode_full.merged[0].null_mean
\newcommand{\mnullSdOne}{0.0003} % results/camera_ready/camera_ready_metrics.json:models.mlp_drift.final.per_mode_full.merged[0].null_std
\newcommand{\mthrOne}{0.9900} % results/camera_ready/camera_ready_metrics.json:models.mlp_drift.final.per_mode_full.merged[0].threshold
\newcommand{\mthrUncOne}{0.9974} % results/camera_ready/camera_ready_metrics.json:models.mlp_drift.final.per_mode_full.merged[0].threshold_uncapped
\newcommand{\nDivGzeroH}{0} % results/camera_ready/camera_ready_metrics.json:models.true_rk4_h.final.flags.n_diverged
\newcommand{\nDivGzeroHSig}{0} % results/camera_ready/camera_ready_metrics.json:models.true_rk4_h_sig.final.flags.n_diverged
\newcommand{\nDivMLPdrift}{0} % results/camera_ready/camera_ready_metrics.json:models.mlp_drift.final.flags.n_diverged
\newcommand{\nDivMLPflow}{0} % results/camera_ready/camera_ready_metrics.json:models.mlp_flow.final.flags.n_diverged
\newcommand{\nDivNODE}{0} % results/camera_ready/camera_ready_metrics.json:models.node.final.flags.n_diverged
\newcommand{\nDivNSDE}{0} % results/camera_ready/camera_ready_metrics.json:models.nsde.final.flags.n_diverged
\newcommand{\nDivRC}{2} % results/camera_ready/camera_ready_metrics.json:models.rc_flowmap.final.flags.n_diverged
\newcommand{\nDivSindy}{0} % results/camera_ready/camera_ready_metrics.json:models.sindy.final.flags.n_diverged
\newcommand{\nDivSindyFull}{0} % results/camera_ready/camera_ready_metrics.json:diagnostic.sindy_full.final.flags.n_diverged
\newcommand{\nDivTransformer}{0} % results/camera_ready/camera_ready_metrics.json:exploratory.transformer.final.flags.n_diverged
\newcommand{\nDivTrue}{0} % results/camera_ready/camera_ready_metrics.json:models.true.final.flags.n_diverged
\newcommand{\nOneWingMLPflow}{0} % results/camera_ready/camera_ready_metrics.json:models.mlp_flow.final.flags.n_one_wing
\newcommand{\nRolloutSeeds}{10} % results/camera_ready/camera_ready_metrics.json:claims.nRolloutSeeds
\newcommand{\nRunsTransformer}{80} % results/camera_ready/camera_ready_metrics.json:exploratory.transformer.n_total
\newcommand{\nScoredGzeroH}{10} % results/camera_ready/camera_ready_metrics.json:models.true_rk4_h.n_scored
\newcommand{\nScoredGzeroHSig}{10} % results/camera_ready/camera_ready_metrics.json:models.true_rk4_h_sig.n_scored
\newcommand{\nScoredMLPdrift}{80} % results/camera_ready/camera_ready_metrics.json:models.mlp_drift.n_scored
\newcommand{\nScoredMLPflow}{80} % results/camera_ready/camera_ready_metrics.json:models.mlp_flow.n_scored
\newcommand{\nScoredNODE}{80} % results/camera_ready/camera_ready_metrics.json:models.node.n_scored
\newcommand{\nScoredNSDE}{80} % results/camera_ready/camera_ready_metrics.json:models.nsde.n_scored
\newcommand{\nScoredRC}{78} % results/camera_ready/camera_ready_metrics.json:models.rc_flowmap.n_scored
\newcommand{\nScoredSindy}{10} % results/camera_ready/camera_ready_metrics.json:models.sindy.n_scored
\newcommand{\nScoredSindyFull}{10} % results/camera_ready/camera_ready_metrics.json:diagnostic.sindy_full.n_scored
\newcommand{\nScoredTransformer}{80} % results/camera_ready/camera_ready_metrics.json:exploratory.transformer.n_scored
\newcommand{\nScoredTrue}{50} % results/camera_ready/camera_ready_metrics.json:models.true.n_scored
\newcommand{\nTrainSeeds}{8} % results/camera_ready/camera_ready_metrics.json:claims.nTrainSeeds
\newcommand{\nXTransformer}{17} % results/camera_ready/camera_ready_metrics.json:exploratory.transformer.final.joint.nX
\newcommand{\nXslowfailTransformer}{2} % results/camera_ready/camera_ready_metrics.json:exploratory.transformer.final.joint.nX_slowfail
\newcommand{\nodeReproMaxFracRange}{0} % results/camera_ready/camera_ready_metrics.json:claims.nodeReproMaxFracRange
\newcommand{\noiseH}{0.01} % results/camera_ready/camera_ready_metrics.json:appendix.inputs.nsde_noise.h
\newcommand{\noiseHatX}{0.9567} % results/camera_ready/camera_ready_metrics.json:appendix.inputs.nsde_noise.sigma_hat_over_s[0]
\newcommand{\noiseHatY}{0.9990} % results/camera_ready/camera_ready_metrics.json:appendix.inputs.nsde_noise.sigma_hat_over_s[1]
\newcommand{\noiseHatZ}{0.9891} % results/camera_ready/camera_ready_metrics.json:appendix.inputs.nsde_noise.sigma_hat_over_s[2]
\newcommand{\noiseN}{80} % results/camera_ready/camera_ready_metrics.json:appendix.inputs.nsde_noise.n
\newcommand{\noisePhiX}{0.9587} % results/camera_ready/camera_ready_metrics.json:appendix.inputs.nsde_noise.sigma_phi_over_s[0]
\newcommand{\noisePhiY}{0.9990} % results/camera_ready/camera_ready_metrics.json:appendix.inputs.nsde_noise.sigma_phi_over_s[1]
\newcommand{\noisePhiZ}{0.9883} % results/camera_ready/camera_ready_metrics.json:appendix.inputs.nsde_noise.sigma_phi_over_s[2]
\newcommand{\noiseRelaxX}{0.9520} % results/camera_ready/camera_ready_metrics.json:appendix.inputs.nsde_noise.relaxation_prediction[0]
\newcommand{\noiseRelaxY}{0.9950} % results/camera_ready/camera_ready_metrics.json:appendix.inputs.nsde_noise.relaxation_prediction[1]
\newcommand{\noiseRelaxZ}{0.9868} % results/camera_ready/camera_ready_metrics.json:appendix.inputs.nsde_noise.relaxation_prediction[2]
\newcommand{\noiseTrueDriftX}{0.9567} % results/camera_ready/camera_ready_metrics.json:appendix.inputs.nsde_noise.true_drift_sigma_hat_over_s[0]
\newcommand{\noiseTrueDriftY}{0.9982} % results/camera_ready/camera_ready_metrics.json:appendix.inputs.nsde_noise.true_drift_sigma_hat_over_s[1]
\newcommand{\noiseTrueDriftZ}{0.9885} % results/camera_ready/camera_ready_metrics.json:appendix.inputs.nsde_noise.true_drift_sigma_hat_over_s[2]
\newcommand{\paNODEOne}{80} % results/camera_ready/camera_ready_metrics.json:models.node.final.per_mode_full.modes[0].a
\newcommand{\paNSDEOne}{80} % results/camera_ready/camera_ready_metrics.json:models.nsde.final.per_mode_full.modes[0].a
\newcommand{\paRCOne}{77} % results/camera_ready/camera_ready_metrics.json:models.rc_flowmap.final.per_mode_full.modes[0].a
\newcommand{\paramAlpha}{0.05} % results/camera_ready/camera_ready_metrics.json:parameters.alpha
\newcommand{\paramBeta}{2.667} % results/camera_ready/camera_ready_metrics.json:parameters.beta
\newcommand{\paramBurnIn}{100} % results/camera_ready/camera_ready_metrics.json:parameters.burn_in_tu
\newcommand{\paramCalibK}{3} % results/camera_ready/camera_ready_metrics.json:parameters.calib_k
\newcommand{\paramDict}{729} % results/camera_ready/camera_ready_metrics.json:parameters.dictionary_size
\newcommand{\paramDtFitMLPdrift}{0.01} % results/camera_ready/camera_ready_metrics.json:parameters.dt_fit.mlp_drift
\newcommand{\paramDtFitMLPflow}{0.01} % results/camera_ready/camera_ready_metrics.json:parameters.dt_fit.mlp_flow
\newcommand{\paramDtFitNODE}{0.01} % results/camera_ready/camera_ready_metrics.json:parameters.dt_fit.node
\newcommand{\paramDtFitNSDE}{0.01} % results/camera_ready/camera_ready_metrics.json:parameters.dt_fit.nsde
\newcommand{\paramDtFitRC}{0.02} % results/camera_ready/camera_ready_metrics.json:parameters.dt_fit.rc_flowmap
\newcommand{\paramDtFitSindy}{0.001} % results/camera_ready/camera_ready_metrics.json:parameters.dt_fit.sindy
\newcommand{\paramDtFitSindyFull}{0.001} % results/camera_ready/camera_ready_metrics.json:parameters.dt_fit.sindy_full
\newcommand{\paramDtFitTransformer}{0.01} % results/camera_ready/camera_ready_metrics.json:exploratory.transformer.appendix.recipe.dt_fit
\newcommand{\paramDtRaw}{0.001} % results/camera_ready/camera_ready_metrics.json:parameters.dt_raw
\newcommand{\paramDtRollMLPdrift}{0.001} % results/camera_ready/camera_ready_metrics.json:parameters.dt_rollout.mlp_drift
\newcommand{\paramDtRollMLPflow}{0.01} % results/camera_ready/camera_ready_metrics.json:parameters.dt_rollout.mlp_flow
\newcommand{\paramDtRollNODE}{0.01} % results/camera_ready/camera_ready_metrics.json:parameters.dt_rollout.node
\newcommand{\paramDtRollNSDE}{0.01} % results/camera_ready/camera_ready_metrics.json:parameters.dt_rollout.nsde
\newcommand{\paramDtRollRC}{0.02} % results/camera_ready/camera_ready_metrics.json:parameters.dt_rollout.rc_flowmap
\newcommand{\paramDtRollSindy}{0.001} % results/camera_ready/camera_ready_metrics.json:parameters.dt_rollout.sindy
\newcommand{\paramDtRollSindyFull}{0.001} % results/camera_ready/camera_ready_metrics.json:parameters.dt_rollout.sindy_full
\newcommand{\paramDtRollTransformer}{0.01} % results/camera_ready/camera_ready_metrics.json:exploratory.transformer.appendix.recipe.dt_rollout
\newcommand{\paramFigQhi}{95} % results/camera_ready/camera_ready_metrics.json:parameters.fig_q_hi
\newcommand{\paramFigQlo}{5} % results/camera_ready/camera_ready_metrics.json:parameters.fig_q_lo
\newcommand{\paramGridN}{20000} % results/camera_ready/camera_ready_metrics.json:parameters.grid_n
\newcommand{\paramHorizonA}{1} % results/camera_ready/camera_ready_metrics.json:parameters.horizon_a
\newcommand{\paramHorizonB}{5} % results/camera_ready/camera_ready_metrics.json:parameters.horizon_b
\newcommand{\paramHorizonC}{20} % results/camera_ready/camera_ready_metrics.json:parameters.horizon_c
\newcommand{\paramKmax}{4} % results/camera_ready/camera_ready_metrics.json:parameters.k_max
\newcommand{\paramMlpBatch}{2048} % results/camera_ready/camera_ready_metrics.json:parameters.recipe.MlpBatch
\newcommand{\paramMlpDepth}{2} % results/camera_ready/camera_ready_metrics.json:parameters.recipe.MlpDepth
\newcommand{\paramMlpEpochs}{25} % results/camera_ready/camera_ready_metrics.json:parameters.recipe.MlpEpochs
\newcommand{\paramMlpInitSeed}{5} % results/camera_ready/camera_ready_metrics.json:parameters.recipe.MlpInitSeed
\newcommand{\paramMlpLr}{10^{-3}} % results/camera_ready/camera_ready_metrics.json:parameters.recipe.MlpLr
\newcommand{\paramMlpPatience}{5} % results/camera_ready/camera_ready_metrics.json:parameters.recipe.MlpPatience
\newcommand{\paramMlpValPct}{10} % results/camera_ready/camera_ready_metrics.json:parameters.recipe.MlpValPct
\newcommand{\paramMlpWidth}{64} % results/camera_ready/camera_ready_metrics.json:parameters.recipe.MlpWidth
\newcommand{\paramMseExp}{-3} % results/camera_ready/camera_ready_metrics.json:parameters.mse_exp
\newcommand{\paramNCalib}{16} % results/camera_ready/camera_ready_metrics.json:parameters.n_calib
\newcommand{\paramNLead}{6} % results/camera_ready/camera_ready_metrics.json:parameters.n_lead
\newcommand{\paramNodeLSel}{5} % results/camera_ready/camera_ready_metrics.json:parameters.recipe.NodeLSel
\newcommand{\paramNodeLrSel}{10^{-3}} % results/camera_ready/camera_ready_metrics.json:parameters.recipe.NodeLrSel
\newcommand{\paramNodeSweepLOne}{5} % results/camera_ready/camera_ready_metrics.json:parameters.recipe.NodeSweepLOne
\newcommand{\paramNodeSweepLTwo}{20} % results/camera_ready/camera_ready_metrics.json:parameters.recipe.NodeSweepLTwo
\newcommand{\paramNodeSweepLrOne}{10^{-3}} % results/camera_ready/camera_ready_metrics.json:parameters.recipe.NodeSweepLrOne
\newcommand{\paramNodeSweepLrTwo}{3\times10^{-4}} % results/camera_ready/camera_ready_metrics.json:parameters.recipe.NodeSweepLrTwo
\newcommand{\paramOne}{1} % results/camera_ready/camera_ready_metrics.json:parameters.one
\newcommand{\paramRcInput}{0.5} % results/camera_ready/camera_ready_metrics.json:parameters.recipe.RcInput
\newcommand{\paramRcLeak}{0.3} % results/camera_ready/camera_ready_metrics.json:parameters.recipe.RcLeak
\newcommand{\paramRcN}{1000} % results/camera_ready/camera_ready_metrics.json:parameters.recipe.RcN
\newcommand{\paramRcRadius}{1} % results/camera_ready/camera_ready_metrics.json:parameters.recipe.RcRadius
\newcommand{\paramRcRidge}{10^{-4}} % results/camera_ready/camera_ready_metrics.json:parameters.recipe.RcRidge
\newcommand{\paramRcWashout}{1000} % results/camera_ready/camera_ready_metrics.json:parameters.recipe.RcWashout
\newcommand{\paramRelThreshExp}{-4} % results/camera_ready/camera_ready_metrics.json:parameters.rel_thresh_exp
\newcommand{\paramRho}{28} % results/camera_ready/camera_ready_metrics.json:parameters.rho
\newcommand{\paramS}{2} % results/camera_ready/camera_ready_metrics.json:parameters.noise_amp
\newcommand{\paramSigma}{10} % results/camera_ready/camera_ready_metrics.json:parameters.sigma
\newcommand{\paramSigmaHidden}{32} % results/camera_ready/camera_ready_metrics.json:parameters.recipe.SigmaHidden
\newcommand{\paramSindyDegree}{2} % results/camera_ready/camera_ready_metrics.json:parameters.recipe.SindyDegree
\newcommand{\paramSindyThreshold}{0.1} % results/camera_ready/camera_ready_metrics.json:parameters.recipe.SindyThreshold
\newcommand{\paramT}{10000} % results/camera_ready/camera_ready_metrics.json:parameters.T_tu
\newcommand{\paramTExp}{4} % results/camera_ready/camera_ready_metrics.json:parameters.T_exp
\newcommand{\paramTSVD}{10^{-4}} % results/camera_ready/camera_ready_metrics.json:parameters.rel_thresh_exp (as 10^{x})
\newcommand{\paramTau}{0.1} % results/camera_ready/camera_ready_metrics.json:parameters.tau
\newcommand{\paramTen}{10} % results/camera_ready/camera_ready_metrics.json:parameters.ten
\newcommand{\paramThrCap}{0.99} % results/camera_ready/camera_ready_metrics.json:parameters.thresh_cap
\newcommand{\paramThrFloor}{0.5} % results/camera_ready/camera_ready_metrics.json:parameters.thresh_floor
\newcommand{\paramTrainTUMLPdrift}{10000} % results/camera_ready/camera_ready_metrics.json:parameters.train_tu.mlp_drift
\newcommand{\paramTrainTUMLPflow}{10000} % results/camera_ready/camera_ready_metrics.json:parameters.train_tu.mlp_flow
\newcommand{\paramTrainTUNODE}{10000} % results/camera_ready/camera_ready_metrics.json:parameters.train_tu.node
\newcommand{\paramTrainTUNSDE}{10000} % results/camera_ready/camera_ready_metrics.json:parameters.train_tu.nsde
\newcommand{\paramTrainTURC}{2000} % results/camera_ready/camera_ready_metrics.json:parameters.train_tu.rc_flowmap
\newcommand{\paramTrainTUSindy}{100} % results/camera_ready/camera_ready_metrics.json:parameters.train_tu.sindy
\newcommand{\paramTrainTUSindyFull}{10000} % results/camera_ready/camera_ready_metrics.json:parameters.train_tu.sindy_full
\newcommand{\paramTrainTUTransformer}{10000} % results/camera_ready/camera_ready_metrics.json:exploratory.transformer.appendix.recipe.train_tu
\newcommand{\paramTwo}{2} % results/camera_ready/camera_ready_metrics.json:parameters.two
\newcommand{\paramZero}{0} % results/camera_ready/camera_ready_metrics.json:parameters.zero
\newcommand{\pcNODEOne}{63} % results/camera_ready/camera_ready_metrics.json:models.node.final.per_mode_full.modes[0].c
\newcommand{\pcNSDEOne}{71} % results/camera_ready/camera_ready_metrics.json:models.nsde.final.per_mode_full.modes[0].c
\newcommand{\pcRCOne}{75} % results/camera_ready/camera_ready_metrics.json:models.rc_flowmap.final.per_mode_full.modes[0].c
\newcommand{\pcSindyFullOne}{0} % results/camera_ready/camera_ready_metrics.json:diagnostic.sindy_full.final.per_mode_full.modes[0].c
\newcommand{\pctBandloo}{2} % results/camera_ready/camera_ready_metrics.json:criteria.band.loo_inside (100 (1 - x), integer %)
\newcommand{\pctChiGzeroH}{0} % results/camera_ready/camera_ready_metrics.json:models.true_rk4_h.final.flags.rate_chi (integer %)
\newcommand{\pctChiGzeroHSig}{30} % results/camera_ready/camera_ready_metrics.json:models.true_rk4_h_sig.final.flags.rate_chi (integer %)
\newcommand{\pctChiMLPdrift}{91} % results/camera_ready/camera_ready_metrics.json:models.mlp_drift.final.flags.rate_chi (integer %)
\newcommand{\pctChiMLPflow}{99} % results/camera_ready/camera_ready_metrics.json:models.mlp_flow.final.flags.rate_chi (integer %)
\newcommand{\pctChiNODE}{9} % results/camera_ready/camera_ready_metrics.json:models.node.final.flags.rate_chi (integer %)
\newcommand{\pctChiNSDE}{9} % results/camera_ready/camera_ready_metrics.json:models.nsde.final.flags.rate_chi (integer %)
\newcommand{\pctChiRC}{8} % results/camera_ready/camera_ready_metrics.json:models.rc_flowmap.final.flags.rate_chi (integer %)
\newcommand{\pctChiSindyFull}{0} % results/camera_ready/camera_ready_metrics.json:diagnostic.sindy_full.final.flags.rate_chi (integer %)
\newcommand{\pctChiTransformer}{13} % results/camera_ready/camera_ready_metrics.json:exploratory.transformer.final.flags.rate_chi (integer %)
\newcommand{\pctChiTrue}{0} % results/camera_ready/camera_ready_metrics.json:models.true.final.flags.rate_chi (integer %)
\newcommand{\pctInvGzeroH}{100} % results/camera_ready/camera_ready_metrics.json:models.true_rk4_h.final.flags.rate_inv (integer %)
\newcommand{\pctInvGzeroHSig}{100} % results/camera_ready/camera_ready_metrics.json:models.true_rk4_h_sig.final.flags.rate_inv (integer %)
\newcommand{\pctInvMLPdrift}{100} % results/camera_ready/camera_ready_metrics.json:models.mlp_drift.final.flags.rate_inv (integer %)
\newcommand{\pctInvMLPflow}{100} % results/camera_ready/camera_ready_metrics.json:models.mlp_flow.final.flags.rate_inv (integer %)
\newcommand{\pctInvNODE}{93} % results/camera_ready/camera_ready_metrics.json:models.node.final.flags.rate_inv (integer %)
\newcommand{\pctInvNSDE}{100} % results/camera_ready/camera_ready_metrics.json:models.nsde.final.flags.rate_inv (integer %)
\newcommand{\pctInvRC}{63} % results/camera_ready/camera_ready_metrics.json:models.rc_flowmap.final.flags.rate_inv (integer %)
\newcommand{\pctInvSindyFull}{20} % results/camera_ready/camera_ready_metrics.json:diagnostic.sindy_full.final.flags.rate_inv (integer %)
\newcommand{\pctInvTransformer}{79} % results/camera_ready/camera_ready_metrics.json:exploratory.transformer.final.flags.rate_inv (integer %)
\newcommand{\pctInvTrue}{4} % results/camera_ready/camera_ready_metrics.json:models.true.final.flags.rate_inv (integer %)
\newcommand{\pctOscGzeroH}{40} % results/camera_ready/camera_ready_metrics.json:models.true_rk4_h.final.flags.rate_osc (integer %)
\newcommand{\pctOscGzeroHSig}{40} % results/camera_ready/camera_ready_metrics.json:models.true_rk4_h_sig.final.flags.rate_osc (integer %)
\newcommand{\pctOscMLPdrift}{100} % results/camera_ready/camera_ready_metrics.json:models.mlp_drift.final.flags.rate_osc (integer %)
\newcommand{\pctOscMLPflow}{100} % results/camera_ready/camera_ready_metrics.json:models.mlp_flow.final.flags.rate_osc (integer %)
\newcommand{\pctOscNODE}{86} % results/camera_ready/camera_ready_metrics.json:models.node.final.flags.rate_osc (integer %)
\newcommand{\pctOscNSDE}{83} % results/camera_ready/camera_ready_metrics.json:models.nsde.final.flags.rate_osc (integer %)
\newcommand{\pctOscRC}{100} % results/camera_ready/camera_ready_metrics.json:models.rc_flowmap.final.flags.rate_osc (integer %)
\newcommand{\pctOscSindyFull}{10} % results/camera_ready/camera_ready_metrics.json:diagnostic.sindy_full.final.flags.rate_osc (integer %)
\newcommand{\pctOscTransformer}{73} % results/camera_ready/camera_ready_metrics.json:exploratory.transformer.final.flags.rate_osc (integer %)
\newcommand{\pctOscTrue}{2} % results/camera_ready/camera_ready_metrics.json:models.true.final.flags.rate_osc (integer %)
\newcommand{\pctOutAcfGzeroH}{60} % results/camera_ready/camera_ready_metrics.json:models.true_rk4_h.final.invariant.acf_distance (100 (1 - k_inside / n_scored), integer %)
\newcommand{\pctOutAcfGzeroHSig}{70} % results/camera_ready/camera_ready_metrics.json:models.true_rk4_h_sig.final.invariant.acf_distance (100 (1 - k_inside / n_scored), integer %)
\newcommand{\pctOutAcfMLPdrift}{100} % results/camera_ready/camera_ready_metrics.json:models.mlp_drift.final.invariant.acf_distance (100 (1 - k_inside / n_scored), integer %)
\newcommand{\pctOutAcfMLPflow}{93} % results/camera_ready/camera_ready_metrics.json:models.mlp_flow.final.invariant.acf_distance (100 (1 - k_inside / n_scored), integer %)
\newcommand{\pctOutAcfNODE}{56} % results/camera_ready/camera_ready_metrics.json:models.node.final.invariant.acf_distance (100 (1 - k_inside / n_scored), integer %)
\newcommand{\pctOutAcfNSDE}{69} % results/camera_ready/camera_ready_metrics.json:models.nsde.final.invariant.acf_distance (100 (1 - k_inside / n_scored), integer %)
\newcommand{\pctOutAcfRC}{27} % results/camera_ready/camera_ready_metrics.json:models.rc_flowmap.final.invariant.acf_distance (100 (1 - k_inside / n_scored), integer %)
\newcommand{\pctOutAcfSindyFull}{0} % results/camera_ready/camera_ready_metrics.json:diagnostic.sindy_full.final.invariant.acf_distance (100 (1 - k_inside / n_scored), integer %)
\newcommand{\pctOutAcfTransformer}{29} % results/camera_ready/camera_ready_metrics.json:exploratory.transformer.final.invariant.acf_distance (100 (1 - k_inside / n_scored), integer %)
\newcommand{\pctOutAcfTrue}{0} % results/camera_ready/camera_ready_metrics.json:models.true.final.invariant.acf_distance (100 (1 - k_inside / n_scored), integer %)
\newcommand{\pctOutFPlusGzeroH}{0} % results/camera_ready/camera_ready_metrics.json:models.true_rk4_h.final.invariant.f_plus (100 (1 - k_inside / n_scored), integer %)
\newcommand{\pctOutFPlusGzeroHSig}{0} % results/camera_ready/camera_ready_metrics.json:models.true_rk4_h_sig.final.invariant.f_plus (100 (1 - k_inside / n_scored), integer %)
\newcommand{\pctOutFPlusMLPdrift}{96} % results/camera_ready/camera_ready_metrics.json:models.mlp_drift.final.invariant.f_plus (100 (1 - k_inside / n_scored), integer %)
\newcommand{\pctOutFPlusMLPflow}{100} % results/camera_ready/camera_ready_metrics.json:models.mlp_flow.final.invariant.f_plus (100 (1 - k_inside / n_scored), integer %)
\newcommand{\pctOutFPlusNODE}{63} % results/camera_ready/camera_ready_metrics.json:models.node.final.invariant.f_plus (100 (1 - k_inside / n_scored), integer %)
\newcommand{\pctOutFPlusNSDE}{61} % results/camera_ready/camera_ready_metrics.json:models.nsde.final.invariant.f_plus (100 (1 - k_inside / n_scored), integer %)
\newcommand{\pctOutFPlusRC}{6} % results/camera_ready/camera_ready_metrics.json:models.rc_flowmap.final.invariant.f_plus (100 (1 - k_inside / n_scored), integer %)
\newcommand{\pctOutFPlusSindyFull}{20} % results/camera_ready/camera_ready_metrics.json:diagnostic.sindy_full.final.invariant.f_plus (100 (1 - k_inside / n_scored), integer %)
\newcommand{\pctOutFPlusTransformer}{64} % results/camera_ready/camera_ready_metrics.json:exploratory.transformer.final.invariant.f_plus (100 (1 - k_inside / n_scored), integer %)
\newcommand{\pctOutFPlusTrue}{0} % results/camera_ready/camera_ready_metrics.json:models.true.final.invariant.f_plus (100 (1 - k_inside / n_scored), integer %)
\newcommand{\pctOutSwitchGzeroH}{100} % results/camera_ready/camera_ready_metrics.json:models.true_rk4_h.final.invariant.switch_ratio (100 (1 - k_inside / n_scored), integer %)
\newcommand{\pctOutSwitchGzeroHSig}{100} % results/camera_ready/camera_ready_metrics.json:models.true_rk4_h_sig.final.invariant.switch_ratio (100 (1 - k_inside / n_scored), integer %)
\newcommand{\pctOutSwitchMLPdrift}{100} % results/camera_ready/camera_ready_metrics.json:models.mlp_drift.final.invariant.switch_ratio (100 (1 - k_inside / n_scored), integer %)
\newcommand{\pctOutSwitchMLPflow}{100} % results/camera_ready/camera_ready_metrics.json:models.mlp_flow.final.invariant.switch_ratio (100 (1 - k_inside / n_scored), integer %)
\newcommand{\pctOutSwitchNODE}{75} % results/camera_ready/camera_ready_metrics.json:models.node.final.invariant.switch_ratio (100 (1 - k_inside / n_scored), integer %)
\newcommand{\pctOutSwitchNSDE}{79} % results/camera_ready/camera_ready_metrics.json:models.nsde.final.invariant.switch_ratio (100 (1 - k_inside / n_scored), integer %)
\newcommand{\pctOutSwitchRC}{58} % results/camera_ready/camera_ready_metrics.json:models.rc_flowmap.final.invariant.switch_ratio (100 (1 - k_inside / n_scored), integer %)
\newcommand{\pctOutSwitchSindyFull}{0} % results/camera_ready/camera_ready_metrics.json:diagnostic.sindy_full.final.invariant.switch_ratio (100 (1 - k_inside / n_scored), integer %)
\newcommand{\pctOutSwitchTransformer}{53} % results/camera_ready/camera_ready_metrics.json:exploratory.transformer.final.invariant.switch_ratio (100 (1 - k_inside / n_scored), integer %)
\newcommand{\pctOutSwitchTrue}{4} % results/camera_ready/camera_ready_metrics.json:models.true.final.invariant.switch_ratio (100 (1 - k_inside / n_scored), integer %)
\newcommand{\pctOutWOneGzeroH}{0} % results/camera_ready/camera_ready_metrics.json:models.true_rk4_h.final.invariant.w1_x (100 (1 - k_inside / n_scored), integer %)
\newcommand{\pctOutWOneGzeroHSig}{0} % results/camera_ready/camera_ready_metrics.json:models.true_rk4_h_sig.final.invariant.w1_x (100 (1 - k_inside / n_scored), integer %)
\newcommand{\pctOutWOneMLPdrift}{100} % results/camera_ready/camera_ready_metrics.json:models.mlp_drift.final.invariant.w1_x (100 (1 - k_inside / n_scored), integer %)
\newcommand{\pctOutWOneMLPflow}{100} % results/camera_ready/camera_ready_metrics.json:models.mlp_flow.final.invariant.w1_x (100 (1 - k_inside / n_scored), integer %)
\newcommand{\pctOutWOneNODE}{58} % results/camera_ready/camera_ready_metrics.json:models.node.final.invariant.w1_x (100 (1 - k_inside / n_scored), integer %)
\newcommand{\pctOutWOneNSDE}{60} % results/camera_ready/camera_ready_metrics.json:models.nsde.final.invariant.w1_x (100 (1 - k_inside / n_scored), integer %)
\newcommand{\pctOutWOneRC}{3} % results/camera_ready/camera_ready_metrics.json:models.rc_flowmap.final.invariant.w1_x (100 (1 - k_inside / n_scored), integer %)
\newcommand{\pctOutWOneSindyFull}{10} % results/camera_ready/camera_ready_metrics.json:diagnostic.sindy_full.final.invariant.w1_x (100 (1 - k_inside / n_scored), integer %)
\newcommand{\pctOutWOneTransformer}{51} % results/camera_ready/camera_ready_metrics.json:exploratory.transformer.final.invariant.w1_x (100 (1 - k_inside / n_scored), integer %)
\newcommand{\pctOutWOneTrue}{0} % results/camera_ready/camera_ready_metrics.json:models.true.final.invariant.w1_x (100 (1 - k_inside / n_scored), integer %)
\newcommand{\pctSlowGzeroH}{10} % results/camera_ready/camera_ready_metrics.json:models.true_rk4_h.final.flags.rate_slow (integer %)
\newcommand{\pctSlowGzeroHSig}{0} % results/camera_ready/camera_ready_metrics.json:models.true_rk4_h_sig.final.flags.rate_slow (integer %)
\newcommand{\pctSlowMLPdrift}{100} % results/camera_ready/camera_ready_metrics.json:models.mlp_drift.final.flags.rate_slow (integer %)
\newcommand{\pctSlowMLPflow}{98} % results/camera_ready/camera_ready_metrics.json:models.mlp_flow.final.flags.rate_slow (integer %)
\newcommand{\pctSlowNODE}{24} % results/camera_ready/camera_ready_metrics.json:models.node.final.flags.rate_slow (integer %)
\newcommand{\pctSlowNSDE}{23} % results/camera_ready/camera_ready_metrics.json:models.nsde.final.flags.rate_slow (integer %)
\newcommand{\pctSlowRC}{76} % results/camera_ready/camera_ready_metrics.json:models.rc_flowmap.final.flags.rate_slow (integer %)
\newcommand{\pctSlowSindyFull}{0} % results/camera_ready/camera_ready_metrics.json:diagnostic.sindy_full.final.flags.rate_slow (integer %)
\newcommand{\pctSlowTransformer}{19} % results/camera_ready/camera_ready_metrics.json:exploratory.transformer.final.flags.rate_slow (integer %)
\newcommand{\pctSlowTrue}{6} % results/camera_ready/camera_ready_metrics.json:models.true.final.flags.rate_slow (integer %)
\newcommand{\refFreqFive}{9.98} % results/camera_ready/camera_ready_metrics.json:models.mlp_drift.final.per_mode_full.modes[4].freq
\newcommand{\refFreqFour}{9.85} % results/camera_ready/camera_ready_metrics.json:models.mlp_drift.final.per_mode_full.modes[3].freq
\newcommand{\refFreqOne}{8.42} % results/camera_ready/camera_ready_metrics.json:models.mlp_drift.final.per_mode_full.modes[0].freq
\newcommand{\refFreqSix}{16.68} % results/camera_ready/camera_ready_metrics.json:models.mlp_drift.final.per_mode_full.modes[5].freq
\newcommand{\refReThree}{-0.508} % results/camera_ready/camera_ready_metrics.json:models.mlp_drift.final.per_mode_full.modes[2].re
\newcommand{\refReTwo}{-0.366} % results/camera_ready/camera_ready_metrics.json:models.mlp_drift.final.per_mode_full.modes[1].re
\newcommand{\relaxGapXPct}{4.8} % results/camera_ready/camera_ready_metrics.json:claims.relaxGapXPct
\newcommand{\rmseFullX}{0.014} % results/camera_ready/camera_ready_metrics.json:appendix.inputs.drift_rmse.sindy_full[0]
\newcommand{\rmseFullY}{0.070} % results/camera_ready/camera_ready_metrics.json:appendix.inputs.drift_rmse.sindy_full[1]
\newcommand{\rmseFullZ}{0.043} % results/camera_ready/camera_ready_metrics.json:appendix.inputs.drift_rmse.sindy_full[2]
\newcommand{\rmseHundredX}{0.717} % results/camera_ready/camera_ready_metrics.json:appendix.inputs.drift_rmse.sindy100[0]
\newcommand{\rmseHundredY}{0.537} % results/camera_ready/camera_ready_metrics.json:appendix.inputs.drift_rmse.sindy100[1]
\newcommand{\rmseHundredZ}{0.647} % results/camera_ready/camera_ready_metrics.json:appendix.inputs.drift_rmse.sindy100[2]
\newcommand{\rmseMLPdriftX}{2.048} % results/camera_ready/camera_ready_metrics.json:appendix.inputs.drift_rmse.mlp_drift_unified[0]
\newcommand{\rmseMLPdriftY}{4.399} % results/camera_ready/camera_ready_metrics.json:appendix.inputs.drift_rmse.mlp_drift_unified[1]
\newcommand{\rmseMLPdriftZ}{4.213} % results/camera_ready/camera_ready_metrics.json:appendix.inputs.drift_rmse.mlp_drift_unified[2]
\newcommand{\rmseNODEXmed}{0.577} % results/camera_ready/camera_ready_metrics.json:appendix.inputs.drift_rmse.node.median[0]
\newcommand{\rmseNODEYmed}{1.125} % results/camera_ready/camera_ready_metrics.json:appendix.inputs.drift_rmse.node.median[1]
\newcommand{\rmseNODEZmed}{0.961} % results/camera_ready/camera_ready_metrics.json:appendix.inputs.drift_rmse.node.median[2]
\newcommand{\rmseNODEn}{8} % results/camera_ready/camera_ready_metrics.json:appendix.inputs.drift_rmse.node.n_train_seeds
\newcommand{\rmseSecantLeadX}{2.439} % results/camera_ready/camera_ready_metrics.json:appendix.inputs.drift_rmse.secant_leading[0]
\newcommand{\rmseSecantLeadY}{4.740} % results/camera_ready/camera_ready_metrics.json:appendix.inputs.drift_rmse.secant_leading[1]
\newcommand{\rmseSecantLeadZ}{5.056} % results/camera_ready/camera_ready_metrics.json:appendix.inputs.drift_rmse.secant_leading[2]
\newcommand{\rmseSecantX}{2.192} % results/camera_ready/camera_ready_metrics.json:appendix.inputs.drift_rmse.secant_exact[0]
\newcommand{\rmseSecantY}{4.264} % results/camera_ready/camera_ready_metrics.json:appendix.inputs.drift_rmse.secant_exact[1]
\newcommand{\rmseSecantZ}{4.541} % results/camera_ready/camera_ready_metrics.json:appendix.inputs.drift_rmse.secant_exact[2]
\newcommand{\seedCalibOne}{1000} % results/camera_ready/camera_ready_metrics.json:parameters.seeds.calib[0]
\newcommand{\seedCalibTwo}{1015} % results/camera_ready/camera_ready_metrics.json:parameters.seeds.calib[1]
\newcommand{\seedEtwoRolloutOne}{0} % results/camera_ready/camera_ready_metrics.json:parameters.seeds.e2_rollout[0]
\newcommand{\seedEtwoRolloutTwo}{9} % results/camera_ready/camera_ready_metrics.json:parameters.seeds.e2_rollout[1]
\newcommand{\seedEtwoTrainOne}{0} % results/camera_ready/camera_ready_metrics.json:parameters.seeds.e2_train[0]
\newcommand{\seedEtwoTrainTwo}{7} % results/camera_ready/camera_ready_metrics.json:parameters.seeds.e2_train[1]
\newcommand{\seedFloorOne}{5000} % results/camera_ready/camera_ready_metrics.json:parameters.seeds.floor[0]
\newcommand{\seedFloorTwo}{5049} % results/camera_ready/camera_ready_metrics.json:parameters.seeds.floor[1]
\newcommand{\seedGrid}{909} % results/camera_ready/camera_ready_metrics.json:parameters.seeds.grid
\newcommand{\seedGzeroHOne}{40001} % results/camera_ready/camera_ready_metrics.json:parameters.seeds.g0_h[0]
\newcommand{\seedGzeroHThree}{41009} % results/camera_ready/camera_ready_metrics.json:parameters.seeds.g0_h[2]
\newcommand{\seedGzeroHTwo}{41001} % results/camera_ready/camera_ready_metrics.json:parameters.seeds.g0_h[1]
\newcommand{\seedGzeroSigOne}{40002} % results/camera_ready/camera_ready_metrics.json:parameters.seeds.g0_sig[0]
\newcommand{\seedGzeroSigThree}{41019} % results/camera_ready/camera_ready_metrics.json:parameters.seeds.g0_sig[2]
\newcommand{\seedGzeroSigTwo}{41011} % results/camera_ready/camera_ready_metrics.json:parameters.seeds.g0_sig[1]
\newcommand{\seedNodeInitOne}{31000} % results/camera_ready/camera_ready_metrics.json:parameters.seeds.node_init[0]
\newcommand{\seedNodeInitTwo}{31007} % results/camera_ready/camera_ready_metrics.json:parameters.seeds.node_init[1]
\newcommand{\seedNsdeInitOne}{31100} % results/camera_ready/camera_ready_metrics.json:parameters.seeds.nsde_init[0]
\newcommand{\seedNsdeInitTwo}{31107} % results/camera_ready/camera_ready_metrics.json:parameters.seeds.nsde_init[1]
\newcommand{\seedRef}{101} % results/camera_ready/camera_ready_metrics.json:parameters.seeds.ref
\newcommand{\seedTrfInitOne}{32000} % results/camera_ready/camera_ready_metrics.json:exploratory.transformer.appendix.recipe.init_seeds[0]
\newcommand{\seedTrfInitTwo}{32007} % results/camera_ready/camera_ready_metrics.json:exploratory.transformer.appendix.recipe.init_seeds[1]
\newcommand{\seedTruthCert}{40000} % results/camera_ready/camera_ready_metrics.json:parameters.seeds.truth_cert
\newcommand{\seedTruthExtOne}{42001} % results/camera_ready/camera_ready_metrics.json:parameters.seeds.truth_ext[0]
\newcommand{\seedTruthExtTwo}{42049} % results/camera_ready/camera_ready_metrics.json:parameters.seeds.truth_ext[1]
\newcommand{\shareAbsLeadNODEmed}{0.04} % results/camera_ready/camera_ready_metrics.json:appendix.chi0_shares.node.abs_lead.median
\newcommand{\shareAbsLeadNSDEmed}{0.03} % results/camera_ready/camera_ready_metrics.json:appendix.chi0_shares.nsde.abs_lead.median
\newcommand{\shareAbsLeadRCmed}{0.07} % results/camera_ready/camera_ready_metrics.json:appendix.chi0_shares.rc_flowmap.abs_lead.median
\newcommand{\shareAbsLeadSindymed}{0.04} % results/camera_ready/camera_ready_metrics.json:appendix.chi0_shares.sindy.abs_lead.median
\newcommand{\shareAbsLeadTruemed}{0.06} % results/camera_ready/camera_ready_metrics.json:appendix.chi0_shares.true.abs_lead.median
\newcommand{\shareLeadNODEmax}{0.50} % results/camera_ready/camera_ready_metrics.json:appendix.chi0_shares.node.lead.max
\newcommand{\shareLeadNODEmed}{0.24} % results/camera_ready/camera_ready_metrics.json:appendix.chi0_shares.node.lead.median
\newcommand{\shareLeadNODEmin}{-0.36} % results/camera_ready/camera_ready_metrics.json:appendix.chi0_shares.node.lead.min
\newcommand{\shareLeadNSDEmax}{0.57} % results/camera_ready/camera_ready_metrics.json:appendix.chi0_shares.nsde.lead.max
\newcommand{\shareLeadNSDEmed}{0.08} % results/camera_ready/camera_ready_metrics.json:appendix.chi0_shares.nsde.lead.median
\newcommand{\shareLeadNSDEmin}{-0.41} % results/camera_ready/camera_ready_metrics.json:appendix.chi0_shares.nsde.lead.min
\newcommand{\shareLeadRCmax}{0.60} % results/camera_ready/camera_ready_metrics.json:appendix.chi0_shares.rc_flowmap.lead.max
\newcommand{\shareLeadRCmed}{0.29} % results/camera_ready/camera_ready_metrics.json:appendix.chi0_shares.rc_flowmap.lead.median
\newcommand{\shareLeadRCmin}{-0.06} % results/camera_ready/camera_ready_metrics.json:appendix.chi0_shares.rc_flowmap.lead.min
\newcommand{\shareLeadRefmed}{0.34} % results/camera_ready/camera_ready_metrics.json:appendix.chi0_shares.ref.lead.median
\newcommand{\shareLeadSindymax}{0.40} % results/camera_ready/camera_ready_metrics.json:appendix.chi0_shares.sindy.lead.max
\newcommand{\shareLeadSindymed}{-0.06} % results/camera_ready/camera_ready_metrics.json:appendix.chi0_shares.sindy.lead.median
\newcommand{\shareLeadSindymin}{-0.52} % results/camera_ready/camera_ready_metrics.json:appendix.chi0_shares.sindy.lead.min
\newcommand{\shareLeadTruemax}{0.60} % results/camera_ready/camera_ready_metrics.json:appendix.chi0_shares.true.lead.max
\newcommand{\shareLeadTruemed}{0.33} % results/camera_ready/camera_ready_metrics.json:appendix.chi0_shares.true.lead.median
\newcommand{\shareLeadTruemin}{-0.15} % results/camera_ready/camera_ready_metrics.json:appendix.chi0_shares.true.lead.min
\newcommand{\shareNNODE}{10} % results/camera_ready/camera_ready_metrics.json:appendix.chi0_shares.node.n
\newcommand{\shareNNSDE}{10} % results/camera_ready/camera_ready_metrics.json:appendix.chi0_shares.nsde.n
\newcommand{\shareNRC}{10} % results/camera_ready/camera_ready_metrics.json:appendix.chi0_shares.rc_flowmap.n
\newcommand{\shareNRef}{1} % results/camera_ready/camera_ready_metrics.json:appendix.chi0_shares.ref.n
\newcommand{\shareNSindy}{10} % results/camera_ready/camera_ready_metrics.json:appendix.chi0_shares.sindy.n
\newcommand{\shareNTrue}{10} % results/camera_ready/camera_ready_metrics.json:appendix.chi0_shares.true.n
\newcommand{\slowflagmax}{0.76} % results/camera_ready/camera_ready_metrics.json:claims.slowflagmax
\newcommand{\slowflagmin}{0.00} % results/camera_ready/camera_ready_metrics.json:claims.slowflagmin
\newcommand{\stsExhibitFisherp}{2.6\times10^{-7}} % results/camera_ready/camera_ready_metrics.json:slow_threshold_sensitivity.exhibit.fisher_p
\newcommand{\stsGzeroHSigflagDseventeen}{0.00} % results/camera_ready/camera_ready_metrics.json:slow_threshold_sensitivity.models.true_rk4_h_sig.rate_d17
\newcommand{\stsGzeroHSigflagThr}{0.00} % results/camera_ready/camera_ready_metrics.json:slow_threshold_sensitivity.models.true_rk4_h_sig.rate_thr
\newcommand{\stsGzeroHSign}{10} % results/camera_ready/camera_ready_metrics.json:slow_threshold_sensitivity.models.true_rk4_h_sig.n
\newcommand{\stsGzeroHSigpassOne}{1.00} % results/camera_ready/camera_ready_metrics.json:slow_threshold_sensitivity.models.true_rk4_h_sig.pass_thr[0]
\newcommand{\stsGzeroHSigpassTwo}{1.00} % results/camera_ready/camera_ready_metrics.json:slow_threshold_sensitivity.models.true_rk4_h_sig.pass_thr[1]
\newcommand{\stsGzeroHflagDseventeen}{0.10} % results/camera_ready/camera_ready_metrics.json:slow_threshold_sensitivity.models.true_rk4_h.rate_d17
\newcommand{\stsGzeroHflagThr}{0.00} % results/camera_ready/camera_ready_metrics.json:slow_threshold_sensitivity.models.true_rk4_h.rate_thr
\newcommand{\stsGzeroHn}{10} % results/camera_ready/camera_ready_metrics.json:slow_threshold_sensitivity.models.true_rk4_h.n
\newcommand{\stsGzeroHpassOne}{1.00} % results/camera_ready/camera_ready_metrics.json:slow_threshold_sensitivity.models.true_rk4_h.pass_thr[0]
\newcommand{\stsGzeroHpassTwo}{1.00} % results/camera_ready/camera_ready_metrics.json:slow_threshold_sensitivity.models.true_rk4_h.pass_thr[1]
\newcommand{\stsMLPdriftflagDseventeen}{1.00} % results/camera_ready/camera_ready_metrics.json:slow_threshold_sensitivity.models.mlp_drift.rate_d17
\newcommand{\stsMLPdriftflagThr}{1.00} % results/camera_ready/camera_ready_metrics.json:slow_threshold_sensitivity.models.mlp_drift.rate_thr
\newcommand{\stsMLPdriftn}{80} % results/camera_ready/camera_ready_metrics.json:slow_threshold_sensitivity.models.mlp_drift.n
\newcommand{\stsMLPdriftpassOne}{0.00} % results/camera_ready/camera_ready_metrics.json:slow_threshold_sensitivity.models.mlp_drift.pass_thr[0]
\newcommand{\stsMLPdriftpassTwo}{0.00} % results/camera_ready/camera_ready_metrics.json:slow_threshold_sensitivity.models.mlp_drift.pass_thr[1]
\newcommand{\stsMLPflowflagDseventeen}{0.97} % results/camera_ready/camera_ready_metrics.json:slow_threshold_sensitivity.models.mlp_flow.rate_d17
\newcommand{\stsMLPflowflagThr}{0.94} % results/camera_ready/camera_ready_metrics.json:slow_threshold_sensitivity.models.mlp_flow.rate_thr
\newcommand{\stsMLPflown}{80} % results/camera_ready/camera_ready_metrics.json:slow_threshold_sensitivity.models.mlp_flow.n
\newcommand{\stsMLPflowpassOne}{0.10} % results/camera_ready/camera_ready_metrics.json:slow_threshold_sensitivity.models.mlp_flow.pass_thr[0]
\newcommand{\stsMLPflowpassTwo}{0.19} % results/camera_ready/camera_ready_metrics.json:slow_threshold_sensitivity.models.mlp_flow.pass_thr[1]
\newcommand{\stsNODEflagDseventeen}{0.24} % results/camera_ready/camera_ready_metrics.json:slow_threshold_sensitivity.models.node.rate_d17
\newcommand{\stsNODEflagThr}{0.07} % results/camera_ready/camera_ready_metrics.json:slow_threshold_sensitivity.models.node.rate_thr
\newcommand{\stsNODEn}{80} % results/camera_ready/camera_ready_metrics.json:slow_threshold_sensitivity.models.node.n
\newcommand{\stsNODEpassOne}{0.95} % results/camera_ready/camera_ready_metrics.json:slow_threshold_sensitivity.models.node.pass_thr[0]
\newcommand{\stsNODEpassTwo}{0.94} % results/camera_ready/camera_ready_metrics.json:slow_threshold_sensitivity.models.node.pass_thr[1]
\newcommand{\stsNSDEflagDseventeen}{0.23} % results/camera_ready/camera_ready_metrics.json:slow_threshold_sensitivity.models.nsde.rate_d17
\newcommand{\stsNSDEflagThr}{0.09} % results/camera_ready/camera_ready_metrics.json:slow_threshold_sensitivity.models.nsde.rate_thr
\newcommand{\stsNSDEn}{80} % results/camera_ready/camera_ready_metrics.json:slow_threshold_sensitivity.models.nsde.n
\newcommand{\stsNSDEpassOne}{0.95} % results/camera_ready/camera_ready_metrics.json:slow_threshold_sensitivity.models.nsde.pass_thr[0]
\newcommand{\stsNSDEpassTwo}{0.91} % results/camera_ready/camera_ready_metrics.json:slow_threshold_sensitivity.models.nsde.pass_thr[1]
\newcommand{\stsRCflagDseventeen}{0.76} % results/camera_ready/camera_ready_metrics.json:slow_threshold_sensitivity.models.rc_flowmap.rate_d17
\newcommand{\stsRCflagThr}{0.55} % results/camera_ready/camera_ready_metrics.json:slow_threshold_sensitivity.models.rc_flowmap.rate_thr
\newcommand{\stsRCn}{78} % results/camera_ready/camera_ready_metrics.json:slow_threshold_sensitivity.models.rc_flowmap.n
\newcommand{\stsRCnX}{26} % results/camera_ready/camera_ready_metrics.json:slow_threshold_sensitivity.exhibit.RCnX
\newcommand{\stsRCnXslowfail}{14} % results/camera_ready/camera_ready_metrics.json:slow_threshold_sensitivity.exhibit.RCnXslowfail
\newcommand{\stsRCpassOne}{0.45} % results/camera_ready/camera_ready_metrics.json:slow_threshold_sensitivity.models.rc_flowmap.pass_thr[0]
\newcommand{\stsRCpassTwo}{0.62} % results/camera_ready/camera_ready_metrics.json:slow_threshold_sensitivity.models.rc_flowmap.pass_thr[1]
\newcommand{\stsSindyFullflagDseventeen}{0.00} % results/camera_ready/camera_ready_metrics.json:slow_threshold_sensitivity.models.sindy_full.rate_d17
\newcommand{\stsSindyFullflagThr}{0.00} % results/camera_ready/camera_ready_metrics.json:slow_threshold_sensitivity.models.sindy_full.rate_thr
\newcommand{\stsSindyFulln}{10} % results/camera_ready/camera_ready_metrics.json:slow_threshold_sensitivity.models.sindy_full.n
\newcommand{\stsSindyFullpassOne}{1.00} % results/camera_ready/camera_ready_metrics.json:slow_threshold_sensitivity.models.sindy_full.pass_thr[0]
\newcommand{\stsSindyFullpassTwo}{1.00} % results/camera_ready/camera_ready_metrics.json:slow_threshold_sensitivity.models.sindy_full.pass_thr[1]
\newcommand{\stsSindyflagDseventeen}{0.00} % results/camera_ready/camera_ready_metrics.json:slow_threshold_sensitivity.models.sindy.rate_d17
\newcommand{\stsSindyflagThr}{0.00} % results/camera_ready/camera_ready_metrics.json:slow_threshold_sensitivity.models.sindy.rate_thr
\newcommand{\stsSindyn}{10} % results/camera_ready/camera_ready_metrics.json:slow_threshold_sensitivity.models.sindy.n
\newcommand{\stsSindypassOne}{1.00} % results/camera_ready/camera_ready_metrics.json:slow_threshold_sensitivity.models.sindy.pass_thr[0]
\newcommand{\stsSindypassTwo}{1.00} % results/camera_ready/camera_ready_metrics.json:slow_threshold_sensitivity.models.sindy.pass_thr[1]
\newcommand{\stsThrDseventeenOne}{0.9450} % results/camera_ready/camera_ready_metrics.json:slow_threshold_sensitivity.threshold_d17[0]
\newcommand{\stsThrDseventeenTwo}{0.9216} % results/camera_ready/camera_ready_metrics.json:slow_threshold_sensitivity.threshold_d17[1]
\newcommand{\stsThrOne}{0.9257} % results/camera_ready/camera_ready_metrics.json:slow_threshold_sensitivity.threshold[0]
\newcommand{\stsThrTwo}{0.8693} % results/camera_ready/camera_ready_metrics.json:slow_threshold_sensitivity.threshold[1]
\newcommand{\stsTransformerflagDseventeen}{0.19} % results/camera_ready/camera_ready_metrics.json:exploratory.transformer.appendix.sts.rate_d17
\newcommand{\stsTransformerflagThr}{0.09} % results/camera_ready/camera_ready_metrics.json:exploratory.transformer.appendix.sts.rate_thr
\newcommand{\stsTransformern}{80} % results/camera_ready/camera_ready_metrics.json:exploratory.transformer.appendix.sts.n
\newcommand{\stsTransformerpassOne}{0.91} % results/camera_ready/camera_ready_metrics.json:exploratory.transformer.appendix.sts.pass_thr[0]
\newcommand{\stsTransformerpassTwo}{0.97} % results/camera_ready/camera_ready_metrics.json:exploratory.transformer.appendix.sts.pass_thr[1]
\newcommand{\stsTrueflagDseventeen}{0.06} % results/camera_ready/camera_ready_metrics.json:slow_threshold_sensitivity.models.true.rate_d17
\newcommand{\stsTrueflagThr}{0.02} % results/camera_ready/camera_ready_metrics.json:slow_threshold_sensitivity.models.true.rate_thr
\newcommand{\stsTruen}{50} % results/camera_ready/camera_ready_metrics.json:slow_threshold_sensitivity.models.true.n
\newcommand{\stsTruenX}{48} % results/camera_ready/camera_ready_metrics.json:slow_threshold_sensitivity.exhibit.TruenX
\newcommand{\stsTruenXslowfail}{1} % results/camera_ready/camera_ready_metrics.json:slow_threshold_sensitivity.exhibit.TruenXslowfail
\newcommand{\stsTruepassOne}{0.98} % results/camera_ready/camera_ready_metrics.json:slow_threshold_sensitivity.models.true.pass_thr[0]
\newcommand{\stsTruepassTwo}{0.98} % results/camera_ready/camera_ready_metrics.json:slow_threshold_sensitivity.models.true.pass_thr[1]
\newcommand{\stsTruthN}{50} % results/camera_ready/camera_ready_metrics.json:slow_threshold_sensitivity.n_truth
\newcommand{\sweepEpochsFour}{25} % results/camera_ready/camera_ready_metrics.json:appendix.node_sweep.rows[3].n_epochs_run
\newcommand{\sweepEpochsOne}{25} % results/camera_ready/camera_ready_metrics.json:appendix.node_sweep.rows[0].n_epochs_run
\newcommand{\sweepEpochsThree}{14} % results/camera_ready/camera_ready_metrics.json:appendix.node_sweep.rows[2].n_epochs_run
\newcommand{\sweepEpochsTwo}{25} % results/camera_ready/camera_ready_metrics.json:appendix.node_sweep.rows[1].n_epochs_run
\newcommand{\sweepLFour}{20} % results/camera_ready/camera_ready_metrics.json:appendix.node_sweep.rows[3].L
\newcommand{\sweepLOne}{5} % results/camera_ready/camera_ready_metrics.json:appendix.node_sweep.rows[0].L
\newcommand{\sweepLThree}{20} % results/camera_ready/camera_ready_metrics.json:appendix.node_sweep.rows[2].L
\newcommand{\sweepLTwo}{5} % results/camera_ready/camera_ready_metrics.json:appendix.node_sweep.rows[1].L
\newcommand{\sweepLrFour}{0.0003} % results/camera_ready/camera_ready_metrics.json:appendix.node_sweep.rows[3].lr
\newcommand{\sweepLrOne}{0.001} % results/camera_ready/camera_ready_metrics.json:appendix.node_sweep.rows[0].lr
\newcommand{\sweepLrThree}{0.001} % results/camera_ready/camera_ready_metrics.json:appendix.node_sweep.rows[2].lr
\newcommand{\sweepLrTwo}{0.0003} % results/camera_ready/camera_ready_metrics.json:appendix.node_sweep.rows[1].lr
\newcommand{\sweepMseLfiveFour}{1.551} % results/camera_ready/camera_ready_metrics.json:appendix.node_sweep.rows[3].val_mse_L5 (x 1e3)
\newcommand{\sweepMseLfiveOne}{1.543} % results/camera_ready/camera_ready_metrics.json:appendix.node_sweep.rows[0].val_mse_L5 (x 1e3)
\newcommand{\sweepMseLfiveThree}{1.555} % results/camera_ready/camera_ready_metrics.json:appendix.node_sweep.rows[2].val_mse_L5 (x 1e3)
\newcommand{\sweepMseLfiveTwo}{1.547} % results/camera_ready/camera_ready_metrics.json:appendix.node_sweep.rows[1].val_mse_L5 (x 1e3)
\newcommand{\sweepMseLoneFour}{0.548} % results/camera_ready/camera_ready_metrics.json:appendix.node_sweep.rows[3].val_mse_L1 (x 1e3)
\newcommand{\sweepMseLoneOne}{0.547} % results/camera_ready/camera_ready_metrics.json:appendix.node_sweep.rows[0].val_mse_L1 (x 1e3)
\newcommand{\sweepMseLoneThree}{0.549} % results/camera_ready/camera_ready_metrics.json:appendix.node_sweep.rows[2].val_mse_L1 (x 1e3)
\newcommand{\sweepMseLoneTwo}{0.548} % results/camera_ready/camera_ready_metrics.json:appendix.node_sweep.rows[1].val_mse_L1 (x 1e3)
\newcommand{\sweepMseLtwentyFour}{6.239} % results/camera_ready/camera_ready_metrics.json:appendix.node_sweep.rows[3].val_mse_L20 (x 1e3)
\newcommand{\sweepMseLtwentyOne}{6.225} % results/camera_ready/camera_ready_metrics.json:appendix.node_sweep.rows[0].val_mse_L20 (x 1e3)
\newcommand{\sweepMseLtwentyThree}{6.268} % results/camera_ready/camera_ready_metrics.json:appendix.node_sweep.rows[2].val_mse_L20 (x 1e3)
\newcommand{\sweepMseLtwentyTwo}{6.242} % results/camera_ready/camera_ready_metrics.json:appendix.node_sweep.rows[1].val_mse_L20 (x 1e3)
\newcommand{\sweepRmseXFour}{0.800} % results/camera_ready/camera_ready_metrics.json:appendix.node_sweep.rows[3].drift_rmse[0]
\newcommand{\sweepRmseXOne}{0.490} % results/camera_ready/camera_ready_metrics.json:appendix.node_sweep.rows[0].drift_rmse[0]
\newcommand{\sweepRmseXThree}{0.841} % results/camera_ready/camera_ready_metrics.json:appendix.node_sweep.rows[2].drift_rmse[0]
\newcommand{\sweepRmseXTwo}{0.648} % results/camera_ready/camera_ready_metrics.json:appendix.node_sweep.rows[1].drift_rmse[0]
\newcommand{\sweepRmseYFour}{1.229} % results/camera_ready/camera_ready_metrics.json:appendix.node_sweep.rows[3].drift_rmse[1]
\newcommand{\sweepRmseYOne}{0.921} % results/camera_ready/camera_ready_metrics.json:appendix.node_sweep.rows[0].drift_rmse[1]
\newcommand{\sweepRmseYThree}{1.483} % results/camera_ready/camera_ready_metrics.json:appendix.node_sweep.rows[2].drift_rmse[1]
\newcommand{\sweepRmseYTwo}{1.078} % results/camera_ready/camera_ready_metrics.json:appendix.node_sweep.rows[1].drift_rmse[1]
\newcommand{\sweepRmseZFour}{1.260} % results/camera_ready/camera_ready_metrics.json:appendix.node_sweep.rows[3].drift_rmse[2]
\newcommand{\sweepRmseZOne}{0.877} % results/camera_ready/camera_ready_metrics.json:appendix.node_sweep.rows[0].drift_rmse[2]
\newcommand{\sweepRmseZThree}{1.479} % results/camera_ready/camera_ready_metrics.json:appendix.node_sweep.rows[2].drift_rmse[2]
\newcommand{\sweepRmseZTwo}{0.986} % results/camera_ready/camera_ready_metrics.json:appendix.node_sweep.rows[1].drift_rmse[2]
\newcommand{\sweepTrainSFour}{601} % results/camera_ready/camera_ready_metrics.json:appendix.node_sweep.rows[3].train_s (rounded)
\newcommand{\sweepTrainSOne}{177} % results/camera_ready/camera_ready_metrics.json:appendix.node_sweep.rows[0].train_s (rounded)
\newcommand{\sweepTrainSThree}{342} % results/camera_ready/camera_ready_metrics.json:appendix.node_sweep.rows[2].train_s (rounded)
\newcommand{\sweepTrainSTwo}{177} % results/camera_ready/camera_ready_metrics.json:appendix.node_sweep.rows[1].train_s (rounded)
\newcommand{\switchRatioGzeroHSigmax}{0.968} % results/camera_ready/camera_ready_metrics.json:models.true_rk4_h_sig.metrics.switch_ratio.max
\newcommand{\switchRatioGzeroHSigmed}{0.948} % results/camera_ready/camera_ready_metrics.json:models.true_rk4_h_sig.metrics.switch_ratio.median
\newcommand{\switchRatioGzeroHmax}{0.975} % results/camera_ready/camera_ready_metrics.json:models.true_rk4_h.metrics.switch_ratio.max
\newcommand{\switchRatioGzeroHmed}{0.953} % results/camera_ready/camera_ready_metrics.json:models.true_rk4_h.metrics.switch_ratio.median
\newcommand{\switchRatioNODEmed}{1.05} % results/camera_ready/camera_ready_metrics.json:models.node.metrics.switch_ratio.median
\newcommand{\switchRatioNSDEmed}{1.05} % results/camera_ready/camera_ready_metrics.json:models.nsde.metrics.switch_ratio.median
\newcommand{\tabBandhi}{0.309} % results/camera_ready/camera_ready_metrics.json:criteria.band.hi (values table)
\newcommand{\tabBandlo}{0.162} % results/camera_ready/camera_ready_metrics.json:criteria.band.lo (values table)
\newcommand{\tabChiMeanGzeroH}{0.253} % results/camera_ready/camera_ready_metrics.json:models.true_rk4_h.metrics.chi0.mean
\newcommand{\tabChiMeanGzeroHSig}{0.281} % results/camera_ready/camera_ready_metrics.json:models.true_rk4_h_sig.metrics.chi0.mean
\newcommand{\tabChiMeanMLPdrift}{1.045} % results/camera_ready/camera_ready_metrics.json:models.mlp_drift.metrics.chi0.mean
\newcommand{\tabChiMeanMLPflow}{0.636} % results/camera_ready/camera_ready_metrics.json:models.mlp_flow.metrics.chi0.mean
\newcommand{\tabChiMeanNODE}{0.225} % results/camera_ready/camera_ready_metrics.json:models.node.metrics.chi0.mean
\newcommand{\tabChiMeanNSDE}{0.221} % results/camera_ready/camera_ready_metrics.json:models.nsde.metrics.chi0.mean
\newcommand{\tabChiMeanRC}{0.203} % results/camera_ready/camera_ready_metrics.json:models.rc_flowmap.metrics.chi0.mean
\newcommand{\tabChiMeanSindyFull}{0.215} % results/camera_ready/camera_ready_metrics.json:diagnostic.sindy_full.metrics.chi0.mean
\newcommand{\tabChiMeanTransformer}{0.200} % results/camera_ready/camera_ready_metrics.json:exploratory.transformer.metrics.chi0.mean
\newcommand{\tabChiMeanTrue}{0.224} % results/camera_ready/camera_ready_metrics.json:models.true.metrics.chi0.mean
\newcommand{\tabChiSdGzeroH}{0.043} % results/camera_ready/camera_ready_metrics.json:models.true_rk4_h.metrics.chi0.sd
\newcommand{\tabChiSdGzeroHSig}{0.043} % results/camera_ready/camera_ready_metrics.json:models.true_rk4_h_sig.metrics.chi0.sd
\newcommand{\tabChiSdMLPdrift}{1.278} % results/camera_ready/camera_ready_metrics.json:models.mlp_drift.metrics.chi0.sd
\newcommand{\tabChiSdMLPflow}{0.313} % results/camera_ready/camera_ready_metrics.json:models.mlp_flow.metrics.chi0.sd
\newcommand{\tabChiSdNODE}{0.045} % results/camera_ready/camera_ready_metrics.json:models.node.metrics.chi0.sd
\newcommand{\tabChiSdNSDE}{0.039} % results/camera_ready/camera_ready_metrics.json:models.nsde.metrics.chi0.sd
\newcommand{\tabChiSdRC}{0.032} % results/camera_ready/camera_ready_metrics.json:models.rc_flowmap.metrics.chi0.sd
\newcommand{\tabChiSdSindyFull}{0.029} % results/camera_ready/camera_ready_metrics.json:diagnostic.sindy_full.metrics.chi0.sd
\newcommand{\tabChiSdTransformer}{0.033} % results/camera_ready/camera_ready_metrics.json:exploratory.transformer.metrics.chi0.sd
\newcommand{\tabChiSdTrue}{0.034} % results/camera_ready/camera_ready_metrics.json:models.true.metrics.chi0.sd
\newcommand{\tabFloorAcf}{0.089} % results/camera_ready/camera_ready_metrics.json:criteria.floors.acf_distance.limits
\newcommand{\tabFloorFPlushi}{0.517} % results/camera_ready/camera_ready_metrics.json:criteria.floors.f_plus.limits[1]
\newcommand{\tabFloorFPluslo}{0.487} % results/camera_ready/camera_ready_metrics.json:criteria.floors.f_plus.limits[0]
\newcommand{\tabFloorSwitchhi}{1.035} % results/camera_ready/camera_ready_metrics.json:criteria.floors.switch_ratio.limits[1]
\newcommand{\tabFloorSwitchlo}{0.985} % results/camera_ready/camera_ready_metrics.json:criteria.floors.switch_ratio.limits[0]
\newcommand{\tabFloorWOne}{0.284} % results/camera_ready/camera_ready_metrics.json:criteria.floors.w1_x.limits
\newcommand{\tabOvMLPdriftFive}{0.290} % results/camera_ready/camera_ready_metrics.json:models.mlp_drift.final.per_mode_full.modes[4].overlap.median (values table)
\newcommand{\tabOvMLPdriftFour}{0.374} % results/camera_ready/camera_ready_metrics.json:models.mlp_drift.final.per_mode_full.modes[3].overlap.median (values table)
\newcommand{\tabOvMLPdriftOne}{0.607} % results/camera_ready/camera_ready_metrics.json:models.mlp_drift.final.per_mode_full.modes[0].overlap.median (values table)
\newcommand{\tabOvMLPdriftSix}{0.115} % results/camera_ready/camera_ready_metrics.json:models.mlp_drift.final.per_mode_full.modes[5].overlap.median (values table)
\newcommand{\tabOvMLPdriftThree}{0.319} % results/camera_ready/camera_ready_metrics.json:models.mlp_drift.final.per_mode_full.modes[2].overlap.median (values table)
\newcommand{\tabOvMLPdriftTwo}{0.043} % results/camera_ready/camera_ready_metrics.json:models.mlp_drift.final.per_mode_full.modes[1].overlap.median (values table)
\newcommand{\tabOvMLPflowFive}{0.691} % results/camera_ready/camera_ready_metrics.json:models.mlp_flow.final.per_mode_full.modes[4].overlap.median (values table)
\newcommand{\tabOvMLPflowFour}{0.751} % results/camera_ready/camera_ready_metrics.json:models.mlp_flow.final.per_mode_full.modes[3].overlap.median (values table)
\newcommand{\tabOvMLPflowOne}{0.963} % results/camera_ready/camera_ready_metrics.json:models.mlp_flow.final.per_mode_full.modes[0].overlap.median (values table)
\newcommand{\tabOvMLPflowSix}{0.245} % results/camera_ready/camera_ready_metrics.json:models.mlp_flow.final.per_mode_full.modes[5].overlap.median (values table)
\newcommand{\tabOvMLPflowThree}{0.755} % results/camera_ready/camera_ready_metrics.json:models.mlp_flow.final.per_mode_full.modes[2].overlap.median (values table)
\newcommand{\tabOvMLPflowTwo}{0.820} % results/camera_ready/camera_ready_metrics.json:models.mlp_flow.final.per_mode_full.modes[1].overlap.median (values table)
\newcommand{\tabOvNODEFive}{0.935} % results/camera_ready/camera_ready_metrics.json:models.node.final.per_mode_full.modes[4].overlap.median (values table)
\newcommand{\tabOvNODEFour}{0.952} % results/camera_ready/camera_ready_metrics.json:models.node.final.per_mode_full.modes[3].overlap.median (values table)
\newcommand{\tabOvNODEOne}{0.998} % results/camera_ready/camera_ready_metrics.json:models.node.final.per_mode_full.modes[0].overlap.median (values table)
\newcommand{\tabOvNODESix}{0.996} % results/camera_ready/camera_ready_metrics.json:models.node.final.per_mode_full.modes[5].overlap.median (values table)
\newcommand{\tabOvNODEThree}{0.962} % results/camera_ready/camera_ready_metrics.json:models.node.final.per_mode_full.modes[2].overlap.median (values table)
\newcommand{\tabOvNODETwo}{0.976} % results/camera_ready/camera_ready_metrics.json:models.node.final.per_mode_full.modes[1].overlap.median (values table)
\newcommand{\tabOvNSDEFive}{0.929} % results/camera_ready/camera_ready_metrics.json:models.nsde.final.per_mode_full.modes[4].overlap.median (values table)
\newcommand{\tabOvNSDEFour}{0.953} % results/camera_ready/camera_ready_metrics.json:models.nsde.final.per_mode_full.modes[3].overlap.median (values table)
\newcommand{\tabOvNSDEOne}{0.998} % results/camera_ready/camera_ready_metrics.json:models.nsde.final.per_mode_full.modes[0].overlap.median (values table)
\newcommand{\tabOvNSDESix}{0.996} % results/camera_ready/camera_ready_metrics.json:models.nsde.final.per_mode_full.modes[5].overlap.median (values table)
\newcommand{\tabOvNSDEThree}{0.960} % results/camera_ready/camera_ready_metrics.json:models.nsde.final.per_mode_full.modes[2].overlap.median (values table)
\newcommand{\tabOvNSDETwo}{0.980} % results/camera_ready/camera_ready_metrics.json:models.nsde.final.per_mode_full.modes[1].overlap.median (values table)
\newcommand{\tabOvRCFive}{0.856} % results/camera_ready/camera_ready_metrics.json:models.rc_flowmap.final.per_mode_full.modes[4].overlap.median (values table)
\newcommand{\tabOvRCFour}{0.895} % results/camera_ready/camera_ready_metrics.json:models.rc_flowmap.final.per_mode_full.modes[3].overlap.median (values table)
\newcommand{\tabOvRCOne}{0.998} % results/camera_ready/camera_ready_metrics.json:models.rc_flowmap.final.per_mode_full.modes[0].overlap.median (values table)
\newcommand{\tabOvRCSix}{0.997} % results/camera_ready/camera_ready_metrics.json:models.rc_flowmap.final.per_mode_full.modes[5].overlap.median (values table)
\newcommand{\tabOvRCThree}{0.888} % results/camera_ready/camera_ready_metrics.json:models.rc_flowmap.final.per_mode_full.modes[2].overlap.median (values table)
\newcommand{\tabOvRCTwo}{0.923} % results/camera_ready/camera_ready_metrics.json:models.rc_flowmap.final.per_mode_full.modes[1].overlap.median (values table)
\newcommand{\tabOvSindyFullFive}{0.993} % results/camera_ready/camera_ready_metrics.json:diagnostic.sindy_full.final.per_mode_full.modes[4].overlap.median (values table)
\newcommand{\tabOvSindyFullFour}{0.990} % results/camera_ready/camera_ready_metrics.json:diagnostic.sindy_full.final.per_mode_full.modes[3].overlap.median (values table)
\newcommand{\tabOvSindyFullOne}{0.999} % results/camera_ready/camera_ready_metrics.json:diagnostic.sindy_full.final.per_mode_full.modes[0].overlap.median (values table)
\newcommand{\tabOvSindyFullSix}{0.998} % results/camera_ready/camera_ready_metrics.json:diagnostic.sindy_full.final.per_mode_full.modes[5].overlap.median (values table)
\newcommand{\tabOvSindyFullThree}{0.986} % results/camera_ready/camera_ready_metrics.json:diagnostic.sindy_full.final.per_mode_full.modes[2].overlap.median (values table)
\newcommand{\tabOvSindyFullTwo}{0.994} % results/camera_ready/camera_ready_metrics.json:diagnostic.sindy_full.final.per_mode_full.modes[1].overlap.median (values table)
\newcommand{\tabOvTrueFive}{0.991} % results/camera_ready/camera_ready_metrics.json:models.true.final.per_mode_full.modes[4].overlap.median (values table)
\newcommand{\tabOvTrueFour}{0.990} % results/camera_ready/camera_ready_metrics.json:models.true.final.per_mode_full.modes[3].overlap.median (values table)
\newcommand{\tabOvTrueOne}{0.999} % results/camera_ready/camera_ready_metrics.json:models.true.final.per_mode_full.modes[0].overlap.median (values table)
\newcommand{\tabOvTrueSix}{0.998} % results/camera_ready/camera_ready_metrics.json:models.true.final.per_mode_full.modes[5].overlap.median (values table)
\newcommand{\tabOvTrueThree}{0.991} % results/camera_ready/camera_ready_metrics.json:models.true.final.per_mode_full.modes[2].overlap.median (values table)
\newcommand{\tabOvTrueTwo}{0.992} % results/camera_ready/camera_ready_metrics.json:models.true.final.per_mode_full.modes[1].overlap.median (values table)
\newcommand{\tabThrFive}{0.960} % results/camera_ready/camera_ready_metrics.json:models.mlp_drift.final.per_mode_full.modes[4].threshold (values table)
\newcommand{\tabThrFour}{0.968} % results/camera_ready/camera_ready_metrics.json:models.mlp_drift.final.per_mode_full.modes[3].threshold (values table)
\newcommand{\tabThrOne}{0.990} % results/camera_ready/camera_ready_metrics.json:models.mlp_drift.final.per_mode_full.modes[0].threshold (values table)
\newcommand{\tabThrSix}{0.990} % results/camera_ready/camera_ready_metrics.json:models.mlp_drift.final.per_mode_full.modes[5].threshold (values table)
\newcommand{\tabThrThree}{0.922} % results/camera_ready/camera_ready_metrics.json:models.mlp_drift.final.per_mode_full.modes[2].threshold (values table)
\newcommand{\tabThrTwo}{0.945} % results/camera_ready/camera_ready_metrics.json:models.mlp_drift.final.per_mode_full.modes[1].threshold (values table)
\newcommand{\tabVAcfMLPdrift}{3.724} % results/camera_ready/camera_ready_metrics.json:models.mlp_drift.final.invariant.acf_distance.median (values table, 3 decimals)
\newcommand{\tabVAcfMLPflow}{0.735} % results/camera_ready/camera_ready_metrics.json:models.mlp_flow.final.invariant.acf_distance.median (values table, 3 decimals)
\newcommand{\tabVAcfNODE}{0.093} % results/camera_ready/camera_ready_metrics.json:models.node.final.invariant.acf_distance.median (values table, 3 decimals)
\newcommand{\tabVAcfNSDE}{0.102} % results/camera_ready/camera_ready_metrics.json:models.nsde.final.invariant.acf_distance.median (values table, 3 decimals)
\newcommand{\tabVAcfRC}{0.070} % results/camera_ready/camera_ready_metrics.json:models.rc_flowmap.final.invariant.acf_distance.median (values table, 3 decimals)
\newcommand{\tabVAcfSindyFull}{0.038} % results/camera_ready/camera_ready_metrics.json:diagnostic.sindy_full.final.invariant.acf_distance.median (values table, 3 decimals)
\newcommand{\tabVAcfTrue}{0.038} % results/camera_ready/camera_ready_metrics.json:models.true.final.invariant.acf_distance.median (values table, 3 decimals)
\newcommand{\tabVFPlusMLPdrift}{0.477} % results/camera_ready/camera_ready_metrics.json:models.mlp_drift.final.invariant.f_plus.median (values table, 3 decimals)
\newcommand{\tabVFPlusMLPflow}{0.474} % results/camera_ready/camera_ready_metrics.json:models.mlp_flow.final.invariant.f_plus.median (values table, 3 decimals)
\newcommand{\tabVFPlusNODE}{0.504} % results/camera_ready/camera_ready_metrics.json:models.node.final.invariant.f_plus.median (values table, 3 decimals)
\newcommand{\tabVFPlusNSDE}{0.502} % results/camera_ready/camera_ready_metrics.json:models.nsde.final.invariant.f_plus.median (values table, 3 decimals)
\newcommand{\tabVFPlusRC}{0.499} % results/camera_ready/camera_ready_metrics.json:models.rc_flowmap.final.invariant.f_plus.median (values table, 3 decimals)
\newcommand{\tabVFPlusSindyFull}{0.494} % results/camera_ready/camera_ready_metrics.json:diagnostic.sindy_full.final.invariant.f_plus.median (values table, 3 decimals)
\newcommand{\tabVFPlusTrue}{0.500} % results/camera_ready/camera_ready_metrics.json:models.true.final.invariant.f_plus.median (values table, 3 decimals)
\newcommand{\tabVSwitchMLPdrift}{0.008} % results/camera_ready/camera_ready_metrics.json:models.mlp_drift.final.invariant.switch_ratio.median (values table, 3 decimals)
\newcommand{\tabVSwitchMLPflow}{0.725} % results/camera_ready/camera_ready_metrics.json:models.mlp_flow.final.invariant.switch_ratio.median (values table, 3 decimals)
\newcommand{\tabVSwitchNODE}{1.049} % results/camera_ready/camera_ready_metrics.json:models.node.final.invariant.switch_ratio.median (values table, 3 decimals)
\newcommand{\tabVSwitchNSDE}{1.050} % results/camera_ready/camera_ready_metrics.json:models.nsde.final.invariant.switch_ratio.median (values table, 3 decimals)
\newcommand{\tabVSwitchRC}{1.037} % results/camera_ready/camera_ready_metrics.json:models.rc_flowmap.final.invariant.switch_ratio.median (values table, 3 decimals)
\newcommand{\tabVSwitchSindyFull}{1.013} % results/camera_ready/camera_ready_metrics.json:diagnostic.sindy_full.final.invariant.switch_ratio.median (values table, 3 decimals)
\newcommand{\tabVSwitchTrue}{1.011} % results/camera_ready/camera_ready_metrics.json:models.true.final.invariant.switch_ratio.median (values table, 3 decimals)
\newcommand{\tabVWOneMLPdrift}{5.599} % results/camera_ready/camera_ready_metrics.json:models.mlp_drift.final.invariant.w1_x.median (values table, 3 decimals)
\newcommand{\tabVWOneMLPflow}{2.193} % results/camera_ready/camera_ready_metrics.json:models.mlp_flow.final.invariant.w1_x.median (values table, 3 decimals)
\newcommand{\tabVWOneNODE}{0.337} % results/camera_ready/camera_ready_metrics.json:models.node.final.invariant.w1_x.median (values table, 3 decimals)
\newcommand{\tabVWOneNSDE}{0.362} % results/camera_ready/camera_ready_metrics.json:models.nsde.final.invariant.w1_x.median (values table, 3 decimals)
\newcommand{\tabVWOneRC}{0.097} % results/camera_ready/camera_ready_metrics.json:models.rc_flowmap.final.invariant.w1_x.median (values table, 3 decimals)
\newcommand{\tabVWOneSindyFull}{0.092} % results/camera_ready/camera_ready_metrics.json:diagnostic.sindy_full.final.invariant.w1_x.median (values table, 3 decimals)
\newcommand{\tabVWOneTrue}{0.076} % results/camera_ready/camera_ready_metrics.json:models.true.final.invariant.w1_x.median (values table, 3 decimals)
\newcommand{\thrFive}{0.9601} % results/camera_ready/camera_ready_metrics.json:models.mlp_drift.final.per_mode_full.modes[4].threshold
\newcommand{\thrFour}{0.9681} % results/camera_ready/camera_ready_metrics.json:models.mlp_drift.final.per_mode_full.modes[3].threshold
\newcommand{\thrOne}{0.9900} % results/camera_ready/camera_ready_metrics.json:models.mlp_drift.final.per_mode_full.modes[0].threshold
\newcommand{\thrSix}{0.9900} % results/camera_ready/camera_ready_metrics.json:models.mlp_drift.final.per_mode_full.modes[5].threshold
\newcommand{\thrThree}{0.9216} % results/camera_ready/camera_ready_metrics.json:models.mlp_drift.final.per_mode_full.modes[2].threshold
\newcommand{\thrTwo}{0.9450} % results/camera_ready/camera_ready_metrics.json:models.mlp_drift.final.per_mode_full.modes[1].threshold
\newcommand{\thrUncFive}{0.9601} % results/camera_ready/camera_ready_metrics.json:models.mlp_drift.final.per_mode_full.modes[4].threshold_uncapped
\newcommand{\thrUncFour}{0.9681} % results/camera_ready/camera_ready_metrics.json:models.mlp_drift.final.per_mode_full.modes[3].threshold_uncapped
\newcommand{\thrUncOne}{0.9984} % results/camera_ready/camera_ready_metrics.json:models.mlp_drift.final.per_mode_full.modes[0].threshold_uncapped
\newcommand{\thrUncSix}{0.9960} % results/camera_ready/camera_ready_metrics.json:models.mlp_drift.final.per_mode_full.modes[5].threshold_uncapped
\newcommand{\thrUncThree}{0.9216} % results/camera_ready/camera_ready_metrics.json:models.mlp_drift.final.per_mode_full.modes[2].threshold_uncapped
\newcommand{\thrUncTwo}{0.9450} % results/camera_ready/camera_ready_metrics.json:models.mlp_drift.final.per_mode_full.modes[1].threshold_uncapped
\newcommand{\trfContext}{16} % results/camera_ready/camera_ready_metrics.json:exploratory.transformer.appendix.recipe.context
\newcommand{\trfDim}{64} % results/camera_ready/camera_ready_metrics.json:exploratory.transformer.appendix.recipe.d_model
\newcommand{\trfHeads}{4} % results/camera_ready/camera_ready_metrics.json:exploratory.transformer.appendix.recipe.n_heads
\newcommand{\trfLayers}{2} % results/camera_ready/camera_ready_metrics.json:exploratory.transformer.appendix.recipe.n_layers
\newcommand{\trfParams}{68{,}547} % results/camera_ready/camera_ready_metrics.json:exploratory.transformer.appendix.recipe.n_params (thousands separator)
\newcommand{\trfRolloutSteps}{8} % results/camera_ready/camera_ready_metrics.json:exploratory.transformer.appendix.recipe.rollout_steps
\newcommand{\pctn}[1]{\fpeval{round(100*(#1))}}
\newcommand{\pct}[1]{\pctn{#1}\%}

\begin{document}

\maketitle

\begin{abstract}
Machine-learning emulators of chaotic and stochastic systems are usually validated on forecast skill and long-run statistics. 
Neither certifies that an emulator responds correctly to forcing, the property that projection and attribution studies rely on. 
Linear response theory makes this testable: the forced response follows from unperturbed correlations through a generalized fluctuation--dissipation relation, and decomposes over the stochastic Ruelle--Pollicott resonances of the Koopman generator. 
Building on the Koopmanism Response framework, we turn this into a calibrated, mode-resolved test for learned surrogates: each surrogate rollout passes or fails each check, and failure rates are compared with those of independent realizations of the true system.
On stochastic Lorenz-63, a three-variable toy model, we evaluate SINDy, an MLP, a reservoir computer, a neural ODE and a neural SDE with learned diffusion, over up to \NODEn{} rollouts each. A sparse-regression model with the correct library passes every check at rates consistent with the true system.
Invariant-statistics fidelity and response fidelity dissociate in both directions: a quarter of reservoir-computer rollouts pass every invariant-statistics check and match the static susceptibility $\chi(0)$, yet misrepresent the slow relaxation modes, while the neural ODE and SDE rarely meet the invariant-statistics floor but recover those modes in three quarters of rollouts.
As expected of a time-integrated quantity dominated here by fast relaxation, $\chi(0)$ does not separate these cases.
For a fixed network, the training formulation (one-step drift, flow map, or multi-step through the integrator) decides which of these properties it gets right.
\end{abstract}

\section{Introduction \& Motivation}
Machine learning is now a standard tool for modelling dynamical systems, from sparse regression of governing equations~\cite{brunton_discovering_2016} to reservoir computing~\cite{pathak_model-free_2018} and neural differential equations~\cite{chen_neural_2018}. Data-driven weather models such as GraphCast~\cite{lam_learning_2023}, FourCastNet~\cite{pathak_fourcastnet_2022}, Pangu-Weather~\cite{bi_accurate_2023} and Aurora~\cite{bodnar_foundation_2025} rival physics-based forecasts at a fraction of the cost.
Their validation still rests mainly on forecast skill and long-run statistics, increasingly through benchmarks such as WeatherBench~2~\cite{rasp_weatherbench_2024}; case studies and dynamical tests of their physical behaviour are emerging~\cite{charlton-perez_ai_2024,hakim_dynamical_2024,brocker_verification_2026}, but are not yet standard.
Both test fidelity to the unperturbed distribution, not how the statistics change under forcing, which is what projection, attribution and the study of tipping points require~\cite{ghil_physics_2020}: a model can reproduce the invariant measure closely and still misrepresent how it deforms under forcing~\cite{falasca_probing_2025}.

Response theory makes this testable. For a large class of chaotic and stochastic systems, the linear response is a correlation function of the unperturbed dynamics, in the spirit of the generalized fluctuation--dissipation theorem~\cite{marconi_fluctuationdissipation_2008,ruelle_review_2009}, and it decomposes over the stochastic Ruelle--Pollicott (RP) resonances of the Koopman generator and their eigenfunctions~\cite{chekroun_ruellepollicott_2020,lucarini_theoretical_2023}. The Koopmanism Response framework~\cite{lucarini_general_2026,zagli_bridging_2026} computes these objects from a single unperturbed trajectory with extended dynamic mode decomposition (EDMD), with ResDMD residuals guarding against spurious modes~\cite{williams_datadriven_2015,colbrook_residual_2023}, and yields a mode-resolved decomposition of the forced Green's function~\cite{lucarini_interpretable_2026}. ``Did the model learn the physics?'' then becomes a falsifiable question: does a surrogate trained only on unperturbed data carry the correct resonances, eigenfunctions and static susceptibility?

Stochastic Lorenz-63 is a three-variable toy model of convection with no spatial structure, and our surrogates are not operational climate emulators: it is a controlled testbed with a known true response, used to map how network design and training formulation shape response fidelity, complementing forced-response studies of neural climate emulators~\cite{falasca_probing_2025}.
We contribute (i)~a calibrated evaluation protocol for learned stochastic simulators, in which invariant statistics, the static susceptibility and a mode-resolved eigenfunction test are each applied to every rollout and compared with the false-alarm rate of independent true realizations;
(ii)~a formulation ladder (one-step drift regression, flow map, multi-step neural ODE, neural SDE with learned diffusion) sharing architecture and data, which isolates the effects of the training target, the multi-step loss and learned noise, complemented by an exploratory causal transformer, the architecture family of current weather emulators;
and (iii)~evidence that invariant-statistics fidelity and response fidelity dissociate in both directions, which a scalar response summary such as $\chi(0)$ cannot detect here but the mode-resolved test does.

\section{Operator theory}

\label{app:operators}

This section collects the operator-theoretic objects that the response diagnostic of Sec.~\ref{sec:method} is built on, following the Koopmanism Response framework \cite{zagli_bridging_2026}.
Standard references for the underlying stochastic-analysis facts are \cite{pavliotis_stochastic_2014} for the Kolmogorov operator and \cite{chekroun_ruellepollicott_2020} for stochastic Ruelle-Pollicott resonances.

\subsection{Stochastic dynamics and the invariant measure}

We consider an autonomous It\^o SDE at statistical stationarity,
\begin{equation}
\label{eq:sde}
d\mathbf{x} = \mathbf{F}(\mathbf{x})\,dt + \varsigma\,d\mathbf{W}_t,
\qquad \mathbf{x}\in\mathcal{M}\subseteq\mathbb{R}^d,
\end{equation}
with (generally nonlinear) drift $\mathbf{F}$, isotropic additive noise of intensity $\varsigma>0$, and a $d$-dimensional Wiener process $\mathbf{W}_t$.
The one-time statistics are carried by the probability density $\rho(\mathbf{x},t)$, which solves the Fokker-Planck (forward Kolmogorov) equation
\begin{equation}
\label{eq:fokkerplanck}
\partial_t \rho = \mathcal{L}_0\,\rho
:= -\nabla\!\cdot(\mathbf{F}\rho) + \tfrac{\varsigma^2}{2}\,\nabla^2\rho .
\end{equation}
Under hypoellipticity of $\mathcal{L}_0$ and mild growth conditions on $\mathbf{F}$, the process \eqref{eq:sde} admits a unique, ergodic invariant measure $d\mu_0(\mathbf{x})=\rho_0(\mathbf{x})\,d\mathbf{x}$ that is absolutely continuous with respect to Lebesgue and has a smooth density \cite{pavliotis_stochastic_2014}.
Stationary expectations are $\langle f\rangle_{\mu_0}=\int_{\mathcal{M}} f\,d\mu_0$, and the two pillars of standard surrogate validation (attractor geometry, long-run marginals) are exactly statements about $\mu_0$.

\subsection{The Markov semigroup and its generator}

The statistics of observables evolve under the Markov semigroup $\mathbf{P}=(P_t)_{t\ge 0}$,
\begin{equation}
(P_t f)(\mathbf{x}) = \int_{\mathcal{M}} f(\mathbf{y})\,
p(\mathbf{y},t\,|\,\mathbf{x})\,d\mathbf{y}
= \mathbb{E}\!\left[\,f(\mathbf{x}(t))\mid \mathbf{x}(0)=\mathbf{x}\,\right],
\end{equation}
where $p(\mathbf{y},t\,|\,\mathbf{x})$ is the transition density.
On $L^2_{\mu_0}=\{f:\int|f|^2\,d\mu_0<\infty\}$ this is a strongly continuous semigroup, so it can be written in exponential form
\begin{equation}
\label{eq:semigroup}
P_t = e^{t\mathcal{K}_0},
\end{equation}
with (generally unbounded) generator $\mathcal{K}_0$, the Kolmogorov (backward) operator.
It is the second-order differential operator
\begin{equation}
\label{eq:generator}
\mathcal{K}_0 = \mathbf{F}(\mathbf{x})\cdot\nabla
+ \tfrac{\varsigma^2}{2}\,\nabla^2,
\end{equation}
and it advances expectations of observables, $\partial_t\,\overline f =\mathcal{K}_0\overline f$ with $\overline f(\mathbf{x},t)=(P_tf)(\mathbf{x})$.
In the deterministic limit $\varsigma\to 0$, $\mathcal{K}_0=\mathbf{F}\cdot\nabla$ reduces to the usual Koopman generator.
We retain the operator name ``Koopman generator'' for $\mathcal{K}_0$ throughout, following the convention of \cite{zagli_bridging_2026}.

\subsection{Why the noise: deterministic vs.\ stochastic spectra}

The regularization $\varsigma>0$ in \eqref{eq:sde} is operator-theoretic, not cosmetic, and it is what makes the response machinery well posed.
The deterministic Lorenz-63 flow carries a singular SRB measure, and its Koopman operator is unitary on $L^2_{\mu_0}$: the dynamics has already settled onto the attractor, so no decay toward it can be observed and the spectrum sits on the imaginary axis (equivalently, on the unit circle for the time-$\tau$ map).
Adding non-degenerate noise makes $\mathcal{L}_0$ hypoelliptic, yielding a unique smooth $\rho_0$ and a quasi-compact semigroup, so that a discrete set of relaxation modes emerges and the data-driven response reconstruction becomes well posed \cite{chekroun_ruellepollicott_2020,zagli_bridging_2026}.

\subsection{Ruelle-Pollicott resonances and the spectral decomposition}

Assume $\mathbf{P}$ is quasi-compact, i.e.\ the eigenvalues of $\mathcal{K}_0$ with real part above the essential growth bound $r_{\mathrm{ess}}<0$ are discrete with finite multiplicity.
By the spectral mapping theorem, if $\lambda_j$ is an eigenvalue of $\mathcal{K}_0$ with eigenfunction $\varphi_j$, then $e^{\lambda_j t}$ is an eigenvalue of $P_t$ with the same eigenfunction.
These discrete eigenvalues $\{\lambda_j\}_{j\ge 0}$, ordered as $\lambda_0=0>\operatorname{Re}\lambda_1\ge\operatorname{Re}\lambda_2\ge\cdots$, are the (stochastic) Ruelle-Pollicott resonances, and they give a physically meaningful decomposition of the semigroup,
\begin{equation}
\label{eq:spectral}
e^{t\mathcal{K}_0} = \sum_{j\ge 0} e^{\lambda_j t}\,\Pi_j + \mathcal{R}(t),
\end{equation}
where $\Pi_j$ projects onto the eigenspace of $\lambda_j$ and $\mathcal{R}(t)$ collects the essential-spectrum contribution, whose norm decays and becomes negligible as $t\to\infty$ under quasi-compactness.
The resonances are intrinsic relaxation rates: a nonzero imaginary part produces oscillatory correlations, the signature of nonequilibrium dynamics, while the spectral gap $\gamma=-\operatorname{Re}\lambda_1$ controls the slowest decay.
An accumulation of resonances toward $\lambda_0=0$ signals diverging sensitivity (critical slowing down).
This is the structure a learned surrogate must reproduce, and it is precisely what neither forecast skill nor invariant-measure agreement constrains.

For zero-mean observables $f,g\in L^2_{\mu_0}$, the time-lagged correlation function inherits \eqref{eq:spectral},
\begin{equation}
\label{eq:corr}
C_{f,g}(t) = \int_{\mathcal{M}}\!\big(e^{t\mathcal{K}_0}f\big)(\mathbf{x})\,
g^*(\mathbf{x})\,d\mu_0(\mathbf{x})
= \sum_{j\ge 1}\alpha_j(f,g)\,e^{\lambda_j t},
\end{equation}
a sum of exponentially decaying (and possibly oscillating) modes; the rates $\lambda_j$ are intrinsic to the system, whereas the coefficients $\alpha_j(f,g)$ encode the projection of the chosen observables onto the Kolmogorov modes.

\subsection{Extended Dynamic Mode Decomposition (EDMD)}

We estimate the spectral objects in \eqref{eq:spectral}-\eqref{eq:corr} from a single long unperturbed trajectory via EDMD \cite{williams_datadriven_2015}.
Given snapshot pairs $\mathbb{D}=\{(\mathbf{x}_i,\mathbf{y}_i=\mathbf{x}_{i+1})\}_{i=1}^{M}$ sampled at step $\tau$ and a dictionary $\Psi(\mathbf{x})=(\psi_1(\mathbf{x}),\dots,\psi_N(\mathbf{x}))$ spanning $\mathcal{F}_N\subseteq L^2_{\mu_0}$, form the Gram and one-step (action) matrices
\begin{equation}
\label{eq:edmd-matrices}
\mathbf{G}=\tfrac1M\sum_{i}\Psi^*(\mathbf{x}_i)\Psi(\mathbf{x}_i),
\qquad
\mathbf{A}=\tfrac1M\sum_{i}\Psi^*(\mathbf{x}_i)\Psi(\mathbf{y}_i),
\qquad
\mathbf{K}=\mathbf{G}^{+}\mathbf{A},
\end{equation}
so that $\mathbf{G}_{ab}=\langle\psi_a,\psi_b\rangle_{\mu_0}$ and $\mathbf{A}_{ab}=\langle\psi_a,\mathcal{K}_0\psi_b\rangle_{\mu_0}$ (here $\mathbf{G}^{+}$ is a pseudoinverse).
The matrix $\mathbf{K}$ is the Galerkin approximation of the semigroup on $\mathcal{F}_N$.
Writing its eigenpairs as $\mathbf{K}\boldsymbol{\xi}_k=\nu_k\boldsymbol{\xi}_k$, the Kolmogorov eigenvalues and eigenfunctions follow from the spectral mapping theorem,
\begin{equation}
\label{eq:edmd-spec}
\lambda_k=\frac{\ln\nu_k}{\tau},
\qquad
\varphi_k(\mathbf{x})=\Psi(\mathbf{x})\,\boldsymbol{\xi}_k,
\end{equation}
and the Koopman-mode coefficients of any observable $f=\Psi\mathbf{h}\in\mathcal{F}_N$ are $\mathbf{f}=\mathbf{W}'\mathbf{h}$, with $\mathbf{W}$ the matrix of left eigenvectors of $\mathbf{K}$.
Correlations then reconstruct as $C_{f,g}(t)=\sum_k \tilde f_k\,e^{\lambda_k t}$, where the projection coefficients are read off from the EDMD Gram matrix, $\langle\varphi_l,\varphi_k\rangle_{\mu_0}=(\boldsymbol{\Xi}'\mathbf{G}\boldsymbol{\Xi})_{lk}$ with $\boldsymbol{\Xi}=[\boldsymbol{\xi}_1\cdots\boldsymbol{\xi}_N]$.
As $N\to\infty$ with $\{\psi_i\}$ a basis of $L^2_{\mu_0}$, EDMD converges to $\mathcal{K}_0$ in the strong operator topology \cite{korda_convergence_2018}, which justifies the reconstruction of correlation and response operators from data.

In the experiments we take $\Psi$ to be a three-dimensional Fourier dictionary of size $(2K_{\max}+1)^3$ on a bounded rescaling of the attractor domain, and regularize $\mathbf{G}^{+}$ by a truncated SVD (TSVD).
% A Fourier basis is the natural, closed-under-differentiation choice on a torus; on $\mathbb{R}^3$ with $\rho_0$ sharply peaked on the two-lobe attractor it is less well matched, and delay-embedding (Takens) dictionaries adapted to the observed geometry, together with resolvent-based generalized eigenfunctions for the continuous part of the spectrum, are the principled remedy \cite{colbrook_rigged_2025} - we address this in the discussion.

\subsection{Certifying the spectrum: ResDMD residuals}

\label{app:resdmd}

Galerkin truncation can produce spurious (spectrally polluted) modes, so a finite EDMD spectrum must be verified rather than trusted.
ResDMD supplies a data-driven residual for each candidate eigenpair using one additional matrix computed from the same snapshots \cite{colbrook_residual_2023},
\begin{equation}
\label{eq:Lmatrix}
\mathbf{L}=\tfrac1M\sum_i \Psi^*(\mathbf{y}_i)\Psi(\mathbf{y}_i),
\qquad
\mathbf{L}_{ab}=\langle\mathcal{K}_0\psi_a,\mathcal{K}_0\psi_b\rangle_{\mu_0}.
\end{equation}
For a candidate pair $(\lambda,\boldsymbol{\xi})$ with $g=\Psi\boldsymbol{\xi}$, the relative residual is
\begin{equation}
\label{eq:resdmd}
\mathrm{res}(\lambda,\boldsymbol{\xi})
=\frac{\big\|(\mathcal{K}_0-\lambda)\,g\big\|_{\mu_0}}
       {\|g\|_{\mu_0}}
=\left(
\frac{\boldsymbol{\xi}^*\big(\mathbf{L}-\lambda\mathbf{A}^*
-\bar\lambda\mathbf{A}+|\lambda|^2\mathbf{G}\big)\boldsymbol{\xi}}
     {\boldsymbol{\xi}^*\mathbf{G}\boldsymbol{\xi}}
\right)^{1/2},
\end{equation}
which is exact in the $L^2_{\mu_0}$ norm and requires no knowledge of the true spectrum.
Residuals of the leading candidate pairs are stored for every scored rollout and used as a diagnostic, not as an acceptance gate.
This is what lets us treat the reference resonances as ground truth and, symmetrically, what makes the mode-resolved surrogate test meaningful. 
We compare a surrogate's leading eigenfunctions against the reference not by eigenvalue proximity alone but by the $\mu_0$-weighted subspace overlap
\begin{equation}
\label{eq:overlap}
O(\varphi,\tilde\varphi)
=\frac{|\langle\varphi,\tilde\varphi\rangle_{\mu_0}|}
       {\|\varphi\|_{\mu_0}\,\|\tilde\varphi\|_{\mu_0}},
\end{equation}
generalized to the principal angles between eigenspaces for near-degenerate clusters (e.g.\ the near-degenerate $\refFreqFour{}/\refFreqFive{}$ pair).
This overlap is the discriminating quantity in Table~\ref{tab:main}: it is sensitive to eigenfunction structure, which the eigenvalue-only self-consistency check $C_{\varphi_k}(t)\propto e^{\lambda_k t}$ and the invariant statistics both miss.

\section{Response theory}

\label{app:response}

We now perturb the stationary dynamics and recall how linear response theory fixes the forced behaviour entirely from unperturbed statistics, through a generalized fluctuation-dissipation relation and the same resonances $\{\lambda_j\}$ of Sec.~\ref{app:operators}.
This is the content that turns ``did the surrogate learn the physics?'' into a falsifiable spectral question.
The formulation follows \cite{zagli_bridging_2026}.

\subsection{Perturbation and the perturbed Fokker-Planck equation}

Perturb the drift by a small, time-modulated phase-space pattern, $\mathbf{F}\to\mathbf{F}+\varepsilon\,\mathbf{X}(\mathbf{x})\,T(t)$, with $\varepsilon\ll1$, time profile $T(t)$, and spatial profile $\mathbf{X}$.
The density then obeys
\begin{equation}
\label{eq:pert-fp}
\partial_t\rho=\big(\mathcal{L}_0+\varepsilon\,T(t)\,\mathcal{L}_p\big)\rho,
\qquad
\rho(\mathbf{x},0)=\rho_0(\mathbf{x}),
\qquad
\mathcal{L}_p\,\cdot=-\nabla\!\cdot\!\big(\mathbf{X}(\mathbf{x})\,\cdot\big),
\end{equation}
where the initial condition encodes that the forcing switches on after the system has reached stationarity.
Because the drive is nonautonomous, no invariant measure of \eqref{eq:pert-fp} exists in general; we work to first order in $\varepsilon$.

\subsection{Linear response and the generalized FDT}

To first order the expectation of any observable $f$ obeys
\begin{equation}
\label{eq:lrt}
\langle f\rangle_{\mu}(t)
=\langle f\rangle_{\mu_0}+\varepsilon\,\big(G_f\star T\big)(t)+o(\varepsilon),
\end{equation}
with the causal Green's function given by the generalized fluctuation-dissipation theorem
\begin{equation}
\label{eq:green}
G_f(t)=\Theta(t)\int_{\mathcal{M}}
\frac{\mathcal{L}_p\rho_0}{\rho_0}\,
\big(e^{t\mathcal{K}_0}f\big)\,d\mu_0
=\Theta(t)\,C_{f,\Gamma}(t),
\qquad
\Gamma(\mathbf{x})=\frac{\mathcal{L}_p\rho_0}{\rho_0},
\end{equation}
where $\Theta$ is the Heaviside step.
Equation~\eqref{eq:green} is the crux: the forced response is a specific unperturbed correlation function between the observable $f$ and the conjugate observable $\Gamma$.
Away from equilibrium there is no first-principles form for $\rho_0$, so $\Gamma$ is intrinsic and not known a priori; nonetheless \eqref{eq:lrt}-\eqref{eq:green} let us estimate the response purely from unperturbed trajectory data.

\subsection{Mode-resolved response and the static susceptibility}

Substituting the spectral decomposition \eqref{eq:spectral} into \eqref{eq:green} gives the response as a sum over the same resonances,
\begin{equation}
\label{eq:green-spectral}
G_f(t)=\Theta(t)\sum_{j\ge1}\beta_j(f,\Gamma)\,e^{\lambda_j t},
\end{equation}
with coefficients $\beta_j(f,\Gamma)$ that project both $f$ and $\Gamma$ onto the Kolmogorov eigenfunctions.
This is a genuinely interpretable diagnostic: rather than a scalar error, the response decomposes into contributions from distinct modes of natural variability, so a surrogate can be interrogated mode by mode.

We summarize the response of the tracked observable ($f=x$ for Lorenz-63) by its zero-frequency limit, the static susceptibility,
\begin{equation}
\label{eq:chi}
\chi(0)=\int_0^\infty G_x(t)\,dt
=\frac{\delta\langle x\rangle_\infty}{\varepsilon},
\end{equation}
i.e.\ the steady-state shift per unit forcing under a switched-on profile $T(t)=\Theta(t)$.
From \eqref{eq:green-spectral}, $\chi(0)=\sum_k \beta_k(x,\Gamma)\!\int_0^\infty e^{\lambda_k t}\,dt=-\sum_k \beta_k(x,\Gamma)/\lambda_k$, convergent because $\operatorname{Re}\lambda_k<0$.
Evaluated from the sampled EDMD spectrum ($\nu_k=e^{\lambda_k\tau}$), the tail integral becomes a geometric series over the sampling grid $t=n\tau$,
\begin{equation}
\label{eq:chi-discrete}
\chi(0)\;\approx\;\tau\sum_{n\ge0}G_x(n\tau)
\;=\;\tau\sum_k \big[\beta_k\big]\sum_{n\ge0}\nu_k^{\,n}
\;=\;\tau\sum_k \frac{\alpha_k^{*}\,M\,\gamma_k}{1-\nu_k},
\end{equation}
which is the estimator used in Sec.~\ref{sec:method}.
Here $\alpha_k^{*}$ and $\gamma_k$ are the EDMD projection coefficients of the observable $x$ and of the conjugate observable $\Gamma$ onto the $k$-th Kolmogorov eigenfunction, and $M$ denotes the eigenfunction overlap (Gram) matrix $\langle\varphi_l,\varphi_k\rangle$ contracting the two, not the snapshot count $M$ of Sec.~\ref{app:operators}.
Because $\chi(0)$ is a mode sum weighted by $1/(1-\nu_k)$, it inherits nontrivial variance at finite trajectory length, which is why we calibrate a ground-truth band by re-evaluating the instrument on independent noise realizations of the true system.

\subsection{Data-driven estimation of \texorpdfstring{$\Gamma$}{Gamma}: Method 3 (integration by parts)}

\label{app:method3}

The one nontrivial ingredient is the conjugate observable $\Gamma=\mathcal{L}_p\rho_0/\rho_0$.
\cite{zagli_bridging_2026} give three routes to $G_f$.
Methods~1 and~2 build $\Gamma$ by first estimating $\rho_0$ (e.g.\ by binning) and then dividing by it; this is fragile whenever $\mu_0$ concentrates near a (near-)fractal attractor, exactly the Lorenz-63 situation.
Method~3 avoids density estimation altogether.
It seeks the $\mu_0$-weighted least-squares projection of $\Gamma$ onto the dictionary,
\begin{equation}
\label{eq:method3}
\mathbf{c}=\mathbf{G}^{+}\boldsymbol{\Delta},
\qquad
\Delta_i=\langle\psi_i,\Gamma\rangle_{\mu_0}
=\int_{\mathcal{M}}\!\big(\mathbf{X}\!\cdot\!\nabla\psi_i^{*}\big)\,\rho_0\,d\mathbf{x}
\;\approx\;\frac1T\int_0^T
\mathbf{X}\big(\mathbf{x}(t)\big)\!\cdot\!\nabla\psi_i^{*}\big(\mathbf{x}(t)\big)\,dt,
\end{equation}
where the middle equality is an integration by parts that moves the derivative off $\rho_0$ (boundary terms vanish for the smooth, decaying $\rho_0$), and the last step replaces the $\mu_0$-average by a time average along the same unperturbed trajectory used for EDMD.
The moments $\Delta_i$ therefore require no estimate of $\rho_0$ and no division by it; if the dictionary is closed under differentiation (as for Fourier), they use no information beyond EDMD itself, and the method is markedly data-efficient.
The coefficients of $\Gamma$ in the eigenfunction basis follow from $\mathbf{c}$ via the left eigenvectors, feeding directly into \eqref{eq:green-spectral}.

\subsection{Lorenz-63 specialization}

We force the $x$-equation only, $\mathbf{X}(\mathbf{x})=(1,0,0)^{\!\top}$, with $T(t)=\Theta(t)$, and track $f=x$.
The conjugate observable then simplifies to
\begin{equation}
\label{eq:gamma-l63}
\Gamma(\mathbf{x})=-\,\partial_x\ln\rho_0(\mathbf{x}),
\end{equation}
and the Method~3 moments \eqref{eq:method3} reduce to $\Delta_i\approx\frac1T\int_0^T\partial_x\psi_i^{*}(\mathbf{x}(t))\,dt$.
We force $x$ because it is the most strongly coupled variable, driving both $y$ and $z$ through the nonlinear terms, so a single-component forcing maximally probes the cross-variable response and is the cleanest first stress test of the framework.
The reference response is the instrument reconstruction on the held-out true realization (seed \seedRef{}).

\subsection{Why the response is a stringent diagnostic}

Equations \eqref{eq:green}-\eqref{eq:green-spectral} show that the entire forced response is determined by unperturbed objects: the resonances $\{\lambda_j\}$, their eigenfunctions $\varphi_j$, and the projection of the intrinsic conjugate observable $\Gamma$.
None of these is constrained by a short-horizon forecast loss or by agreement of the invariant measure.
A surrogate trained only on unperturbed data can therefore match the attractor, the wing-switching statistics, and even the scalar $\chi(0)$, while still carrying the wrong eigenfunction structure, a discrepancy that only the mode-resolved test of Sec.~\ref{app:resdmd}, applied through this response machinery, resolves.
Response theory thus supplies a rigorous, physically grounded criterion for whether a learned simulator respects dynamical structure rather than merely fitting the attractor.

\section{Method}
\label{sec:method}
\subsection{System and response instrument} For a dissipative stochastic system at stationarity, the linear response to a small forcing is fixed by the unperturbed dynamics through the RP resonances of the Koopman generator and a generalized fluctuation--dissipation relation~\cite{lucarini_general_2026,zagli_bridging_2026} (Secs.~\ref{app:operators},~\ref{app:response}). Following~\cite{lucarini_general_2026}, we study stochastic Lorenz-63, the It\^o SDE $dx = F(x)\,dt + \varsigma\,dW_t$ with $(\sigma,\rho,\beta) = (10, 28, 8/3)$ and isotropic noise $\varsigma = \paramS{}$, integrated by Euler--Maruyama at $\delta t = \paramDtRaw{}$. We force the $x$-equation and track $x$. From one unperturbed trajectory of $T = \paramT{}$ time units sampled at $\tau = \paramTau{}$, we build a Fourier-EDMD approximation of the generator (dictionary size $(2K_{\max}+1)^3 = \paramDict{}$, TSVD relative threshold $\paramTSVD{}$), read off the leading resonances, reconstruct the Green's function $G_x(t)$ with the Method-3 estimator~\cite{zagli_bridging_2026}, and summarize it by the static susceptibility
\begin{equation}
\chi(0)\;=\;\int_0^\infty G_x(t)\,dt\;=\;\frac{\delta\langle x\rangle_\infty}{\varepsilon}\;=\;\tau\sum_{k}\frac{\alpha_k^*\,M\,\gamma_k}{1-\nu_k},
\end{equation}
where $\delta\langle x\rangle_\infty$ is the long-term shift of $\langle x\rangle$, $\nu_k$ are the EDMD eigenvalues and the sum runs over all EDMD modes.

\subsection{Learned surrogates.} All models are trained on segments of the same unperturbed master trajectory and rolled out with additive noise (Sec.~\ref{app:models}). (i)~SINDy~\cite{brunton_discovering_2016,boninsegna_sparse_2018} with a quadratic library containing the true equations (positive control); 
(ii)~an MLP ($2\times64$, tanh) fitted to the drift at $\Delta t_{\mathrm{fit}} = \paramDtFitMLPdrift{}$ (MLP-drift) or to the next state (MLP-flow); 
(iii)~a reservoir computer (RC)~\cite{pathak_model-free_2018} fitted as a flow map;
(iv)~a neural ODE (NODE)~\cite{chen_neural_2018}, the same network trained through a differentiable RK4 step on multi-step windows; 
(v)~a neural SDE (NSDE)~\cite{kidger_neural_2021}, the NODE drift plus a learned diagonal diffusion. All rollouts except the NSDE use the true $\varsigma$, an oracle that favours the surrogates; the NSDE uses no true parameter.

\subsection{Evaluation protocol} Every surrogate passes through the identical instrument and is compared with a held-out true realization by pass/fail checks in four rungs. (1) Invariant statistics: $x$-marginal Wasserstein-1, $x$-autocorrelation distance, wing-switch ratio and lobe balance $f_+$, each of which must lie within the range (``floor'') reached by \floorN{} independent true realizations.
(2) Static susceptibility: $\chi(0)$ must lie within the range of \chiZeroBandn{} true realizations.
(3--4) Mode-resolved test, scored separately on the two slow real modes and the four oscillatory ones (Fig.~\ref{fig:main}c): for each leading eigenfunction we measure the alignment of the surrogate's eigenfunction $\tilde\varphi$ with the reference $\varphi$ by the overlap $O=|\langle\varphi,\tilde\varphi\rangle_{\mu_0}|/(\|\varphi\|_{\mu_0}\|\tilde\varphi\|_{\mu_0})\in[0,1]$, the cosine of the angle between them weighted by the invariant measure ($O=1$: same function up to scale). The mode passes if $O$ exceeds what independent true realizations typically reach, their mean overlap minus three standard deviations (clipped to $[0.5, 0.99]$ for the oscillatory modes).
A rollout fails a rung if it fails any check in it, and a model's failure rate on a rung is the fraction of its rollouts that fail. The failure rate of \Truen{} fresh true realizations, scored as a perfect surrogate, is each rung's false-alarm rate. We call a model's failure rate consistent with the truth when a Fisher exact test against the true realizations gives $p \ge \paramAlpha{}$. Table~\ref{tab:main} reports medians over \nTrainSeeds{} training seeds $\times$ \nRolloutSeeds{} rollout seeds (SINDy: \SindyFulln{} rollouts); per-rollout failure rates are tabulated in Table~\ref{tab:flags}.

\section{Surrogate models}
\label{app:models}

All five surrogates, the exploratory transformer of Sec.~\ref{app:transformer} and the discretization control of Sec.~\ref{app:gzero} are evaluated through the identical instrument of Sec.~\ref{app:operators}: same dictionary ($K_{\max}=\paramKmax{}$, size \paramDict{}), same TSVD threshold, same sampling step $\tau=\paramTau{}$, same forcing direction; the surrogates are trained on unperturbed trajectories only. The transformer is a first step from the formulation ladder toward the architectures of current weather emulators; it was run with a fixed configuration and no tuning, so it is reported separately from the ladder. The training and rollout settings and seeds of every model are listed in Table~\ref{tab:params} (App.~\ref{app:tables}).

\subsection{SINDy}
\label{app:sindy}

We fit the drift by sequentially thresholded least squares (STLSQ) on a polynomial library of degree \paramSindyDegree{}, with threshold \paramSindyThreshold{} on standardized columns.
Under Euler-Maruyama the one-step increment of \eqref{eq:sde} at the raw step satisfies $\mathbb{E}[\Delta x/\Delta t\mid x]=F(x)$, so regressing increments on state recovers the drift with the noise averaged out \cite{boninsegna_sparse_2018}; at a coarser step its conditional mean is the secant drift (Sec.~\ref{app:mlp}).
The library contains the true Lorenz-63 equations, so SINDy is a positive control on the whole pipeline rather than a competitive baseline: if the instrument cannot certify a model whose functional form is exactly right, nothing else it reports is trustworthy.
The positive control is trained on the same \fitTFull{} time units of the master as the neural models, at the raw step. As a robustness check on record length we also fit \fitTHundred{} time units (SINDy-100; this shorter fit was the control before the master trajectory was refitted, and the full-record fit replaced it after the results were seen).
SINDy-100 acquires symmetry-breaking constant terms that skew the lobes ($f_+ = \SINDyhfplus$); it fails invariant statistics in every rollout yet passes both slow modes in every rollout, the converse of the dissociation in Sec.~\ref{sec:exp}.
Table~\ref{tab:appx-sindy} gives both fits' coefficients and Table~\ref{tab:appx-rmse} the drift RMSE on the grid of \gridN{} points drawn from the held-out true realization (seed \gridSeed{}); the MLP-drift row is one representative fit (initialization seed \paramMlpInitSeed{}).
% GENERATED by scripts/camera_ready/make_appendix_tables.py from results/camera_ready/camera_ready_metrics.json via numbers.tex -- do not edit.
% HAND-EDITED 2026-10-06 (shorter caption/layout for the extended version): port to the generator before the next `make paper`, or this edit is overwritten.
\begin{table}[t]
\centering
\scriptsize
\begin{minipage}[t]{0.5\linewidth}
\centering
\caption{\textbf{SINDy coefficients}, true and fitted (large terms: above the STLSQ threshold).}
\label{tab:appx-sindy}
\setlength{\tabcolsep}{3pt}
\begin{tabular}{@{}lccc@{}}
\toprule
term & true & SINDy-100 & SINDy (full) \\
\midrule
$\sigma$ ($y$ in $\dot x$) & \coefSigmaTrue{} & \coefSigmaHundred{} & \coefSigmaFull{} \\
$-$($x$ in $\dot x$) & \coefMinusXTrue{} & \coefMinusXHundred{} & \coefMinusXFull{} \\
$\rho$ ($x$ in $\dot y$) & \coefRhoTrue{} & \coefRhoHundred{} & \coefRhoFull{} \\
$\beta$ ($-z$ in $\dot z$) & \coefBetaTrue{} & \coefBetaHundred{} & \coefBetaFull{} \\
const in $\dot x$ & \coefConstDxTrue{} & \coefConstDxHundred{} & \coefConstDxFull{} \\
const in $\dot y$ & \coefConstDyTrue{} & \coefConstDyHundred{} & \coefConstDyFull{} \\
const in $\dot z$ & \coefConstDzTrue{} & \coefConstDzHundred{} & \coefConstDzFull{} \\
active terms & \coefActiveTrue{} & \coefActiveHundred{} & \coefActiveFull{} \\
large terms & \coefBigTrue{} & \coefBigHundred{} & \coefBigFull{} \\
span (t.u.) / pairs & -- & \fitTHundred{} / \fitPairsHundred{} & \fitTFull{} / \fitPairsFull{} \\
$f_+$ of training window & -- & \fitFPlusHundred{} & \fitFPlusFull{} \\
\bottomrule
\end{tabular}
\end{minipage}\hfill
\begin{minipage}[t]{0.47\linewidth}
\centering
\caption{\textbf{Drift RMSE} against the true drift on the \gridN{}-point grid.}
\label{tab:appx-rmse}
\setlength{\tabcolsep}{3pt}
\begin{tabular}{@{}lccc@{}}
\toprule
& $x$ & $y$ & $z$ \\
\midrule
SINDy (full data) & \rmseFullX{} & \rmseFullY{} & \rmseFullZ{} \\
SINDy-100 & \rmseHundredX{} & \rmseHundredY{} & \rmseHundredZ{} \\
NODE (median, \rmseNODEn{} fits) & \rmseNODEXmed{} & \rmseNODEYmed{} & \rmseNODEZmed{} \\
MLP-drift & \rmseMLPdriftX{} & \rmseMLPdriftY{} & \rmseMLPdriftZ{} \\
secant bias (exact) & \rmseSecantX{} & \rmseSecantY{} & \rmseSecantZ{} \\
secant bias $(\Delta t/\paramTwo{})J_FF$ & \rmseSecantLeadX{} & \rmseSecantLeadY{} & \rmseSecantLeadZ{} \\
\bottomrule
\end{tabular}
\end{minipage}
\end{table}

\subsection{MLP}
\label{app:mlp}

The MLP is a map $\mathbb{R}^3\to\mathbb{R}^3$ with \paramMlpDepth{} hidden tanh layers of width \paramMlpWidth{} and standardized inputs and outputs, trained with Adam (learning rate $\paramMlpLr{}$, batch \paramMlpBatch{})~\cite{kingma_adam_2017} on the mean-squared error, with early stopping on a random \paramMlpValPct{}\% validation split (patience \paramMlpPatience{}, at most \paramMlpEpochs{} epochs).
It is trained on \paramTrainTUMLPdrift{} time units of the master at $\Delta t_{\mathrm{fit}}=\paramDtFitMLPdrift{}$, either on the increment $\Delta x/\Delta t$ (MLP-drift) or on the next state (MLP-flow); the two share data, architecture and seed, so the target is the only difference.
A pre-refit sweep over width and depth did not reduce the drift error.
At $\Delta t_{\mathrm{fit}}$ the conditional mean of $\Delta x/\Delta t$ is not $F$ but the secant drift $F+(\Delta t/2)\,J_F F+O(\Delta t^2)$ (the It\^o term $\tfrac{\varsigma^2}{2}\nabla^2F$ vanishes because each component of the Lorenz-63 drift is at most bilinear), and the drift MLP learns it: its error against $F$ matches the secant bias of its target to within \MLPdriftSecantMaxDevPct{}\% per component (Table~\ref{tab:appx-rmse}).
The MLP-drift is integrated with Euler--Maruyama at $\delta t=\paramDtRollMLPdrift{}$.

\subsection{Reservoir Computing}
\label{app:rc}

The reservoir is a leaky echo-state network~\cite{jaeger__2001} with $N=\paramRcN{}$ units, spectral radius \paramRcRadius{}, leak rate \paramRcLeak{}, input scaling \paramRcInput{}, washout \paramRcWashout{} steps and ridge readout regularization $\paramRcRidge{}$, trained on \paramTrainTURC{} time units of the master at $\Delta t_{\mathrm{fit}}=\paramDtFitRC{}$ and advanced as a native flow map with injected noise.
Like every other surrogate, it is trained on the Euler--Maruyama master, i.e.\ on noisy data from the stochastic system; its hyperparameters were not retuned for it.
The configuration was selected on valid-prediction horizon against a reference trajectory of the deterministic Lorenz-63 system, subject to an attractor gate (bounded, bimodal, plausible wing-switch ratio); response metrics were never used for selection.
Two implementation faults had to be fixed first, both of which silently degrade the attractor rather than failing loudly: an off-by-one in the closed-loop rollout seed, and a washout discontinuity in which warm-up ended at the last training state while the rollout restarted from the first.

\subsection{Neural ODE and neural SDE}
\label{app:neural}

The NODE uses the MLP network of Sec.~\ref{app:mlp} as a drift $F_\theta$, trained through a differentiable RK4 step at $h=\Delta t_{\mathrm{fit}}=\paramDtFitNODE{}$ on the same data with the same validation fraction: from each window start the model is unrolled for $L$ steps, the mean-squared error against the master is taken in standardized coordinates and minimized with Adam, with early stopping on the held-out multi-step validation MSE.
We swept $L\in\{\paramNodeSweepLOne{},\paramNodeSweepLTwo{}\}$ and learning rates $\{\paramNodeSweepLrOne{},\paramNodeSweepLrTwo{}\}$ and selected by validation MSE at a common horizon $L=\paramHorizonC{}$ ($L=\paramNodeLSel{}$, learning rate $\paramNodeLrSel{}$; Table~\ref{tab:appx-sweep}); the NODE is rolled out with the same RK4 step and oracle noise.
The NSDE freezes the NODE drift and adds a diagonal diffusion $\sigma_\phi(x)$, a network with one hidden layer of \paramSigmaHidden{} tanh units and a softplus output, fitted by the Gaussian negative log-likelihood of the frozen drift's one-step residuals; it is rolled out with its learned noise and uses no true parameter.
The NSDE's learned noise matches the one-step noise measured with the true drift to within \NSDEsigmaRelDiffPctMax{}\% per component (largest in $x$); in $x$ both lie \NSDEgapXPct\% below $\varsigma$, close to the \relaxGapXPct\% predicted by in-step relaxation under the linear damping alone (Table~\ref{tab:appx-noise}).
Both models use \nTrainSeeds{} training seeds $\times$ \nRolloutSeeds{} rollout seeds; the pilot runs used to set up the pipeline are excluded from every table.
Refitting the seed-0 NODE and NSDE from scratch changes every reported metric by at most \nodeReproMaxFracRange{} of its range across seeds: on the same GPU the refits and rollouts were bit-identical; reproducibility across hardware was not tested.

\subsection{Exploratory: a causal transformer flow map}
\label{app:transformer}
Transformer weather emulators such as Pangu-Weather~\cite{bi_accurate_2023}, FourCastNet~\cite{pathak_fourcastnet_2022} and Aurora~\cite{bodnar_foundation_2025} are, in the taxonomy of Sec.~\ref{sec:method}, learned one-step flow maps rolled out autoregressively, with attention mixing spatial tokens.
A three-variable system has no spatial tokens, so the analogue we can test here is attention over time.
As an exploratory model, with a configuration fixed before the run and no tuning, we add a decoder-only causal transformer~\cite{vaswani_attention_2017,radford_language_2019, zhou_informer_2021} that maps the last \trfContext{} states to the next-state increment (embedding dimension \trfDim{}, \trfLayers{} layers, \trfHeads{} heads, \trfParams{} parameters).
It is trained on the MLP-flow data with a teacher-forced one-step loss, then fine-tuned on \trfRolloutSteps{}-step autoregressive rollouts, as recent emulators are~\cite{lam_learning_2023,bodnar_foundation_2025}, and rolled out with the same oracle noise as MLP-flow.
Because the true state is fully observed and Markov, the context cannot improve the one-step prediction: the model is a sequence model with memory, like RC, not a stand-in for spatial-attention emulators.

It passed the smoke-test gate fixed in advance, and all \nRunsTransformer{} rollouts over \nTrainSeeds{} training $\times$ \nRolloutSeeds{} rollout seeds were scored, none diverging. Its pattern matches NODE and NSDE (Table~\ref{tab:flags}): it fails invariant statistics in \pct{\flagInvTransformer} of rollouts but the slow modes in only \pct{\flagSlowTransformer}, and the oscillatory modes in \pct{\flagOscTransformer}, mostly through the individual near-degenerate eigenvectors while recovering their joint subspace (merged-subspace failure rate \pct{\flagOscMergedTransformer}; Table~\ref{tab:appx-modes}). Of its rollouts, \nXTransformer{} lie inside all four invariant-statistics floors and the $\chi(0)$ band, and \nXslowfailTransformer{} of these fail the slow-mode test. It differs from MLP-flow in both architecture and the rollout fine-tuning, so its improvement over MLP-flow cannot be attributed to either alone; that it, like the NODE, was trained on multi-step rollouts is consistent with the training formulation, rather than the architecture, setting which properties a model gets right.

\section{Experiments}
\subsection{Comparison across models}
\label{sec:exp}
% GENERATED by scripts/camera_ready/make_fig1.py from results/camera_ready/camera_ready_metrics.json via numbers.tex -- do not edit.
\begin{figure}[!t]
\centering
\includegraphics[width=\textwidth]{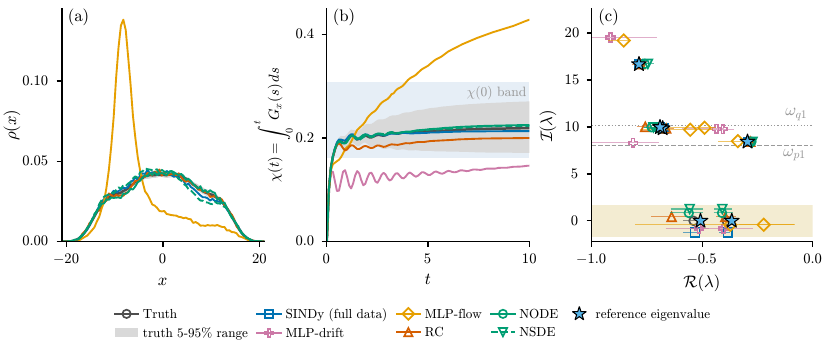}
\caption{\textbf{Invariant statistics, static susceptibility and resonances, per model.} (a)~$x$-marginal density: truth median and \paramFigQlo{}--\paramFigQhi{}\% range over \Truen{} realizations; per model, the rollout with median $\chi(0)$ (MLP-drift not shown). (b)~Accumulated response $\chi(t)=\int_0^t G_x(s)\,ds$, which tends to $\chi(0)$: truth median and \paramFigQlo{}--\paramFigQhi{}\% range; per model, the median over rollouts; blue band: the $\chi(0)$ band of Table~\ref{tab:main}. (c)~Leading Koopman/RP resonances ($\mathrm{Im}\,\lambda\ge\paramZero{}$): reference (stars) and, per model, the median of each leading group over rollouts, with its \paramFigQlo{}--\paramFigQhi{}\% range in $\mathrm{Re}\,\lambda$; real modes spread vertically for visibility. Dashed and dotted lines: the shortest-UPO twirling frequency $\omega_{p1}$ and the intra-lobe turnover frequency $\omega_{q1}$~\cite{lucarini_general_2026}.}
\label{fig:main}
\end{figure}

% GENERATED by scripts/camera_ready/make_table1_values.py from results/camera_ready/camera_ready_metrics.json via numbers.tex -- do not edit.
% HAND-EDITED 2026-10-06 (shorter caption/layout for the extended version): port to the generator before the next `make paper`, or this edit is overwritten.
\begin{table}[t]
\centering
\footnotesize
\definecolor{flaghi}{HTML}{B2182B}
\caption{\textbf{Invariant statistics, $\chi(0)$ and eigenfunction overlaps, per model.} Medians over scored rollouts; $\chi(0)$ as mean$_{(\mathrm{sd})}$. Threshold row: the invariant-statistics floors, the $\chi(0)$ band and the overlap threshold of each leading mode (columns labelled by reference eigenvalue or frequency). Red: median outside the threshold; bold: best model per column. $n$: scored rollouts (RC: \nDivRC{} diverged).}
\label{tab:main}
\setlength{\tabcolsep}{1.3pt}
\begin{tabular*}{\linewidth}{@{\extracolsep{\fill}}lcccccccccccc}
\toprule
 & & \multicolumn{4}{c}{Invariant statistics} & & \multicolumn{2}{c}{Slow ($\mathrm{Re}\,\lambda$)} & \multicolumn{4}{c}{Oscillatory ($\omega$)} \\
\cmidrule(lr){3-6}\cmidrule(lr){8-9}\cmidrule(lr){10-13}
Model & $n$ & $W_1$ & ACF & switch & $f_+$ & $\chi(0)$ & $\refReTwo{}$ & $\refReThree{}$ & $\refFreqOne{}$ & $\refFreqFour{}$ & $\refFreqFive{}$ & $\refFreqSix{}$ \\
\midrule
{\scriptsize\itshape threshold} & & {\scriptsize $\le$\tabFloorWOne{}} & {\scriptsize $\le$\tabFloorAcf{}} & {\scriptsize \tabFloorSwitchlo{}--\tabFloorSwitchhi{}} & {\scriptsize \tabFloorFPluslo{}--\tabFloorFPlushi{}} & {\scriptsize \tabBandlo{}--\tabBandhi{}} & {\scriptsize $\ge$\tabThrTwo{}} & {\scriptsize $\ge$\tabThrThree{}} & {\scriptsize $\ge$\tabThrOne{}} & {\scriptsize $\ge$\tabThrFour{}} & {\scriptsize $\ge$\tabThrFive{}} & {\scriptsize $\ge$\tabThrSix{}} \\
\midrule
Truth & \nScoredTrue{} & \tabVWOneTrue{} & \tabVAcfTrue{} & \tabVSwitchTrue{} & \tabVFPlusTrue{} & $\tabChiMeanTrue{}_{(\tabChiSdTrue{})}$ & \tabOvTrueTwo{} & \tabOvTrueThree{} & \tabOvTrueOne{} & \tabOvTrueFour{} & \tabOvTrueFive{} & \tabOvTrueSix{} \\
\midrule
SINDy & \nScoredSindyFull{} & \textbf{\tabVWOneSindyFull{}} & \textbf{\tabVAcfSindyFull{}} & \textbf{\tabVSwitchSindyFull{}} & \tabVFPlusSindyFull{} & $\tabChiMeanSindyFull{}_{(\tabChiSdSindyFull{})}$ & \textbf{\tabOvSindyFullTwo{}} & \textbf{\tabOvSindyFullThree{}} & \textbf{\tabOvSindyFullOne{}} & \textbf{\tabOvSindyFullFour{}} & \textbf{\tabOvSindyFullFive{}} & \textbf{\tabOvSindyFullSix{}} \\
MLP-drift & \nScoredMLPdrift{} & \textcolor{flaghi}{\tabVWOneMLPdrift{}} & \textcolor{flaghi}{\tabVAcfMLPdrift{}} & \textcolor{flaghi}{\tabVSwitchMLPdrift{}} & \textcolor{flaghi}{\tabVFPlusMLPdrift{}} & \textcolor{flaghi}{$\tabChiMeanMLPdrift{}_{(\tabChiSdMLPdrift{})}$} & \textcolor{flaghi}{\tabOvMLPdriftTwo{}} & \textcolor{flaghi}{\tabOvMLPdriftThree{}} & \textcolor{flaghi}{\tabOvMLPdriftOne{}} & \textcolor{flaghi}{\tabOvMLPdriftFour{}} & \textcolor{flaghi}{\tabOvMLPdriftFive{}} & \textcolor{flaghi}{\tabOvMLPdriftSix{}} \\
MLP-flow & \nScoredMLPflow{} & \textcolor{flaghi}{\tabVWOneMLPflow{}} & \textcolor{flaghi}{\tabVAcfMLPflow{}} & \textcolor{flaghi}{\tabVSwitchMLPflow{}} & \textcolor{flaghi}{\tabVFPlusMLPflow{}} & \textcolor{flaghi}{$\tabChiMeanMLPflow{}_{(\tabChiSdMLPflow{})}$} & \textcolor{flaghi}{\tabOvMLPflowTwo{}} & \textcolor{flaghi}{\tabOvMLPflowThree{}} & \textcolor{flaghi}{\tabOvMLPflowOne{}} & \textcolor{flaghi}{\tabOvMLPflowFour{}} & \textcolor{flaghi}{\tabOvMLPflowFive{}} & \textcolor{flaghi}{\tabOvMLPflowSix{}} \\
RC & \nScoredRC{} & \tabVWOneRC{} & \tabVAcfRC{} & \textcolor{flaghi}{\tabVSwitchRC{}} & \textbf{\tabVFPlusRC{}} & $\tabChiMeanRC{}_{(\tabChiSdRC{})}$ & \textcolor{flaghi}{\tabOvRCTwo{}} & \textcolor{flaghi}{\tabOvRCThree{}} & \tabOvRCOne{} & \textcolor{flaghi}{\tabOvRCFour{}} & \textcolor{flaghi}{\tabOvRCFive{}} & \tabOvRCSix{} \\
NODE & \nScoredNODE{} & \textcolor{flaghi}{\tabVWOneNODE{}} & \textcolor{flaghi}{\tabVAcfNODE{}} & \textcolor{flaghi}{\tabVSwitchNODE{}} & \tabVFPlusNODE{} & {\boldmath$\tabChiMeanNODE{}_{(\tabChiSdNODE{})}$} & \tabOvNODETwo{} & \tabOvNODEThree{} & \tabOvNODEOne{} & \textcolor{flaghi}{\tabOvNODEFour{}} & \textcolor{flaghi}{\tabOvNODEFive{}} & \tabOvNODESix{} \\
NSDE & \nScoredNSDE{} & \textcolor{flaghi}{\tabVWOneNSDE{}} & \textcolor{flaghi}{\tabVAcfNSDE{}} & \textcolor{flaghi}{\tabVSwitchNSDE{}} & \tabVFPlusNSDE{} & $\tabChiMeanNSDE{}_{(\tabChiSdNSDE{})}$ & \tabOvNSDETwo{} & \tabOvNSDEThree{} & \tabOvNSDEOne{} & \textcolor{flaghi}{\tabOvNSDEFour{}} & \textcolor{flaghi}{\tabOvNSDEFive{}} & \tabOvNSDESix{} \\
\bottomrule
\end{tabular*}
\end{table}

% GENERATED by scripts/camera_ready/make_table1.py from results/camera_ready/camera_ready_metrics.json via numbers.tex -- do not edit.
% HAND-EDITED 2026-10-06 (shorter caption/layout for the extended version): port to the generator before the next `make paper`, or this edit is overwritten.
\begin{table}[t]
\centering
\small
\definecolor{flaghi}{HTML}{B2182B}\definecolor{flaglo}{HTML}{707070}
\caption{\textbf{Per-rollout failure rates} (\% of scored rollouts; redder = fails more often), per rung and per invariant-statistics check, with the thresholds of Table~\ref{tab:main}. Truth: fresh true realizations scored as a perfect surrogate, i.e.\ the false-alarm rate; being inside the $\chi(0)$ band by construction, its leave-one-out out-rate is given in parentheses. $\chi(0)$: mean$_{(\mathrm{sd})}$. Bottom rows: the exploratory transformer (Sec.~\ref{app:transformer}) and the true drift rolled out on the NODE and NSDE integrator paths (G0, Sec.~\ref{app:gzero}). div.: diverged rollouts (non-finite, or a coordinate beyond $\gateFactor{}\times$ its true maximum), left out of the rates.}
\label{tab:flags}
\setlength{\tabcolsep}{3.4pt}
\begin{tabular*}{\linewidth}{@{\extracolsep{\fill}}lcccccccccc}
\toprule
 & & \multicolumn{5}{c}{Invariant statistics (\%)} & \multicolumn{2}{c}{$\chi(0)$} & \multicolumn{2}{c}{Modes (\%)} \\
\cmidrule(lr){3-7}\cmidrule(lr){8-9}\cmidrule(lr){10-11}
Model & $n$ (div.) & any & $W_1$ & ACF & switch & $f_+$ & mean$_{(\mathrm{sd})}$ & out (\%) & slow & oscill. \\
\midrule
Truth & \nScoredTrue{}\,(\nDivTrue{}) & \textcolor{flaghi!4!flaglo}{\pctInvTrue{}} & \textcolor{flaghi!0!flaglo}{\pctOutWOneTrue{}} & \textcolor{flaghi!0!flaglo}{\pctOutAcfTrue{}} & \textcolor{flaghi!4!flaglo}{\pctOutSwitchTrue{}} & \textcolor{flaghi!0!flaglo}{\pctOutFPlusTrue{}} & $\tabChiMeanTrue{}_{(\tabChiSdTrue{})}$ & \textcolor{flaghi!0!flaglo}{\pctChiTrue{}}\,(\pctBandloo{}) & \textcolor{flaghi!6!flaglo}{\pctSlowTrue{}} & \textcolor{flaghi!2!flaglo}{\pctOscTrue{}} \\
\midrule
SINDy & \nScoredSindyFull{}\,(\nDivSindyFull{}) & \textcolor{flaghi!20!flaglo}{\pctInvSindyFull{}} & \textcolor{flaghi!10!flaglo}{\pctOutWOneSindyFull{}} & \textcolor{flaghi!0!flaglo}{\pctOutAcfSindyFull{}} & \textcolor{flaghi!0!flaglo}{\pctOutSwitchSindyFull{}} & \textcolor{flaghi!20!flaglo}{\pctOutFPlusSindyFull{}} & $\tabChiMeanSindyFull{}_{(\tabChiSdSindyFull{})}$ & \textcolor{flaghi!0!flaglo}{\pctChiSindyFull{}} & \textcolor{flaghi!0!flaglo}{\pctSlowSindyFull{}} & \textcolor{flaghi!10!flaglo}{\pctOscSindyFull{}} \\
MLP-drift & \nScoredMLPdrift{}\,(\nDivMLPdrift{}) & \textcolor{flaghi!100!flaglo}{\pctInvMLPdrift{}} & \textcolor{flaghi!100!flaglo}{\pctOutWOneMLPdrift{}} & \textcolor{flaghi!100!flaglo}{\pctOutAcfMLPdrift{}} & \textcolor{flaghi!100!flaglo}{\pctOutSwitchMLPdrift{}} & \textcolor{flaghi!96!flaglo}{\pctOutFPlusMLPdrift{}} & $\tabChiMeanMLPdrift{}_{(\tabChiSdMLPdrift{})}$ & \textcolor{flaghi!91!flaglo}{\pctChiMLPdrift{}} & \textcolor{flaghi!100!flaglo}{\pctSlowMLPdrift{}} & \textcolor{flaghi!100!flaglo}{\pctOscMLPdrift{}} \\
MLP-flow & \nScoredMLPflow{}\,(\nDivMLPflow{}) & \textcolor{flaghi!100!flaglo}{\pctInvMLPflow{}} & \textcolor{flaghi!100!flaglo}{\pctOutWOneMLPflow{}} & \textcolor{flaghi!93!flaglo}{\pctOutAcfMLPflow{}} & \textcolor{flaghi!100!flaglo}{\pctOutSwitchMLPflow{}} & \textcolor{flaghi!100!flaglo}{\pctOutFPlusMLPflow{}} & $\tabChiMeanMLPflow{}_{(\tabChiSdMLPflow{})}$ & \textcolor{flaghi!99!flaglo}{\pctChiMLPflow{}} & \textcolor{flaghi!98!flaglo}{\pctSlowMLPflow{}} & \textcolor{flaghi!100!flaglo}{\pctOscMLPflow{}} \\
RC & \nScoredRC{}\,(\nDivRC{}) & \textcolor{flaghi!63!flaglo}{\pctInvRC{}} & \textcolor{flaghi!3!flaglo}{\pctOutWOneRC{}} & \textcolor{flaghi!27!flaglo}{\pctOutAcfRC{}} & \textcolor{flaghi!58!flaglo}{\pctOutSwitchRC{}} & \textcolor{flaghi!6!flaglo}{\pctOutFPlusRC{}} & $\tabChiMeanRC{}_{(\tabChiSdRC{})}$ & \textcolor{flaghi!8!flaglo}{\pctChiRC{}} & \textcolor{flaghi!76!flaglo}{\pctSlowRC{}} & \textcolor{flaghi!100!flaglo}{\pctOscRC{}} \\
NODE & \nScoredNODE{}\,(\nDivNODE{}) & \textcolor{flaghi!93!flaglo}{\pctInvNODE{}} & \textcolor{flaghi!58!flaglo}{\pctOutWOneNODE{}} & \textcolor{flaghi!56!flaglo}{\pctOutAcfNODE{}} & \textcolor{flaghi!75!flaglo}{\pctOutSwitchNODE{}} & \textcolor{flaghi!63!flaglo}{\pctOutFPlusNODE{}} & $\tabChiMeanNODE{}_{(\tabChiSdNODE{})}$ & \textcolor{flaghi!9!flaglo}{\pctChiNODE{}} & \textcolor{flaghi!24!flaglo}{\pctSlowNODE{}} & \textcolor{flaghi!86!flaglo}{\pctOscNODE{}} \\
NSDE & \nScoredNSDE{}\,(\nDivNSDE{}) & \textcolor{flaghi!100!flaglo}{\pctInvNSDE{}} & \textcolor{flaghi!60!flaglo}{\pctOutWOneNSDE{}} & \textcolor{flaghi!69!flaglo}{\pctOutAcfNSDE{}} & \textcolor{flaghi!79!flaglo}{\pctOutSwitchNSDE{}} & \textcolor{flaghi!61!flaglo}{\pctOutFPlusNSDE{}} & $\tabChiMeanNSDE{}_{(\tabChiSdNSDE{})}$ & \textcolor{flaghi!9!flaglo}{\pctChiNSDE{}} & \textcolor{flaghi!23!flaglo}{\pctSlowNSDE{}} & \textcolor{flaghi!83!flaglo}{\pctOscNSDE{}} \\
\midrule
Transformer (expl.) & \nScoredTransformer{}\,(\nDivTransformer{}) & \textcolor{flaghi!79!flaglo}{\pctInvTransformer{}} & \textcolor{flaghi!51!flaglo}{\pctOutWOneTransformer{}} & \textcolor{flaghi!29!flaglo}{\pctOutAcfTransformer{}} & \textcolor{flaghi!53!flaglo}{\pctOutSwitchTransformer{}} & \textcolor{flaghi!64!flaglo}{\pctOutFPlusTransformer{}} & $\tabChiMeanTransformer{}_{(\tabChiSdTransformer{})}$ & \textcolor{flaghi!13!flaglo}{\pctChiTransformer{}} & \textcolor{flaghi!19!flaglo}{\pctSlowTransformer{}} & \textcolor{flaghi!73!flaglo}{\pctOscTransformer{}} \\
\midrule
G0, NODE path & \nScoredGzeroH{}\,(\nDivGzeroH{}) & \textcolor{flaghi!100!flaglo}{\pctInvGzeroH{}} & \textcolor{flaghi!0!flaglo}{\pctOutWOneGzeroH{}} & \textcolor{flaghi!60!flaglo}{\pctOutAcfGzeroH{}} & \textcolor{flaghi!100!flaglo}{\pctOutSwitchGzeroH{}} & \textcolor{flaghi!0!flaglo}{\pctOutFPlusGzeroH{}} & $\tabChiMeanGzeroH{}_{(\tabChiSdGzeroH{})}$ & \textcolor{flaghi!0!flaglo}{\pctChiGzeroH{}} & \textcolor{flaghi!10!flaglo}{\pctSlowGzeroH{}} & \textcolor{flaghi!40!flaglo}{\pctOscGzeroH{}} \\
G0, NSDE path & \nScoredGzeroHSig{}\,(\nDivGzeroHSig{}) & \textcolor{flaghi!100!flaglo}{\pctInvGzeroHSig{}} & \textcolor{flaghi!0!flaglo}{\pctOutWOneGzeroHSig{}} & \textcolor{flaghi!70!flaglo}{\pctOutAcfGzeroHSig{}} & \textcolor{flaghi!100!flaglo}{\pctOutSwitchGzeroHSig{}} & \textcolor{flaghi!0!flaglo}{\pctOutFPlusGzeroHSig{}} & $\tabChiMeanGzeroHSig{}_{(\tabChiSdGzeroHSig{})}$ & \textcolor{flaghi!30!flaglo}{\pctChiGzeroHSig{}} & \textcolor{flaghi!0!flaglo}{\pctSlowGzeroHSig{}} & \textcolor{flaghi!40!flaglo}{\pctOscGzeroHSig{}} \\
\bottomrule
\end{tabular*}
\end{table}

The false-alarm rates are low: fresh true realizations fail the invariant statistics in \pct{\flagInvTrue} of cases, the slow modes in \pct{\flagSlowTrue} and the oscillatory modes in \pct{\flagOscTrue}, and a true realization left out of the $\chi(0)$ band falls outside it in \pct{\chiZeroBandloo} of cases (Table~\ref{tab:flags}). SINDy trained on the full record recovers the structural coefficients to within \coefMaxRelErrPct{}\% and fails at rates consistent with the truth on every rung, with no mode pass decided by the 0.99 cap: the instrument certifies a correct model rather than only rejecting wrong ones.

The reservoir computer reproduces the invariant measure closely ($W_1$ inside its floor in \inWOneRC{} of \RCn{} rollouts; Fig.~\ref{fig:main}a), and \RCnX{} of its rollouts lie inside all four invariant-statistics floors and the $\chi(0)$ band. Of these, \RCnXslowfail{} fail the slow-mode test, against \TruenXslowfail{} of the \TruenX{} truths that pass the same checks (Fisher $p = \RCexhibitFisherp{}$; Fig.~\ref{fig:appx-modes}b maps this dissociation rollout by rollout). The failure is not a disguised switching error: RC fails the slow modes in \pct{\RCslowfailInSwitch} of rollouts inside the switch-ratio floor and \pct{\RCslowfailOutSwitch} outside it ($p = \RCfisherSwitchSlowp$). Spectrally, RC places its two leading oscillatory resonances within \RCfreqOffsetMaxPct{}\% of the reference frequencies, but its second slow mode decays faster than in any true realization in \pct{\RCrealTwoFasterFrac} of rollouts (median $\mathrm{Re}\,\lambda = \eigRealTwoReRCmed$ against the true range $[\eigRealTwoReTruemin, \eigRealTwoReTruemax]$; Fig.~\ref{fig:main}c), and its slow eigenfunctions pass in under a third of rollouts.

The MLP baselines fail every rung (Table~\ref{tab:main}). The drift MLP learns the $O(\Delta t)$ secant drift rather than $F$ (Sec.~\ref{app:mlp}), and its Euler rollout collapses onto one wing in \MLPdriftonewing{} of \MLPdriftn{} rollouts. Fitted as a flow map, the same network restores two-wing switching (\nOneWingMLPflow{} of \nScoredMLPflow{} rollouts one-wing) but not the measure (median $f_+=\MLPflowfplus{}$); its slowest real mode decays too slowly ($\MLPflowslow$ against the true $\TrueRealOneRemed$), so only \pct{\MLPflowacc} of its $\chi(0)$ has accrued by $t = 1$, against a median \pct{\chiaccTruth} for the truths, and $\chi(0)$ overshoots. Trained through the integrator on multi-step windows, the same network (NODE) and its stochastic extension (NSDE) show the reverse of RC's pattern: both over-switch in every rollout (median switch ratio \switchRatioNODEmed{} and \switchRatioNSDEmed{}) and fail invariant statistics in \pct{\NODEinvflag} and \pct{\NSDEinvflag} of rollouts, yet the slow modes in only \pct{\NODEslowflag} and \pct{\NSDEslowflag}; the over-switching is learned, not numerical (Sec.~\ref{app:gzero}). An exploratory causal transformer fine-tuned on short rollouts, as recent weather emulators are, shows the same pattern (Sec.~\ref{app:transformer}).

For SINDy, RC, NODE and NSDE the $\chi(0)$ failure rate never exceeds \pct{\chiflagmax}, while their slow-mode failure rates range from \pct{\slowflagmin} to \pct{\slowflagmax}. This is expected for a time-integrated quantity: for forcing and observing $x$, $\chi(0)$ is set mostly by fast relaxation: across the true realizations a median \pct{\chiaccTruth} of it has accrued by $t = 1$ (Fig.~\ref{fig:main}b), before the slow modes (e-folding times \eFoldRealTwoRef{} and \eFoldRealOneRef) have relaxed. Forcings or observables coupled to lobe occupation may load on exactly these modes; the resonances and eigenfunctions, unlike $\chi(0)$, do not depend on that choice.

\subsection{Instrument diagnostics and calibration}
\label{app:diagnostics}

Several checks quantify how far the finite-dictionary, finite-trajectory instrument departs from its exact counterpart, and they bound what the comparisons of Sec.~\ref{sec:exp} can resolve; their tables are collected in App.~\ref{app:tables}.
With $\Gamma=-\partial_x\ln\rho_0$ and forcing along $x$, the value of the Green's function at zero lag is fixed analytically: from \eqref{eq:green},
\begin{equation}
\label{eq:g0-identity}
G_x(0)=\int_{\mathcal M} x\,\big(-\partial_x\ln\rho_0\big)\,d\mu_0
=-\int_{\mathcal M} x\,\partial_x\rho_0\,d\mathbf{x}
=\int_{\mathcal M}\rho_0\,\partial_x x\,d\mathbf{x}=1,
\end{equation}
the boundary term vanishing for the smooth, decaying $\rho_0$.
Nothing in the pipeline normalizes $G_x$, so the measured deviation from unity is a free, per-series estimate of dictionary-truncation error: $|G_x(0)-1|$ is at most \maxGdevTruthPct{}\% over the truths and at most \maxGdevAllPct{}\% over all scored rollouts.

The invariant-statistics floors are the maxima of each distance over the \floorN{} independent true realizations scored against the held-out reference realization (seed \seedRef{}) rather than the master trajectory, which is every model's training data; no floor fails its leave-one-out null test (Table~\ref{tab:appx-null}). The $\chi(0)$ band is the range of the \bandPooledn{} pooled true realizations, because the range of the independent-IC ensemble alone left \bandNullkOut{} of the \bandNullnTest{} truths scored through the runner outside (one-sided binomial $P=\bandNullp{}$ against its leave-one-out rate); Table~\ref{tab:appx-bands} gives the in-band counts under each alternative band.

$\chi(0)$ is evaluated in closed form through \eqref{eq:chi-discrete}, not by truncating $\tau\sum_n G_x(n\tau)$ at a finite lag. The slowest real eigenvalue varies little across the \Truen{} truths (coefficient of variation \TrueRealOneRecv{}), while the share of $\chi(0)$ carried by the leading resonances varies strongly (coefficient of variation \chiShareLeadTruthcv{} over \chiShareTruthn{} truths; Table~\ref{tab:appx-shares}).
Since $\chi(0)=-\sum_k\beta_k/\lambda_k$ and the $\lambda_k$ barely move, the finite-length variance enters mainly through the projection amplitudes $\beta_k$, hence through the Method-3 estimate of $\Gamma$ (Sec.~\ref{app:method3}), which is a time average along a single trajectory.

\subsection{Discretization control}
\label{app:gzero}

To separate learning error from the error of the NODE/NSDE integrator, we roll out the true drift through the same RK4 path at $h=\Delta t_{\mathrm{fit}}=\paramDtFitNODE{}$ with oracle noise, on both code paths (NODE: constant $\varsigma$; NSDE: the diffusion network replaced by the constant $\varsigma$), \Gzeron{} draws each (G0 rows of Table~\ref{tab:flags}).
The switch ratio falls below one in every draw (median \switchRatioGzeroHmed{} and \switchRatioGzeroHSigmed{}, maximum \switchRatioGzeroHmax{} and \switchRatioGzeroHSigmax{}), whereas NODE and NSDE lie above one in \NODEswitchAboveOne{} of \NODEn{} and \NSDEswitchAboveOne{} of \NSDEn{} rollouts, so their over-switching is learned rather than numerical.
The step shifts the median $\chi(0)$ by \Gzerochishift{} and \GzeroSigchishift{} from the truth median (\chiZeroTruemed{}) and leaves the slow modes largely intact (slow-mode failure rate at most \pct{\Gzeroslowflag}).

\subsection{Per-mode overlaps and threshold sensitivity}
\label{sec:modes}
% GENERATED by scripts/camera_ready/make_appendix_tables.py from results/camera_ready/camera_ready_metrics.json via numbers.tex, restyled (shaded cells) -- regenerate both together.
% HAND-EDITED 2026-10-06 (shorter caption/layout for the extended version): port to the generator before the next `make paper`, or this edit is overwritten.
\begin{table*}[t]
\centering
\scriptsize
\definecolor{flaghi}{HTML}{B2182B}
\caption{\textbf{Per-mode eigenfunction overlaps.} Each cell: median overlap $O$ over scored rollouts (top) and the percentage of rollouts passing the mode (bottom); redder = more rollouts fail. Threshold row: the applied threshold, uncapped value in parentheses. $\dagger$: more than half of the passes are decided by the \paramThrCap{} cap. Merged: the joint subspace of the near-degenerate pair (\mFreqAOne{} and \mFreqBOne{}; null $\mnullMeanOne{}\pm\mnullSdOne{}$ over the calibration truths). Bottom block: SINDy-100 (Sec.~\ref{app:sindy}) and the true drift on the NODE/NSDE integrator paths (Sec.~\ref{app:gzero}).}
\label{tab:appx-modes}
\setlength{\tabcolsep}{2.5pt}
\renewcommand{\arraystretch}{1.15}
\begin{tabular}{@{}lcccccccc@{}}
\toprule
 & & \multicolumn{2}{c}{Slow ($\mathrm{Re}\,\lambda$)} & \multicolumn{4}{c}{Oscillatory ($\omega$)} & \\
\cmidrule(lr){3-4}\cmidrule(lr){5-8}
Model & $n$ & \refReTwo{} & \refReThree{} & \refFreqOne{} & \refFreqFour{} & \refFreqFive{} & \refFreqSix{} & merged \\
\midrule
{\itshape threshold} & & \begin{tabular}[t]{@{}c@{}}\thrTwo{}\\{\color{black!65}(\thrUncTwo{})}\end{tabular} & \begin{tabular}[t]{@{}c@{}}\thrThree{}\\{\color{black!65}(\thrUncThree{})}\end{tabular} & \begin{tabular}[t]{@{}c@{}}\thrOne{}\\{\color{black!65}(\thrUncOne{})}\end{tabular} & \begin{tabular}[t]{@{}c@{}}\thrFour{}\\{\color{black!65}(\thrUncFour{})}\end{tabular} & \begin{tabular}[t]{@{}c@{}}\thrFive{}\\{\color{black!65}(\thrUncFive{})}\end{tabular} & \begin{tabular}[t]{@{}c@{}}\thrSix{}\\{\color{black!65}(\thrUncSix{})}\end{tabular} & \begin{tabular}[t]{@{}c@{}}\mthrOne{}\\{\color{black!65}(\mthrUncOne{})}\end{tabular} \\
\midrule
Truth & \nScoredTrue{} & \cellcolor{flaghi!2!white}\begin{tabular}[t]{@{}c@{}}0.992\\{\color{black!65}96\%}\end{tabular} & \cellcolor{flaghi!3!white}\begin{tabular}[t]{@{}c@{}}0.991\\{\color{black!65}94\%}\end{tabular} & \begin{tabular}[t]{@{}c@{}}0.999\\{\color{black!65}100\%}\end{tabular} & \cellcolor{flaghi!1!white}\begin{tabular}[t]{@{}c@{}}0.990\\{\color{black!65}98\%}\end{tabular} & \begin{tabular}[t]{@{}c@{}}0.991\\{\color{black!65}100\%}\end{tabular} & \begin{tabular}[t]{@{}c@{}}0.998\\{\color{black!65}100\%}\end{tabular} & \begin{tabular}[t]{@{}c@{}}0.998\\{\color{black!65}100\%}\end{tabular} \\[1pt]
\midrule
SINDy & \nScoredSindyFull{} & \begin{tabular}[t]{@{}c@{}}0.994\\{\color{black!65}100\%}\end{tabular} & \begin{tabular}[t]{@{}c@{}}0.986\\{\color{black!65}100\%}\end{tabular} & \begin{tabular}[t]{@{}c@{}}0.999\\{\color{black!65}100\%}\end{tabular} & \cellcolor{flaghi!5!white}\begin{tabular}[t]{@{}c@{}}0.990\\{\color{black!65}90\%}\end{tabular} & \cellcolor{flaghi!5!white}\begin{tabular}[t]{@{}c@{}}0.993\\{\color{black!65}90\%}\end{tabular} & \begin{tabular}[t]{@{}c@{}}0.998\\{\color{black!65}100\%}\end{tabular} & \begin{tabular}[t]{@{}c@{}}0.998\\{\color{black!65}100\%}\end{tabular} \\[1pt]
MLP-drift & \nScoredMLPdrift{} & \cellcolor{flaghi!55!white}\begin{tabular}[t]{@{}c@{}}0.043\\{\color{black!65}0\%}\end{tabular} & \cellcolor{flaghi!55!white}\begin{tabular}[t]{@{}c@{}}0.319\\{\color{black!65}0\%}\end{tabular} & \cellcolor{flaghi!55!white}\begin{tabular}[t]{@{}c@{}}0.607\\{\color{black!65}0\%}\end{tabular} & \cellcolor{flaghi!55!white}\begin{tabular}[t]{@{}c@{}}0.373\\{\color{black!65}0\%}\end{tabular} & \cellcolor{flaghi!55!white}\begin{tabular}[t]{@{}c@{}}0.290\\{\color{black!65}0\%}\end{tabular} & \cellcolor{flaghi!55!white}\begin{tabular}[t]{@{}c@{}}0.115\\{\color{black!65}0\%}\end{tabular} & \cellcolor{flaghi!55!white}\begin{tabular}[t]{@{}c@{}}0.231\\{\color{black!65}0\%}\end{tabular} \\[1pt]
MLP-flow & \nScoredMLPflow{} & \cellcolor{flaghi!51!white}\begin{tabular}[t]{@{}c@{}}0.820\\{\color{black!65}8\%}\end{tabular} & \cellcolor{flaghi!52!white}\begin{tabular}[t]{@{}c@{}}0.754\\{\color{black!65}5\%}\end{tabular} & \cellcolor{flaghi!55!white}\begin{tabular}[t]{@{}c@{}}0.963\\{\color{black!65}0\%}\end{tabular} & \cellcolor{flaghi!55!white}\begin{tabular}[t]{@{}c@{}}0.751\\{\color{black!65}0\%}\end{tabular} & \cellcolor{flaghi!55!white}\begin{tabular}[t]{@{}c@{}}0.691\\{\color{black!65}0\%}\end{tabular} & \cellcolor{flaghi!55!white}\begin{tabular}[t]{@{}c@{}}0.245\\{\color{black!65}0\%}\end{tabular} & \cellcolor{flaghi!55!white}\begin{tabular}[t]{@{}c@{}}0.486\\{\color{black!65}0\%}\end{tabular} \\[1pt]
RC & \nScoredRC{} & \cellcolor{flaghi!40!white}\begin{tabular}[t]{@{}c@{}}0.923\\{\color{black!65}27\%}\end{tabular} & \cellcolor{flaghi!39!white}\begin{tabular}[t]{@{}c@{}}0.887\\{\color{black!65}29\%}\end{tabular} & \cellcolor{flaghi!1!white}\begin{tabular}[t]{@{}c@{}}0.998\\{\color{black!65}99\%}$^\dagger$\end{tabular} & \cellcolor{flaghi!55!white}\begin{tabular}[t]{@{}c@{}}0.895\\{\color{black!65}0\%}\end{tabular} & \cellcolor{flaghi!55!white}\begin{tabular}[t]{@{}c@{}}0.856\\{\color{black!65}0\%}\end{tabular} & \cellcolor{flaghi!1!white}\begin{tabular}[t]{@{}c@{}}0.997\\{\color{black!65}97\%}\end{tabular} & \cellcolor{flaghi!1!white}\begin{tabular}[t]{@{}c@{}}0.996\\{\color{black!65}97\%}$^\dagger$\end{tabular} \\[1pt]
NODE & \nScoredNODE{} & \cellcolor{flaghi!8!white}\begin{tabular}[t]{@{}c@{}}0.976\\{\color{black!65}86\%}\end{tabular} & \cellcolor{flaghi!12!white}\begin{tabular}[t]{@{}c@{}}0.962\\{\color{black!65}78\%}\end{tabular} & \begin{tabular}[t]{@{}c@{}}0.998\\{\color{black!65}100\%}$^\dagger$\end{tabular} & \cellcolor{flaghi!38!white}\begin{tabular}[t]{@{}c@{}}0.952\\{\color{black!65}30\%}\end{tabular} & \cellcolor{flaghi!43!white}\begin{tabular}[t]{@{}c@{}}0.935\\{\color{black!65}21\%}\end{tabular} & \begin{tabular}[t]{@{}c@{}}0.996\\{\color{black!65}100\%}$^\dagger$\end{tabular} & \cellcolor{flaghi!10!white}\begin{tabular}[t]{@{}c@{}}0.996\\{\color{black!65}82\%}$^\dagger$\end{tabular} \\[1pt]
NSDE & \nScoredNSDE{} & \cellcolor{flaghi!6!white}\begin{tabular}[t]{@{}c@{}}0.980\\{\color{black!65}89\%}\end{tabular} & \cellcolor{flaghi!12!white}\begin{tabular}[t]{@{}c@{}}0.960\\{\color{black!65}78\%}\end{tabular} & \begin{tabular}[t]{@{}c@{}}0.998\\{\color{black!65}100\%}$^\dagger$\end{tabular} & \cellcolor{flaghi!39!white}\begin{tabular}[t]{@{}c@{}}0.953\\{\color{black!65}29\%}\end{tabular} & \cellcolor{flaghi!42!white}\begin{tabular}[t]{@{}c@{}}0.929\\{\color{black!65}24\%}\end{tabular} & \begin{tabular}[t]{@{}c@{}}0.996\\{\color{black!65}100\%}$^\dagger$\end{tabular} & \cellcolor{flaghi!5!white}\begin{tabular}[t]{@{}c@{}}0.996\\{\color{black!65}90\%}$^\dagger$\end{tabular} \\[1pt]
\midrule
Transformer (expl.) & \nScoredTransformer{} & \cellcolor{flaghi!9!white}\begin{tabular}[t]{@{}c@{}}0.981\\{\color{black!65}84\%}\end{tabular} & \cellcolor{flaghi!8!white}\begin{tabular}[t]{@{}c@{}}0.974\\{\color{black!65}86\%}\end{tabular} & \begin{tabular}[t]{@{}c@{}}0.997\\{\color{black!65}100\%}$^\dagger$\end{tabular} & \cellcolor{flaghi!39!white}\begin{tabular}[t]{@{}c@{}}0.949\\{\color{black!65}29\%}\end{tabular} & \cellcolor{flaghi!32!white}\begin{tabular}[t]{@{}c@{}}0.954\\{\color{black!65}42\%}\end{tabular} & \begin{tabular}[t]{@{}c@{}}0.996\\{\color{black!65}100\%}$^\dagger$\end{tabular} & \cellcolor{flaghi!8!white}\begin{tabular}[t]{@{}c@{}}0.996\\{\color{black!65}85\%}$^\dagger$\end{tabular} \\[1pt]
\midrule
SINDy-100 & \nScoredSindy{} & \begin{tabular}[t]{@{}c@{}}0.977\\{\color{black!65}100\%}\end{tabular} & \begin{tabular}[t]{@{}c@{}}0.969\\{\color{black!65}100\%}\end{tabular} & \begin{tabular}[t]{@{}c@{}}0.995\\{\color{black!65}100\%}$^\dagger$\end{tabular} & \cellcolor{flaghi!50!white}\begin{tabular}[t]{@{}c@{}}0.904\\{\color{black!65}10\%}\end{tabular} & \cellcolor{flaghi!55!white}\begin{tabular}[t]{@{}c@{}}0.748\\{\color{black!65}0\%}\end{tabular} & \begin{tabular}[t]{@{}c@{}}0.995\\{\color{black!65}100\%}$^\dagger$\end{tabular} & \begin{tabular}[t]{@{}c@{}}0.993\\{\color{black!65}100\%}$^\dagger$\end{tabular} \\[1pt]
$G_0$, NODE path & \nScoredGzeroH{} & \begin{tabular}[t]{@{}c@{}}0.988\\{\color{black!65}100\%}\end{tabular} & \cellcolor{flaghi!5!white}\begin{tabular}[t]{@{}c@{}}0.975\\{\color{black!65}90\%}\end{tabular} & \begin{tabular}[t]{@{}c@{}}0.994\\{\color{black!65}100\%}$^\dagger$\end{tabular} & \cellcolor{flaghi!17!white}\begin{tabular}[t]{@{}c@{}}0.978\\{\color{black!65}70\%}\end{tabular} & \cellcolor{flaghi!17!white}\begin{tabular}[t]{@{}c@{}}0.967\\{\color{black!65}70\%}\end{tabular} & \begin{tabular}[t]{@{}c@{}}0.993\\{\color{black!65}100\%}$^\dagger$\end{tabular} & \cellcolor{flaghi!50!white}\begin{tabular}[t]{@{}c@{}}0.988\\{\color{black!65}10\%}$^\dagger$\end{tabular} \\[1pt]
$G_0$, NSDE path & \nScoredGzeroHSig{} & \begin{tabular}[t]{@{}c@{}}0.990\\{\color{black!65}100\%}\end{tabular} & \begin{tabular}[t]{@{}c@{}}0.978\\{\color{black!65}100\%}\end{tabular} & \begin{tabular}[t]{@{}c@{}}0.993\\{\color{black!65}100\%}$^\dagger$\end{tabular} & \cellcolor{flaghi!22!white}\begin{tabular}[t]{@{}c@{}}0.982\\{\color{black!65}60\%}\end{tabular} & \cellcolor{flaghi!17!white}\begin{tabular}[t]{@{}c@{}}0.971\\{\color{black!65}70\%}\end{tabular} & \begin{tabular}[t]{@{}c@{}}0.993\\{\color{black!65}100\%}$^\dagger$\end{tabular} & \cellcolor{flaghi!38!white}\begin{tabular}[t]{@{}c@{}}0.744\\{\color{black!65}30\%}$^\dagger$\end{tabular} \\[1pt]
\bottomrule
\end{tabular}
\end{table*}

% GENERATED by scripts/camera_ready/make_fig1.py --appendix from results/camera_ready/camera_ready_metrics.json via numbers.tex -- do not edit.
% HAND-EDITED 2026-10-06 (shorter caption/layout for the extended version): port to the generator before the next `make paper`, or this edit is overwritten.
\begin{figure}[t]
\centering
\includegraphics[width=\textwidth]{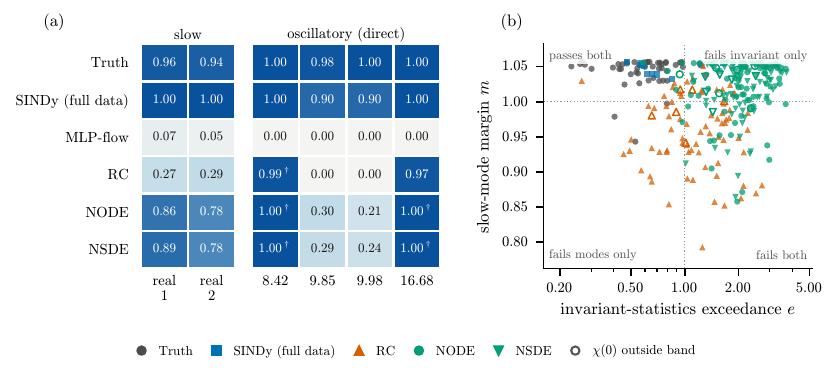}
\caption{\textbf{Per-mode passes and the dissociation map.} (a)~Fraction of rollouts passing each leading mode; $\dagger$ as in Table~\ref{tab:appx-modes}; MLP-drift passes no mode (not shown). (b)~One point per scored rollout: $e$ is the largest exceedance over the four invariant-statistics floors ($W_1$, ACF: value over floor; switch ratio, $f_+$: distance from the interval centre over its half-width) and $m$ the smaller of the two slow-mode overlaps over their thresholds. Rollouts with $e\le\paramOne{}$ lie inside every floor, those with $m\ge\paramOne{}$ pass both slow modes; hollow markers: $\chi(0)$ outside the band. MLP models not shown.}
\label{fig:appx-modes}
\end{figure}

Table~\ref{tab:appx-modes} and Fig.~\ref{fig:appx-modes}a resolve the mode rungs into the individual modes. RC, NODE and NSDE miss the individual eigenvectors of the near-degenerate oscillatory pair in most rollouts while recovering their joint subspace; the MLPs recover neither, and only SINDy resolves the split. The thresholds of the $\omega\approx\refFreqOne{}$ and $\refFreqSix{}$ pairs sit at the \paramThrCap{} cap (uncapped \thrUncOne{} and \thrUncSix{}), so passes on these modes can be decided by the cap, as for \pcRCOne{} of RC's \paRCOne{} passes at $\refFreqOne{}$ and \pcNODEOne{} and \pcNSDEOne{} of the \paNODEOne{} and \paNSDEOne{} NODE and NSDE passes there, against \pcSindyFullOne{} for SINDy. Rebuilding the two slow-mode thresholds from the \stsTruthN{} fresh truths (leave-one-out minimum: \stsThrOne{} and \stsThrTwo{} instead of \stsThrDseventeenOne{} and \stsThrDseventeenTwo{}) lowers RC's slow-mode failure rate from \pct{\RCslowflag} to \pct{\RCslowflagTruthThr} and NODE's from \pct{\NODEslowflag} to \pct{\NODEslowflagTruthThr}, and the RC exhibit becomes \stsRCnXslowfail{} of \stsRCnX{} rollouts against \stsTruenXslowfail{} of \stsTruenX{} truths (Fisher $p=\stsExhibitFisherp{}$; Table~\ref{tab:appx-slowthr}).

\section{Discussion}
The test is calibrated, not exhaustive: one low-dimensional system, one forcing direction and one observable. Thresholds come from \paramNCalib{} calibration truths; rebuilding the slow-mode thresholds from the \Truen{} fresh truths lowers RC's slow-mode failure rate from \pct{\RCslowflag} to \pct{\RCslowflagTruthThr}, and the 0.99 cap decides many oscillatory passes (Sec.~\ref{sec:modes}). Except for the NSDE, rollouts use the true noise intensity; as a robustness check, SINDy fitted to only \fitTHundred{} time units acquires symmetry-breaking terms (Sec.~\ref{app:sindy}).

Neither invariant statistics nor $\chi(0)$ certifies the slow dynamics: a quarter of RC's rollouts pass both with the wrong slow modes. A one-step drift regressed at a coarse step builds the $O(\Delta t)$ secant bias into the model and breaks the attractor. One-step flow maps can match the attractor (RC) yet misplace the slow, weakly damped modes, plausibly because these change little over one step and barely enter a one-step loss. Training on multi-step windows through the integrator (NODE, NSDE) recovers the slow modes in three quarters of rollouts but over-switches; the rollout-fine-tuned transformer shows the same split. Multi-step losses, now common in weather emulators~\cite{lam_learning_2023,bodnar_foundation_2025}, thus appear to carry the slow dynamics, and could be paired with invariant-measure constraints or spectral priors rather than with more flexible architectures. As slow modes govern transitions between the quasi-invariant sets of a metastable system and, for forcings that couple to them, its response, a mode-resolved check, needing one unperturbed trajectory per model, should accompany invariant statistics in validating learned simulators. 
A natural next step is to test whether these patterns persist in spatially extended systems and for state-of-the-art architectures of current emulators.

{\small
\bibliographystyle{plainnat}
\bibliography{references}
}

\appendix
\section{Supplementary tables}
\label{app:tables}
This appendix collects the settings and calibration tables behind Secs.~\ref{app:models} and~\ref{sec:exp}: model parameters and seeds (Table~\ref{tab:params}), the NODE sweep and the NSDE learned noise (Tables~\ref{tab:appx-sweep} and~\ref{tab:appx-noise}), the share of $\chi(0)$ carried by the leading resonances (Table~\ref{tab:appx-shares}), the null calibration of the floors and the $\chi(0)$ band (Tables~\ref{tab:appx-null} and~\ref{tab:appx-bands}) and the slow-mode threshold sensitivity (Table~\ref{tab:appx-slowthr}).
% GENERATED by scripts/camera_ready/make_appendix_tables.py from results/camera_ready/camera_ready_metrics.json via numbers.tex -- do not edit.
\begin{table}[t]
\centering
\scriptsize
\caption{\textbf{Parameters.} System, instrument, seeds (top) and per-model training and rollout (bottom). Every non-NSDE surrogate is rolled out with the true noise amplitude $s$ (oracle noise).}
\label{tab:params}
\begin{tabular}{@{}lp{0.78\linewidth}@{}}
\toprule
System & $\sigma,\rho,\beta = \paramSigma{}, \paramRho{}, \paramBeta{}$; noise $s=\paramS{}$; forcing on $x$ \\
Truth scheme & Euler--Maruyama, $dt_{\rm raw}=\paramDtRaw{}$; burn-in \paramBurnIn{} t.u.; $T=\paramTen{}^{\paramTExp{}}$ t.u. \\
Instrument & $\tau=\paramTau{}$, $K_{\max}=\paramKmax{}$ (\paramDict{} functions), TSVD $\paramTen{}^{\paramRelThreshExp{}}$, $n_{\rm lead}=\paramNLead{}$ \\
Mode thresholds & clip(mean $-\,\paramCalibK{}$ std, \paramThrFloor{}, \paramThrCap{}) over \paramNCalib{} calibration truths; grid \paramGridN{} points \\
Seeds & ref \seedRef{}; grid \seedGrid{}; calibration \seedCalibOne{}--\seedCalibTwo{}; floor truths \seedFloorOne{}--\seedFloorTwo{}; runner truths \seedTruthCert{}, \seedTruthExtOne{}--\seedTruthExtTwo{} \\
 & E2 train \seedEtwoTrainOne{}--\seedEtwoTrainTwo{} $\times$ rollout \seedEtwoRolloutOne{}--\seedEtwoRolloutTwo{}; NODE init \seedNodeInitOne{}--\seedNodeInitTwo{}; NSDE init \seedNsdeInitOne{}--\seedNsdeInitTwo{}; G0 \seedGzeroHOne{}, \seedGzeroHTwo{}--\seedGzeroHThree{} and \seedGzeroSigOne{}, \seedGzeroSigTwo{}--\seedGzeroSigThree{}; Transformer init \seedTrfInitOne{}--\seedTrfInitTwo{} \\
\bottomrule
\end{tabular}
\par\medskip
\begin{tabular}{@{}lcccl@{}}
\toprule
Model & $\Delta t_{\rm fit}$ & rollout step & train span (t.u.) & rollout noise \\
\midrule
SINDy (full data) & \paramDtFitSindyFull{} & \paramDtRollSindyFull{} & \paramTrainTUSindyFull{} & \texttt{x += dt\,F(x) + s\,sqrt(dt)\,$\xi$} \\
SINDy-100 & \paramDtFitSindy{} & \paramDtRollSindy{} & \paramTrainTUSindy{} & \texttt{x += dt\,F(x) + s\,sqrt(dt)\,$\xi$} \\
MLP-drift & \paramDtFitMLPdrift{} & \paramDtRollMLPdrift{} & \paramTrainTUMLPdrift{} & \texttt{x += dt\,F\_theta(x) + s\,sqrt(dt)\,$\xi$} \\
MLP-flow & \paramDtFitMLPflow{} & \paramDtRollMLPflow{} & \paramTrainTUMLPflow{} & \texttt{x = MLP(x) + s\,sqrt(dt)\,$\xi$} \\
RC & \paramDtFitRC{} & \paramDtRollRC{} & \paramTrainTURC{} & \texttt{x = ESN(x, r) + s\,sqrt(dt)\,$\xi$} \\
NODE & \paramDtFitNODE{} & \paramDtRollNODE{} & \paramTrainTUNODE{} & RK4 step of $F_\theta$ + $s\sqrt{h}\,\xi$ (oracle $s$) \\
NSDE & \paramDtFitNSDE{} & \paramDtRollNSDE{} & \paramTrainTUNSDE{} & RK4 step of $F_\theta$ + $\sigma_\phi(x)\sqrt{h}\,\xi$ (learned) \\
\midrule
Transformer (expl.) & \paramDtFitTransformer{} & \paramDtRollTransformer{} & \paramTrainTUTransformer{} & $x_n+\Delta_\theta(x_{n-K+\paramOne{}},\dots,x_n)$ + $s\sqrt{h}\,\xi$ (oracle $s$) \\
\bottomrule
\end{tabular}
\end{table}

% GENERATED by scripts/camera_ready/make_appendix_tables.py from results/camera_ready/camera_ready_metrics.json via numbers.tex -- do not edit.
\begin{table}[t]
\centering
\scriptsize
\caption{\textbf{NODE sweep} (one seed): validation MSE ($\times \paramTen{}^{\paramMseExp{}}$, standardized) at horizons $L=\paramHorizonA{},\paramHorizonB{},\paramHorizonC{}$; selection by the lowest MSE at the common horizon $L=\paramHorizonC{}$ (first row). Drift RMSE on the grid is reported only.}
\label{tab:appx-sweep}
\setlength{\tabcolsep}{3pt}
\begin{tabular}{@{}cccccccc@{}}
\toprule
$L$ & lr & MSE$_{\paramHorizonA{}}$ & MSE$_{\paramHorizonB{}}$ & MSE$_{\paramHorizonC{}}$ & epochs & train (s) & drift RMSE $x/y/z$ \\
\midrule
\sweepLOne{} & \sweepLrOne{} & \sweepMseLoneOne{} & \sweepMseLfiveOne{} & \sweepMseLtwentyOne{} & \sweepEpochsOne{} & \sweepTrainSOne{} & \sweepRmseXOne{}/\sweepRmseYOne{}/\sweepRmseZOne{} \\
\sweepLTwo{} & \sweepLrTwo{} & \sweepMseLoneTwo{} & \sweepMseLfiveTwo{} & \sweepMseLtwentyTwo{} & \sweepEpochsTwo{} & \sweepTrainSTwo{} & \sweepRmseXTwo{}/\sweepRmseYTwo{}/\sweepRmseZTwo{} \\
\sweepLThree{} & \sweepLrThree{} & \sweepMseLoneThree{} & \sweepMseLfiveThree{} & \sweepMseLtwentyThree{} & \sweepEpochsThree{} & \sweepTrainSThree{} & \sweepRmseXThree{}/\sweepRmseYThree{}/\sweepRmseZThree{} \\
\sweepLFour{} & \sweepLrFour{} & \sweepMseLoneFour{} & \sweepMseLfiveFour{} & \sweepMseLtwentyFour{} & \sweepEpochsFour{} & \sweepTrainSFour{} & \sweepRmseXFour{}/\sweepRmseYFour{}/\sweepRmseZFour{} \\
\bottomrule
\end{tabular}
\end{table}
\begin{table}[t]
\centering
\scriptsize
\caption{\textbf{NSDE learned noise} relative to the true $s$, at $h=\noiseH{}$; $a=(\sigma,\paramOne{},\beta)$.}
\label{tab:appx-noise}
\setlength{\tabcolsep}{3pt}
\begin{tabular}{@{}lccc@{}}
\toprule
& $x$ & $y$ & $z$ \\
\midrule
NSDE mean $\sigma_\phi/s$ (median of \noiseN{}) & \noisePhiX{} & \noisePhiY{} & \noisePhiZ{} \\
NSDE constant $\hat\sigma/s$ & \noiseHatX{} & \noiseHatY{} & \noiseHatZ{} \\
true drift, step $h$, constant $\hat\sigma/s$ & \noiseTrueDriftX{} & \noiseTrueDriftY{} & \noiseTrueDriftZ{} \\
$\sqrt{(\paramOne{}-e^{-\paramTwo{}ah})/(\paramTwo{}ah)}$ & \noiseRelaxX{} & \noiseRelaxY{} & \noiseRelaxZ{} \\
\bottomrule
\end{tabular}
\end{table}

% GENERATED by scripts/camera_ready/make_appendix_tables.py from results/camera_ready/camera_ready_metrics.json via numbers.tex -- do not edit.
\begin{table}[t]
\centering
\scriptsize
\caption{\textbf{Where $\chi(0)$ sits.} Signed share of $\chi(0)$ carried by the leading resonances (median [min, max]), their share of the summed absolute mode contributions, and the fraction of $\chi(0)$ accumulated by $t=\paramOne{}$.}
\label{tab:appx-shares}
\setlength{\tabcolsep}{3pt}
\begin{tabular}{@{}lcccc@{}}
\toprule
& $n$ & leading share & leading $|$share$|$ & accumulated by $t$\,=\,\paramOne{} \\
\midrule
reference & \shareNRef{} & \shareLeadRefmed{} & -- & \accOneRefmed{} \\
truth & \shareNTrue{} & \shareLeadTruemed{}\,[\shareLeadTruemin{},\,\shareLeadTruemax{}] & \shareAbsLeadTruemed{} & \accOneTruemed{}\,[\accOneTruemin{},\,\accOneTruemax{}] \\
SINDy-100 & \shareNSindy{} & \shareLeadSindymed{}\,[\shareLeadSindymin{},\,\shareLeadSindymax{}] & \shareAbsLeadSindymed{} & \accOneSindymed{}\,[\accOneSindymin{},\,\accOneSindymax{}] \\
RC & \shareNRC{} & \shareLeadRCmed{}\,[\shareLeadRCmin{},\,\shareLeadRCmax{}] & \shareAbsLeadRCmed{} & \accOneRCmed{}\,[\accOneRCmin{},\,\accOneRCmax{}] \\
NODE & \shareNNODE{} & \shareLeadNODEmed{}\,[\shareLeadNODEmin{},\,\shareLeadNODEmax{}] & \shareAbsLeadNODEmed{} & \accOneNODEmed{}\,[\accOneNODEmin{},\,\accOneNODEmax{}] \\
NSDE & \shareNNSDE{} & \shareLeadNSDEmed{}\,[\shareLeadNSDEmin{},\,\shareLeadNSDEmax{}] & \shareAbsLeadNSDEmed{} & \accOneNSDEmed{}\,[\accOneNSDEmin{},\,\accOneNSDEmax{}] \\
\bottomrule
\end{tabular}
\end{table}

% GENERATED by scripts/camera_ready/make_appendix_tables.py from results/camera_ready/camera_ready_metrics.json via numbers.tex -- do not edit.
% HAND-EDITED 2026-10-06 (shorter caption/layout for the extended version): port to the generator before the next `make paper`, or this edit is overwritten.
\begin{table*}[t]
\centering
\scriptsize
\caption{\textbf{Null calibration.} Leave-one-out out-rates on the independent true ensemble (floors, E1b band) and the pooled band; out-of-range counts of the \bandNullnTest{} runner truths, with the one-sided binomial test against the leave-one-out rate. Only the E1b band fails ($p<\paramAlpha{}$), so the band is the pooled range over \bandPooledn{} truths. Runner vs.\ E1b $\chi(0)$ variances: $F=\bandVarF{}$, $p=\bandVarp{}$.}
\label{tab:appx-null}
\setlength{\tabcolsep}{3pt}
\begin{tabular}{@{}lccccccc@{}}
\toprule
& $W_1$ & ACF & switch & $f_+$ & all four & band (E1b) & band (pooled) \\
\midrule
leave-one-out out-rate & \floorWOnelooOut{} & \floorAcflooOut{} & \floorSwitchlooOut{} & \floorFPluslooOut{} & \floorJointLooOut{} & \bandEoneblooOut{} & \bandPooledlooOut{} \\
runner truths outside & \floorWOnekOut{} & \floorAcfkOut{} & \floorSwitchkOut{} & \floorFPluskOut{} & -- & \bandNullkOut{} & -- \\
one-sided binomial $p$ & \floorWOnep{} & \floorAcfp{} & \floorSwitchp{} & \floorFPlusp{} & -- & \bandNullp{} & -- \\
\bottomrule
\end{tabular}
\end{table*}
\begin{table*}[t]
\centering
\scriptsize
\caption{\textbf{$\chi(0)$ in-band counts under each candidate band}, and diverged rollouts ($|x|,|y|,|z|>\gateBoundX{},\gateBoundY{},\gateBoundZ{}$ or non-finite), excluded from all rates.}
\label{tab:appx-bands}
\setlength{\tabcolsep}{3pt}
\begin{tabular}{@{}lcccccc@{}}
\toprule
Model & $n$ & pooled $[\bandPooledlo{},\bandPooledhi{}]$ & E1b $[\bandEoneblo{},\bandEonebhi{}]$ & $N$\,=\,\bandNEighteenSn{} $[\bandNEighteenSlo{},\bandNEighteenShi{}]$ & runner $[\bandRunnerlo{},\bandRunnerhi{}]$ & diverged \\
\midrule
Truth (perfect surrogate) & \nScoredTrue{} & \inBandPooledTrue{} & \inBandEonebTrue{} & \inBandNEighteenSTrue{} & \inBandRunnerTrue{} & \nDivTrue{} \\
SINDy (full data) & \nScoredSindyFull{} & \inBandPooledSindyFull{} & \inBandEonebSindyFull{} & \inBandNEighteenSSindyFull{} & \inBandRunnerSindyFull{} & \nDivSindyFull{} \\
MLP-drift & \nScoredMLPdrift{} & \inBandPooledMLPdrift{} & \inBandEonebMLPdrift{} & \inBandNEighteenSMLPdrift{} & \inBandRunnerMLPdrift{} & \nDivMLPdrift{} \\
MLP-flow & \nScoredMLPflow{} & \inBandPooledMLPflow{} & \inBandEonebMLPflow{} & \inBandNEighteenSMLPflow{} & \inBandRunnerMLPflow{} & \nDivMLPflow{} \\
RC & \nScoredRC{} & \inBandPooledRC{} & \inBandEonebRC{} & \inBandNEighteenSRC{} & \inBandRunnerRC{} & \nDivRC{} \\
NODE & \nScoredNODE{} & \inBandPooledNODE{} & \inBandEonebNODE{} & \inBandNEighteenSNODE{} & \inBandRunnerNODE{} & \nDivNODE{} \\
NSDE & \nScoredNSDE{} & \inBandPooledNSDE{} & \inBandEonebNSDE{} & \inBandNEighteenSNSDE{} & \inBandRunnerNSDE{} & \nDivNSDE{} \\
SINDy-100 & \nScoredSindy{} & \inBandPooledSindy{} & \inBandEonebSindy{} & \inBandNEighteenSSindy{} & \inBandRunnerSindy{} & \nDivSindy{} \\
G0 (NODE path) & \nScoredGzeroH{} & \inBandPooledGzeroH{} & \inBandEonebGzeroH{} & \inBandNEighteenSGzeroH{} & \inBandRunnerGzeroH{} & \nDivGzeroH{} \\
G0 (NSDE path) & \nScoredGzeroHSig{} & \inBandPooledGzeroHSig{} & \inBandEonebGzeroHSig{} & \inBandNEighteenSGzeroHSig{} & \inBandRunnerGzeroHSig{} & \nDivGzeroHSig{} \\
\bottomrule
\end{tabular}
\end{table*}

% GENERATED by scripts/camera_ready/make_appendix_tables.py from results/camera_ready/camera_ready_metrics.json via numbers.tex -- do not edit.
% HAND-EDITED 2026-10-06 (shorter caption/layout for the extended version): port to the generator before the next `make paper`, or this edit is overwritten.
\begin{table}[t]
\centering
\scriptsize
\caption{\textbf{Slow-mode threshold sensitivity.} Thresholds rebuilt from the \stsTruthN{} fresh truths by a leave-one-out minimum ($\lambda\approx\refReTwo{}$: \stsThrOne{} instead of \stsThrDseventeenOne{}; $\lambda\approx\refReThree{}$: \stsThrTwo{} instead of \stsThrDseventeenTwo{}): each truth is judged against the other truths, every model against all of them. Stored overlaps are re-thresholded, nothing is re-run; Table~\ref{tab:flags} keeps the calibration-set thresholds.}
\label{tab:appx-slowthr}
\setlength{\tabcolsep}{3pt}
\begin{tabular}{@{}lccccc@{}}
\toprule
& & \multicolumn{2}{c}{slow modes failed (\%)} & \multicolumn{2}{c}{pass, rebuilt thresholds (\%)} \\
\cmidrule(lr){3-4}\cmidrule(lr){5-6}
& $n$ & calibration set & rebuilt & $\lambda\approx\refReTwo{}$ & $\lambda\approx\refReThree{}$ \\
\midrule
Truth (perfect surrogate) & \stsTruen{} & \pctn{\stsTrueflagDseventeen} & \pctn{\stsTrueflagThr} & \pctn{\stsTruepassOne} & \pctn{\stsTruepassTwo} \\
SINDy (full data) & \stsSindyFulln{} & \pctn{\stsSindyFullflagDseventeen} & \pctn{\stsSindyFullflagThr} & \pctn{\stsSindyFullpassOne} & \pctn{\stsSindyFullpassTwo} \\
MLP-drift & \stsMLPdriftn{} & \pctn{\stsMLPdriftflagDseventeen} & \pctn{\stsMLPdriftflagThr} & \pctn{\stsMLPdriftpassOne} & \pctn{\stsMLPdriftpassTwo} \\
MLP-flow & \stsMLPflown{} & \pctn{\stsMLPflowflagDseventeen} & \pctn{\stsMLPflowflagThr} & \pctn{\stsMLPflowpassOne} & \pctn{\stsMLPflowpassTwo} \\
RC & \stsRCn{} & \pctn{\stsRCflagDseventeen} & \pctn{\stsRCflagThr} & \pctn{\stsRCpassOne} & \pctn{\stsRCpassTwo} \\
NODE & \stsNODEn{} & \pctn{\stsNODEflagDseventeen} & \pctn{\stsNODEflagThr} & \pctn{\stsNODEpassOne} & \pctn{\stsNODEpassTwo} \\
NSDE & \stsNSDEn{} & \pctn{\stsNSDEflagDseventeen} & \pctn{\stsNSDEflagThr} & \pctn{\stsNSDEpassOne} & \pctn{\stsNSDEpassTwo} \\
SINDy-100 & \stsSindyn{} & \pctn{\stsSindyflagDseventeen} & \pctn{\stsSindyflagThr} & \pctn{\stsSindypassOne} & \pctn{\stsSindypassTwo} \\
G0 (NODE path) & \stsGzeroHn{} & \pctn{\stsGzeroHflagDseventeen} & \pctn{\stsGzeroHflagThr} & \pctn{\stsGzeroHpassOne} & \pctn{\stsGzeroHpassTwo} \\
G0 (NSDE path) & \stsGzeroHSign{} & \pctn{\stsGzeroHSigflagDseventeen} & \pctn{\stsGzeroHSigflagThr} & \pctn{\stsGzeroHSigpassOne} & \pctn{\stsGzeroHSigpassTwo} \\
\midrule
Transformer (expl.) & \stsTransformern{} & \pctn{\stsTransformerflagDseventeen} & \pctn{\stsTransformerflagThr} & \pctn{\stsTransformerpassOne} & \pctn{\stsTransformerpassTwo} \\
\bottomrule
\end{tabular}
\end{table}

\end{document}